\documentclass[12pt,reqno, table]{amsart}
\usepackage{amsfonts, amsmath, amssymb, color}
\usepackage{amsthm}
\usepackage{mathrsfs}
\usepackage{booktabs}
\usepackage{graphicx}
\usepackage{enumitem}
\usepackage{bm}
\usepackage{float}
\usepackage{hhline}
\usepackage{soul, color, xcolor}
\usepackage{lipsum}
\usepackage{lscape}
\usepackage{subfig}
\usepackage{textcomp}
\usepackage{tikz}
\usetikzlibrary{decorations.markings}
\usetikzlibrary{decorations.pathmorphing,shapes,arrows,positioning}
\usetikzlibrary{arrows.meta}
\usepackage{xcolor}
\usepackage{young}
\usepackage{youngtab}
\usepackage{ytableau}
\usepackage{rotating}
\usepackage[colorlinks=true, linkcolor=blue, citecolor=magenta, urlcolor=blue, backref=page]{hyperref}
\usepackage{scalefnt}
\usetikzlibrary{calc}
\usepackage{a4wide}
\usepackage{pdfpages}

\usepackage{young}
\usepackage{youngtab}
\usepackage{ytableau}
\usepackage{rotating} 
\usepackage{soul}
\usepackage{tabularx} 
\usepackage{array}    
\usepackage{xcolor} 

\usepackage{listings}
\usepackage{upquote}

\definecolor{codebg}{HTML}{F7F7F8}
\definecolor{codeframe}{HTML}{D9D9E3}
\definecolor{codekw}{HTML}{1F5AA6}
\definecolor{codestr}{HTML}{A31515}
\definecolor{codecom}{HTML}{6A737D}
\definecolor{codenum}{HTML}{888888}
\definecolor{codeemph}{HTML}{7A3E9D}

\lstdefinestyle{pythonstyle}{
  language=Python,
  backgroundcolor=\color{codebg},
  basicstyle=\ttfamily\footnotesize,
  keywordstyle=\color{codekw}\bfseries,
  commentstyle=\color{codecom}\itshape,
  stringstyle=\color{codestr},
  numberstyle=\scriptsize\color{codenum},
  numbers=left,
  numbersep=10pt,
  stepnumber=1,
  showstringspaces=false,
  showspaces=false,
  showtabs=false,
  tabsize=4,
  breaklines=true,
  breakatwhitespace=false,
  keepspaces=true,
  columns=fullflexible,
  upquote=true,
  frame=single,
  framerule=0.4pt,
  rulecolor=\color{codeframe},
  framesep=6pt,
  xleftmargin=14pt,
  framexleftmargin=10pt,
  aboveskip=1.2ex,
  belowskip=1.2ex,
  captionpos=b
}

\hypersetup{colorlinks,linkcolor=blue,urlcolor=cyan,citecolor=magenta}

\allowdisplaybreaks[4]
\newcommand{\Yd}[2]{%
  \begin{tikzpicture}[scale=.35,baseline=(current bounding box.center)]
    \begin{scope}
      \clip #1;
      \draw [color=black!25] (0,0) grid #2;
    \end{scope}
    \draw[thick] #1;
  \end{tikzpicture}%
}

\newcommand{\YDthreeone}{\Yd{(0,0) -| (1,1) -| (3,2) -- (0,2) -- cycle}{(3,2)}}

\newcommand{\YDtwoonone}{\Yd{(0,0) -| (1,2) -| (2,3) -- (0,3) -- cycle}{(2,3)}}

\newcommand{\YDtwtwo}{\Yd{(0,0) rectangle (2,2)}{(2,2)}}

\newcommand{\YDtwo}{\Yd{(0,0) rectangle (2,1)}{(2,1)}}

\newcommand{\YDonesone}{\Yd{(0,0) rectangle (1,2)}{(1,2)}}

\newcommand{\YDfourcol}{\Yd{(0,0) rectangle (1,4)}{(1,4)}}

\newcommand{\YDtweone}{\Yd{(0,0) -| (1,1) -| (2,2) -- (0,2) -- cycle}{(2,2)}}

\newcommand{\Base}[2]{%
  \begin{scope}
    \clip #1;
    \draw [color=black!25] (0,0) grid #2;
  \end{scope}
  \draw[thick] #1;
}

\newcommand{\YDTop}{
\begin{tikzpicture}[scale = 0.4]
\Base{(0,0) -| (1,1) -| (2,2) -| (3,4) -| (6,5) -| (7,8) -- (0,8) -- cycle}{(7,8)}
\node at (-.7,7.5) {1}; \node at (-.7,6.5) {3}; \node at (-.7,5.5) {2}; \node at (-.7,4.5) {3};
\node at (-.7,3.5) {2}; \node at (-.7,2.5) {1}; \node at (-.7,1.5) {2}; \node at (-.7,0.5) {3};
\end{tikzpicture}
}

\newcommand{\YDMidL}{%
\begin{tikzpicture}[scale = 0.4]
\Base{(0,0) -| (1,3) -| (3,4) -| (6,8) -- (0,8) -- cycle}{(6,8)}
\node at (-.7,7.5) {3}; \node at (-.7,6.5) {2}; \node at (-.7,5.5) {1}; \node at (-.7,4.5) {3};
\node at (-.7,3.5) {2}; \node at (-.7,2.5) {2}; \node at (-.7,1.5) {1}; \node at (-.7,0.5) {3};
\end{tikzpicture}%
}

\newcommand{\YDMidR}{%
\begin{tikzpicture}[scale = 0.4]
\Base{(0,0) -| (1,1) -| (2,2) -| (3,4) -| (4,5) -| (5,6) -| (6,8) -- (0,8) -- cycle}{(6,8)}
\node at (-.7,7.5) {3}; \node at (-.7,6.5) {2}; \node at (-.7,5.5) {3}; \node at (-.7,4.5) {1};
\node at (-.7,3.5) {2}; \node at (-.7,2.5) {1}; \node at (-.7,1.5) {2}; \node at (-.7,0.5) {3};
\end{tikzpicture}%
}

\newcommand{\YDBotI}{%
\begin{tikzpicture}[scale = 0.4]
\Base{(0,0) -| (1,3) -| (3,4) -| (5,7) -| (6,8) -- (0,8) -- cycle}{(6,8)}
\node at (-.7,7.5) {3}; \node at (-.7,6.5) {1}; \node at (-.7,5.5) {3}; \node at (-.7,4.5) {2};
\node at (-.7,3.5) {2}; \node at (-.7,2.5) {2}; \node at (-.7,1.5) {1}; \node at (-.7,0.5) {3};
\end{tikzpicture}%
}

\newcommand{\YDBotII}{%
\begin{tikzpicture}[scale = 0.4]
\Base{(0,0) -| (3,1) -| (6,5) -- (0,5) -- cycle}{(6,5)}
\node at (-.7,4.5) {3}; \node at (-.7,3.5) {2}; \node at (-.7,2.5) {1}; \node at (-.7,1.5) {3}; \node at (-.7,0.5) {2};
\end{tikzpicture}%
}

\newcommand{\YDBotIII}{%
\begin{tikzpicture}[scale = 0.4]
\Base{(0,0) -| (1,1) -| (2,2) -| (3,4) -| (4,7) -| (6,8) -- (0,8) -- cycle}{(6,8)}
\node at (-.7,7.5) {3}; \node at (-.7,6.5) {3}; \node at (-.7,5.5) {2}; \node at (-.7,4.5) {1};
\node at (-.7,3.5) {2}; \node at (-.7,2.5) {1}; \node at (-.7,1.5) {2}; \node at (-.7,0.5) {3};
\end{tikzpicture}%
}

\newcommand{\YDBotIV}{%
\begin{tikzpicture}[scale = 0.4]
\Base{(0,0) -| (3,3) -| (4,4) -| (5,5) -| (6,7) -- (0,7) -- cycle}{(6,7)}
\node at (-.7,6.5) {3}; \node at (-.7,5.5) {2}; \node at (-.7,4.5) {3}; \node at (-.7,3.5) {1};
\node at (-.7,2.5) {2}; \node at (-.7,1.5) {1}; \node at (-.7,0.5) {1};
\end{tikzpicture}%
}

\newtheorem{mainthm}{Theorem}           

\newtheorem{thm}{Theorem}[section]
\newtheorem{lem}[thm]{Lemma}
\newtheorem{prop}[thm]{Proposition}
\newtheorem{cor}[thm]{Corollary}
\newtheorem{rem}[thm]{Remark}
\theoremstyle{plain}
\newtheorem{defn}[thm]{Definition}
\newtheorem{exmp}[thm]{Example}

\newcommand{\Rmnum}[1]{\expandafter\@slowromancap\romannumeral #1@}
\makeatother
\newcommand{\la}{\lambda}
\newcommand{\brho}{\bm{\rho}}
\newcommand{\btau}{\bm{\tau}}
\newcommand{\balp}{\bm{\alpha}}
\newcommand{\bnu}{\bm{\nu}}
\newcommand{\bbeta}{\bm{\beta}}
\newcommand{\bgam}{\bm{\gamma}}
\newcommand{\bdel}{\bm{\delta}}
\newcommand{\bvar}{\bm{\varepsilon}}
\newcommand{\bmu}{\bm{\mu}}
\newcommand{\bla}{\bm{\lambda}}

\newcommand{\mbtr}{\mathrm{mbtr}}
\newcommand{\wt}{\mathrm{wt}}
\newcommand{\T}{\mathrm{T}}
\newcommand{\m}{\mathtt{m}}
\newcommand{\Tr}{\mathrm{Tr}}

\numberwithin{equation}{section} \errorcontextlines=0

\numberwithin{equation}{section} \errorcontextlines=0

\def\D{\mathcal{D}}

\def\R{\mathcal{R}}

\def\I{\mathcal{I}}

\numberwithin{equation}{section}

\begin{document}


\title{Murnaghan--Nakayama rule for the cyclotomic Hecke algebra and applications} 
\author{Naihuan Jing}
\address{Department of Mathematics, North Carolina State University, Raleigh, NC 27695, USA}
\email{jing@ncsu.edu}
\author{Ning Liu}
\address{ Beijing International Center for Mathematical Research, Peking University, Beijing 100871, China}
\email{mathliu123@outlook.com}





\begin{abstract}
We establish a Murnaghan--Nakayama rule for the irreducible characters of the cyclotomic Hecke algebra $\mathscr H_{m,n}(q,u)$ on Shoji's standard elements. Combined with Shoji's determinacy result, our formula provides a direct combinatorial route to the full irreducible character table of $\mathscr H_{m,n}(q,u)$. Our construction is based on our recent multi-parameter Murnaghan--Nakayama rule for Macdonald polynomials and specializes uniformly to several previously known formulas, including those for the complex reflection group of type $G(m,1,n)$ and the Iwahori--Hecke algebras of types $A$ and $B$. In a dual framework, using the vertex operator realization of Schur functions, we also derive a complementary iterative formula for irreducible characters on upper multipartitions, which may be viewed as a dual Murnaghan--Nakayama rule.

As applications, we obtain a Regev-type formula and a L\"ubeck--Prasad--Adin--Roichman-type formula for cyclotomic Hecke algebras, extending the corresponding formulas for the Iwahori--Hecke algebra of type $A$ and the complex reflection group, respectively. We further introduce the notion of multiple bitrace for cyclotomic Hecke algebras and give a general combinatorial formula for the multiple bitrace. As a specialization, this yields the second orthogonality relation for irreducible characters of the complex reflection group $W_{m,n}$. For practical computation, we also include in an appendix a SageMath implementation of our Murnaghan--Nakayama rule, which computes individual character values and the full character table.
\end{abstract}
\maketitle
\tableofcontents




\section{Introduction}
\subsection{Background}\label{ss:Backgroud}
The Murnaghan--Nakayama rule \cite{Mur,Nak1,Nak2} serves as a combinatorial algorithm for calculating the irreducible characters of symmetric groups $\mathfrak{S}_n$. The Iwahori--Hecke algebra $\mathscr{H}_{n}(q)$ in type $A$, which is the $q$-deformation of $\mathfrak{S}_n$, is a unital associative algebra over the field of rational functions $\mathbb{C}(q)$. Ram \cite{Ram} used the quantum Schur--Weyl duality \cite{Jim} to establish a Frobenius-type character formula for $\mathscr{H}_{n}(q)$, linking it to one-row Hall-Littlewood functions and Schur symmetric functions (see also \cite{KV, KW}). He also introduced the $q$-Murnaghan--Nakayama rule, an iterative combinatorial method for computing irreducible characters using smaller partitions. An alternative approach was explored in \cite{Pfe1}. Recent work \cite{JL1} employed the vertex operator realization of Schur functions to revisit the $q$-Murnaghan--Nakayama rule.
This approach has led to the
formulation of a Murnaghan--Nakayama rule for the Hecke--Clifford algebra \cite{JL2}, utilizing the vertex operator realization of Schur $Q$-functions. 

At the end of \cite{Pfe2}, Pfeiffer noted the potential for deriving irreducible character formulas for the cyclotomic Hecke algebra $\mathscr{H}_{m,n}(q,\bm{u})$ over some specific elements, yet evidence was lacking to show that these character values could yield a complete character table. This raises the question of which elements can determine the characters and how to derive the character values for others. In 1998, inspired by Greene's derivation of the Murnaghan--Nakayama formula via Young’s seminormal representations \cite{Gre}, Halverson and Ram \cite{HR} established an analogue of the Murnaghan--Nakayama rule for $\mathscr{H}_{m,n}(q,\bm{u})$, focusing on certain ``standard elements". They claimed, without proof, that these values would completely determine the characters. Subsequently, a result by Mak \cite{Mak} confirms that the characters of $\mathscr{H}_{m,n}(q,\bm{u})$ can be determined by the values of the ``standard elements" through a highly complicated process.

In 2000, Shoji \cite{Sho} utilized Schur--Weyl reciprocity between $\mathscr{H}_{m,n}(q,\bm{u})$ and a Levi subalgebra of the quantized enveloping algebra $\mathscr{U}_{q}(gl_m)$ \cite{SS} to provide a Frobenius-type formula for the characters of $\mathscr{H}_{m,n}(q,\bm{u})$. He introduced a new presentation and standard elements for $\mathscr{H}_{m,n}(q,\bm{u})$, proving that the characters are completely determined by their values on these elements \cite[Proposition 7.5]{Sho}, addressing Pfeiffer's earlier question. In contrast to the ``standard elements'' mentioned above, Shoji's version of ``standard elements'' has simpler forms.

In 2013, Regev \cite{Reg} obtained an 
elegant formula, now usually referred to as the Regev formula, for the characters of super representations of the symmetric group. This result was derived through the super analogue of the classical Schur--Weyl duality (sometimes called the Schur-Sergeev duality), combined with the combinatorial theory of Lie superalgebras, which was initially developed by Sergeev \cite{Ser84} and later expanded in detail by Berele and Regev \cite{BR}. Taylor \cite{Tay} provided a new proof of Regev's work by utilizing the Murnaghan--Nakayama formula for the skew characters of symmetric groups. Building on Taylor’s work, Zhao \cite{Zhao3} extended the (quantum) Regev formula to the Iwahori--Hecke algebra of type $A$ by employing the corresponding Murnaghan--Nakayama rule for skew characters. An alternative proof of the quantum Regev formula was recently presented in \cite{JL1} using the vertex operator realization of Schur functions. 

In \cite{LP}, L\"ubeck and Prasad established a character identity that connects the irreducible character values of the Weyl group $\mathscr{W}_{2,n}$ to those of the symmetric group $\mathfrak{S}_{2n}$ using the Frobenius character formula for symmetric groups. Subsequently, Adin and Roichman \cite{AR} provided a combinatorial proof of this identity and extended it to the setting of wreath products. We collectively refer to these results as the L\"ubeck--Prasad--Adin--Roichman formula.

The second orthogonality relation for the irreducible characters of a finite group plays a crucial role in computing these characters. The second orthogonality relation for the irreducible characters of a split semisimple algebra also makes sense. An explicit form of the second orthogonality relation for the irreducible characters of the Brauer algebra was given in \cite{Ram97}. The $q$-analogue of the second orthogonality relation for the Hecke algebra of type $A$ formulated in terms of the bitrace was determined in \cite{HLR, JL1}. 
A spin analogue of the bitrace was recently developed for the Hecke--Clifford algebra in \cite{JL2}.

\subsection{Aims}

While the earlier Murnaghan--Nakayama-type formulas of Halverson--Ram are fundamental, they are not fully suited to our present purposes. In particular, Halverson and Ram computed characters only on certain ``standard elements'' and did not prove in that paper that these values determine the full character table (although they indicated a method for doing so later). Moreover, their approach relies on explicit seminormal/Hoefsmit-type representations together with double-centralizer (Clifford-theoretic) arguments and a bitrace reduction. This framework is structurally powerful, but less convenient for the Frobenius-formula-based combinatorial recursions and applications that we pursue here. Mak later resolved the determinacy issue from a structural viewpoint via a reducibility theorem for trace functions on Ariki--Koike algebras, but this is of a different nature and does not directly provide an explicit iterative combinatorial character formula. By contrast, Shoji's Frobenius-type formula and his choice of standard elements offer a natural starting point for such a development.

This paper has two main goals. First, we formulate a Murnaghan--Nakayama rule for the cyclotomic Hecke algebra based on Shoji's Frobenius-type formula, together with an effective operational procedure for computing character values. Compared with earlier formulations, our rule has a more canonical form and specializes uniformly to several known formulas, including those for the complex reflection group of type $G(m,1,n)$ and the Iwahori--Hecke algebras of types $A$ and $B$. In addition, combined with Shoji's result \cite[Proposition 7.5]{Sho}, it provides a simpler route to the complete character table of $\mathscr{H}_{m,n}(q,\bm{u})$. Second, in a dual framework based on the vertex-operator realization of Schur functions, we derive an iterative formula on upper multipartitions, which may be viewed as a dual Murnaghan--Nakayama rule for the cyclotomic Hecke algebra.

Our second aim is to demonstrate the scope of these rules through several applications. We obtain a Regev-type formula for the cyclotomic Hecke algebra, establish a L\"ubeck--Prasad--Adin--Roichman-type formula that computes character values via partitions rather than multipartitions, and introduce the notion of multiple bitrace together with a general combinatorial formula. These applications illustrate the practical strength of the new Murnaghan--Nakayama rule and do not seem to be readily accessible from the earlier formulas. They also extend the corresponding formulas discussed in Subsection \ref{ss:Backgroud}.

\subsection{Main results} 
We now summarize the main results of the paper. In the Hecke algebra of type $A$, the Murnaghan--Nakayama rule is naturally expressed in terms of generalized ribbons. In the cyclotomic setting, the corresponding combinatorial objects become {\em generalized multi-ribbons}, namely $m$-tuples of generalized ribbons. When $m=1$, the cyclotomic Hecke algebra reduces to the Hecke algebra of type $A$, and generalized multi-ribbons reduce exactly to generalized ribbons. In this sense, our construction is a genuine extension of the type $A$ picture.

More broadly, our result fits into a unified pattern relating character formulas, symmetric functions, and skew-diagram combinatorics. Table \ref{t:M-N} summarizes Murnaghan--Nakayama rules for several algebras. The present work completes the last line of the table and simultaneously generalizes the cases of the symmetric group, Hecke algebras of types $A$ and $B$, and the complex reflection group $\mathscr{W}_{m,n}$.

\begin{table}
    \centering  
    \renewcommand{\arraystretch}{2}
{\small\begin{tabularx}{\textwidth}{>{\centering\arraybackslash}X>{\centering\arraybackslash}X>{\centering\arraybackslash}X}
\toprule
 \bf{Algebras} & \bf{Symmetric functions} & \bf{Skew diagrams} \\ \midrule
symmetric group $\mathfrak{S}_n$ & Schur function $s_{\la}$ & ribbon \\ \midrule
Hecke algebra (type $A$) $\mathscr{H}_n$ & Schur function $s_{\la}$ & generalized ribbon \\ \midrule
Hecke algebra (type $B$) $\mathscr{B}_n$ & double Schur function $s_{(\la,\mu)}$ & double generalized ribbon \\ \midrule
Sergeev algebra $\mathscr{C}_n\rtimes \mathfrak{S}_n$ & Schur $Q$-function $Q_{\la}$ & double strip \\ \midrule
Hecke--Clifford algebra $\mathscr{C}_n\rtimes \mathscr{H}_n$ & Schur $Q$-function $Q_{\la}$ & generalized double strip \\ \midrule
complex reflection group $\mathscr{W}_{m,n}$ & multi-Schur function $s_{\bm{\la}}$ & multi-ribbon \\ \midrule
{\em cyclotomic Hecke algebra} $\mathscr{H}_{m,n}(q,\bm{u})$ & {\em multi-Schur function} $s_{\bm{\la}}$ & {\em generalized multi-ribbon} \\ \bottomrule
\end{tabularx}}
\caption{Murnaghan--Nakayama rules for various algebras}\label{t:M-N}
\end{table}

Our first main result is a recursive Murnaghan--Nakayama rule for irreducible characters of the cyclotomic Hecke algebra, formulated directly on Shoji's standard elements.

\medskip
\begin{mainthm}[Corollary \ref{t:m-n2}]\label{t:A}
    Let $\bm{\la}$ and $\bm{\mu}$ be two multipartitions of $n$, and let $\chi_{\bm{\mu}}^{\bm{\la}}$ be the value of the irreducible character indexed by $\bm{\la}$ at the standard element $g(\bm{\mu})$. Then
\begin{align}
\chi_{\bm{\mu}}^{\bm{\la}}=\sum_{\bm{\nu}}\frac{q^{\mu^{(r)}_j}}{q-q^{-1}}\wt(\bm{\la/\nu};q^{-2},\bm{u},r)\chi^{\bm{\nu}}_{\bm{\mu}\setminus \mu^{(r)}_j},
\end{align}
where the sum runs over all multipartitions $\bm{\nu}\subset \bm{\la}$ such that $\bm{\la/\nu}$ is a $\mu^{(r)}_j$-generalized multi-ribbon. Here $\bm{\mu}\setminus \mu^{(r)}_j$ is obtained from $\bm{\mu}$ by deleting the $j$-th row of the $r$-th partition, and $\wt(\bm{\la/\nu};q^{-2},\bm{u},r)$ is the explicit weight defined in \eqref{e:defwt}.
\end{mainthm}

\medskip

Iterating Theorem \ref{t:A} yields a tableau-type expansion formula: the character value $\chi_{\bm{\mu}}^{\bm{\la}}$ can be written as a weighted sum over generalized multi-ribbon tableaux of shape $\bm{\la}$ and content $\bm{\mu}$ (Corollary \ref{t:m-n3}). This provides a fully combinatorial global formula, complementing the recursive version above.

In addition to the above (``lowering'') recursion, we establish a dual iterative formula in the opposite direction (Corollary \ref{c:dualMN2}), derived from the vertex operator realization of Schur functions. This dual formula expresses a character as a linear combination of characters indexed by multipartitions with strictly smaller upper data, and therefore gives a complementary computational mechanism to the Murnaghan--Nakayama rule.

We next turn to applications. The first application is a Regev-type formula for the permutation super representation of $\mathscr{H}_{m,n}(q,\bm{u})$.

\medskip
\begin{mainthm}[Theorem \ref{t:Regev}]\label{t:B}
    Let $\bm{\mu}$ be a $m$-multipartition of $n$. Then the value of the character of the permutation super representation on the standard element $g(\bm{\mu})$ is given by
\begin{align}\label{e:D}
\begin{split}
    &\chi^{\bm{u},q}_{\bm{k}|\bm{l};n}(g(\bm{\mu}))\\
=&q^{n-l(\bm{\mu})}\prod_{r=1}^m\prod_{j=1}^{l(\mu^{(r)})}\sum_{(\bm{\alpha};\bm{\beta})\in \mathfrak{C}(\mu^{(r)}_j;\bm{k}\mid\bm{l})} \bm{u}^r_{\overrightarrow{l}(\bm{\alpha};\bm{\beta})}(1-q^{-2})^{l(\balp;\bbeta)-1}(-q^{-2})^{|\bbeta|-l(\bbeta)}\prod_{i=1}^m\binom{k_i}{l(\alpha^{(i)})}\binom{l_i}{l(\beta^{(i)})}.
\end{split}
\end{align}
Here $\mathfrak{C}(\mu^{(r)}_j;\bm{k}\mid\bm{l})$ is defined in Sec. \ref{s:Regev}, $\overrightarrow{l}(\bm{\alpha};\bm{\beta}):=\max\{1\leq i\leq m\mid l(\alpha^{(i)})+l(\beta^{(i)})>0\}$, and $l(\balp;\bbeta):=l(\balp)+l(\bbeta)$.
\end{mainthm}
\medskip

A special case of Theorem \ref{t:B} yields a hook-sum formula (Cor. \ref{t:special-hook}) for irreducible characters indexed by multipartitions all of whose components are hooks. We also provide an alternative proof of this formula in Appendix \ref{app}.

Our second application is a cyclotomic Hecke algebra analogue of the L\"ubeck--Prasad--Adin--Roichman formula. It expresses irreducible characters indexed by multipartitions in terms of data attached to the corresponding partition with empty $m$-core.

\medskip
\begin{mainthm}[Theorem \ref{t:AR}]\label{t:C}
    Let $\bla=(\la^{(1)},\la^{(2)},\cdots,\la^{(m)})$ and $\bmu\in \mathcal{P}_{n,m}$, and let $\la$ be the unique partition with empty $m$-core and $m$-quotient $\bla$. Then
\begin{align}\label{e:E}
    \chi^{\bla}_{\bmu}=q^{n-l(\bmu)}\sum_{r=1}^m\sum_{\btau\in\mathcal{P}_m(\mu^{(r)})}\eta^{\la}_{\btau;r},
\end{align}
where $\eta^{\la}_{\btau;r}$ is the partition-based quantity introduced in Definition \ref{d:eta'}.
\end{mainthm}
\medskip

Under the specialization $\bm{u}=\bm{\zeta}=(1,\zeta,\zeta^2,\cdots,\zeta^{m-1})$, $q=1$, and $\bmu=(\varnothing,\cdots,\varnothing,\mu)$, Theorem E specializes to the Adin-Roichman formula for the complex reflection group $\mathscr{W}_{m,n}$.

Finally, we introduce the notion of {\em multiple bitrace} for $\mathscr{H}_{m,n}(q,\bm{u})$ and derive an explicit combinatorial formula. This may be viewed as a cyclotomic analogue of the bitrace formulas known in types $A$ and related settings.

\medskip
\begin{mainthm}[Theorem \ref{t:mbtr}]\label{t:D}
    Let $\bm{\mu}$ and $\bm{\nu}$ be two multipartitions of $n$. Then
\begin{align}\label{e:F}
        \mbtr(\bm{\mu},\bm{\nu})=\frac{q^{2n}}{(q-q^{-1})^{l(\bm{\mu})+l(\bm{\nu})}}\sum_{\bm{M},\bm{N},\bm{A}}\bm{u}_{\bm{M}}\bm{u}_{\bm{N}}\prod_{i=1}^m\wt_{q^{-2}}(A^{(i)}),
\end{align}
where the sum runs over pairs of $m$-multimatrices $(\bm{M},\bm{N})$ and $m$-tuples of matrices $\bm{A}$ satisfying the conditions stated in Theorem \ref{t:mbtr}. Here $\bm{u}_{\bm{M}}$ and $\wt_q(A^{(i)})$ are the corresponding weight polynomials attached to $\bm{M}$ and $A^{(i)}$, respectively.
\end{mainthm}
\medskip

By specialization, the left-hand side and right-hand side of \eqref{e:F} reduce to $\sum_{\bm{\la}\in \mathcal{P}_{n,m}}\phi^{\bm{\la}}_{\bm{\mu}}\overline{\phi^{\bm{\la}}_{\bm{\nu}}}$ and $\delta_{\bmu,\bm{\nu}}m^{l(\bmu)}\prod_{i=1}^m z_{\mu^{(i)}}$, respectively, and thus Theorem F recovers the second orthogonality relation for irreducible characters of $\mathscr{W}_{m,n}$:
\begin{align}
    \sum_{\bm{\la}\in \mathcal{P}_{n,m}}\phi^{\bm{\la}}_{\bm{\mu}}\overline{\phi^{\bm{\la}}_{\bm{\nu}}}=
    \begin{cases}
      m^{l(\bmu)}\prod\limits_{i=1}^m z_{\mu^{(i)}}, &\text{if $\bmu=\bm{\nu}$},\\
      0, &\text{otherwise},
    \end{cases}
\end{align}
where $z_{\lambda}:=\prod_{i\geq 1}i^{m_i(\lambda)}m_i(\lambda)!$ for a partition $\lambda$.

\subsection{Organization} 
The layout of this paper is as follows. In the next section, we will briefly review key terminologies related to multipartitions and symmetric functions indexed by multipartitions, along with a recapitulation of Shoji's Frobenius formula. In Sec. \ref{s:M-N}, we will present our Murnaghan--Nakayama rule, accompanied by several straightforward examples. Sec. \ref{s:special} focuses on three special cases of the Murnaghan--Nakayama rule, demonstrating how our rule recovers existing formulas for the complex reflection group of type $G(m, 1, n)$, as well as for Hecke algebras of type $A$ and type $B$. In Sec. \ref{s:OV}, we provide an alternative computational formula (dual Murnaghan--Nakayama rule) for the irreducible characters in a dual picture, utilizing the vertex operator realization of Schur functions. In Sec. \ref{s:Regev}, we present the first application of our Murnaghan--Nakayama rule: the Regev formula for the cyclotomic Hecke algebra. We extend the work of L\"ubeck--Prasad--Adin--Roichman to the setting of the cyclotomic Hecke algebra in Sec. \ref{s:LPAR}. Sec. \ref{s:mbtr} concentrates on the multiple bitrace of irreducible characters, followed by the presentation of a combinatorial formula for it. In Appendix \ref{app:sagecode}, we show how to use the Sagemath program \cite{Sage} to compute the irreducible character tables $(\chi^{\bm\la}_{\bm\mu})_{\bla,\bmu\in\mathcal P_{n,m}}$ via our Murnaghan--Nakayama rule.

\section{Preliminaries}

\subsection{Multipartitions and corresponding symmetric functions}
The terminologies and notations about classical {\it partitions}, such as {\it length, weight, diagram}, etc, can be found in \cite[Chapter I]{Mac}. Now we introduce the corresponding notions and definitions in the context of {\em multipartitions}. A $m$-multipartition (resp. $m$-multicomposition) of $n$ is an ordered tuple $\bm{\la}=(\la^{(1)},\la^{(2)},\cdots,\la^{(m)})$ of partitions (resp. compositions) $\la^{(i)}$ such that $|\bm{\la}|=\sum_{i=1}^{m}|\la^{(i)}|=n$. We denote by $\mathcal{P}_{n,m}$ the set of all $m$-multipartitions of $n$. Let $C_{n,m}$ denote the set of $m$-tuples of nonnegative integers summing to $n$, i.e., an element in $C_{n,m}$ has the form $(c_1, c_2, \cdots, c_m)$ with $c_i$ nonnegative integer such that $\sum_{i=1}^{m}c_i=n$. 
Occasionally, we will omit $m$, just say $\bla$ is a multipartition (or multicomposition), for short unless special emphasis. We usually write $\bm{\la}$ in bold to emphasize that $\bm{\la}$ is a multipartition. Define $l(\bm{\la})=\sum_{i\geq1}^ml(\la^{(i)})$. The {\it diagram} of the multipartition $\bm{\la}$ is naturally the sequence of the consecutive Young diagrams of $\la^{(i)}$, $i=1,2,\cdots,m$. Here is an example for a $5$-multipartition of $24$ with length $11$.
\begin{align}\label{e:la}
\big((3,2,1),(3,3,1,1),\varnothing,(2,2,2),(4)\big)=\left( \begin{array}{ccccc}
\begin{tikzpicture}[scale = 0.33]
  \begin{scope}
    \clip (0,0) -| (1,1) -| (2,2) -| (3,3) -- (0,3) -- (0,0);
    \draw [color=black!25] (0,0) grid (3,3);
  \end{scope}
  \draw [thick] (0,0) -| (1,1) -| (2,2) -| (3,3) -- (0,3) -- (0,0);
\end{tikzpicture},&
\begin{tikzpicture}[scale = 0.33]
  \begin{scope}
    \clip (0,0) -| (1,2) -| (3,4) -- (0,4) -- (0,0);
    \draw [color=black!25] (0,0) grid (3,4);
  \end{scope}
  \draw [thick] (0,0) -| (1,2) -| (3,4) -| (0,4) -- (0,0);
\end{tikzpicture},&
    \varnothing,&
    \begin{tikzpicture}[scale = 0.33]
  \begin{scope}
    \clip (0,0) -| (2,3) -- (0,3) -- (0,0);
    \draw [color=black!25] (0,0) grid (2,3);
  \end{scope}
  \draw [thick] (0,0) -| (2,3)  -- (0,3) -- (0,0);
\end{tikzpicture},&
\begin{tikzpicture}[scale = 0.33]
  \begin{scope}
    \clip (0,0) -| (4,1) -- (0,1) -- (0,0);
    \draw [color=black!25] (0,0) grid (4,1);
  \end{scope}
  \draw [thick] (0,0) -| (4,1) -- (0,1) -- (0,0);
\end{tikzpicture}
\end{array}  \right).
\end{align}


For two arbitrary multipartitions $\bm{\la},\bm{\mu}\in \mathcal{P}_{n,m}$, we say $\bm{\mu}\subset\bm{\la}$ if $\mu^{(i)}\subset\la^{(i)}$ for all $1\leq i\leq m$. We define the {\it $m$-multi skew diagram} of $\bm{\la}/\bm{\mu}$ (denoted by $\D(\bla/\bmu)$) as the sequence of the consecutive skew diagrams of $\la^{(i)}/\mu^{(i)}$. For instance, we choose $\bm{\la}$ in the above figure and $\bm{\mu}=((2,1,1),(2,1,1),\varnothing,(2,1),(2))$. Clearly, $\bmu\subset\bla$. The skew diagram of $\bla/\bmu$ is presented as below:
\begin{align}\label{e:la/mu}
\D(\bm{\la/\mu})=\left( \begin{array}{ccccc}
\begin{tikzpicture}[scale = 0.35]
  \begin{scope}
    \clip (0,0) -| (1,1) -| (2,2) -| (3,3) -- (0,3) -- (0,0);
    \draw [color=black!25] (0,0) grid (3,3);
  \end{scope}
  \draw [thick] (0,0) -| (1,1) -| (2,2) -| (3,3) -- (0,3) -- (0,0);
      \node [draw, circle, fill = white, inner sep = 1.5pt] at (1.5,1.5) { };
      \node [draw, circle, fill = white, inner sep = 1.5pt] at (2.5,2.5) { };
\end{tikzpicture},&
\begin{tikzpicture}[scale = 0.35]
  \begin{scope}
    \clip (0,0) -| (1,2) -| (3,4) -- (0,4) -- (0,0);
    \draw [color=black!25] (0,0) grid (3,4);
  \end{scope}
  \draw [thick] (0,0) -| (1,2) -| (3,4) -- (0,4) -- (0,0);
     \draw [thick, rounded corners] (1.5,2.5) -- (2.5,2.5) -- (2.5,3.5);
       \draw [color=black,fill=black,thick] (2.5,3.5) circle (.6ex); 
      \node [draw, circle, fill = white, inner sep = 1.5pt] at (0.5,0.5) { };
      \node [draw, circle, fill = white, inner sep = 1.5pt] at (1.5,2.5) { };
\end{tikzpicture},&
\varnothing,&
 \begin{tikzpicture}[scale = 0.35]
  \begin{scope}
    \clip (0,0) -| (2,3) -- (0,3) -- (0,0);
    \draw [color=black!25] (0,0) grid (2,3);
  \end{scope}
  \draw [thick] (0,0) -| (2,3)  -- (0,3) -- (0,0);
  \draw [thick, rounded corners] (0.5,0.5) -- (1.5,0.5) -- (1.5,1.5);
       \draw [color=black,fill=black,thick] (1.5,1.5) circle (.6ex); 
      \node [draw, circle, fill = white, inner sep = 1.5pt] at (0.5,0.5) { };
\end{tikzpicture},&
\begin{tikzpicture}[scale = 0.35]
  \begin{scope}
    \clip (0,0) -| (4,1) -- (0,1) -- (0,0);
    \draw [color=black!25] (0,0) grid (4,1);
  \end{scope}
  \draw [thick] (0,0) -| (4,1) -- (0,1) -- (0,0);
  \draw [thick, rounded corners] (2.5,0.5) -- (3.5,0.5);
       \draw [color=black,fill=black,thick] (3.5,0.5) circle (.6ex); 
      \node [draw, circle, fill = white, inner sep = 1.5pt] at (2.5,0.5) { };
\end{tikzpicture}
\end{array}  \right).
\end{align}

Let $\Lambda_{\mathbb{Q}}$ be the ring of symmetric functions in the $x_n$ ($n\in\mathbb N$) over $\mathbb{Q}$. As is well known, the ring $\Lambda_{\mathbb{Q}}$ has several linear bases indexed by partitions such as power sum functions $p_{\la}(x)$, elementary symmetric functions $e_{\la}(x)$, complete symmetric functions $h_{\la}(x)$ and Schur functions $s_{\la}(x)$, etc., and the reader 
is referred to \cite[Chap. I]{Mac} for their definitions. There exists the canonical Hall inner product on $\Lambda_{\mathbb{Q}}$ defined by:
\begin{align}\label{e:z}
\langle p_{\la}(x), p_{\mu}(x) \rangle=\delta_{\la,\mu}z_{\la},
\end{align}
where $z_{\lambda}=\prod_{i\geq 1}i^{m_i(\lambda)}m_i(\lambda)!$, under which the Schur functions forms an orthonormal basis. 

Let us consider the tensor product $\Lambda_{\mathbb{Q}}^{\otimes m}$. Now we extend these notions to the setting of $m$-multipartitions.
Let $x^{(i)}_j$ ($1\leq i\leq m$, $j\geq 1$) be $m$ sets of variables.  For each integer $k\geq 1$ and $1\leq i\leq m$, let
\begin{align}
p_{k}^{(i)}(x)=\sum_{j=1}^{m}\zeta^{-ij}p_{k}(x^{(j)}),
\end{align}
where $p_k(x^{(j)})$ are the ordinary $k$-th power sum symmetric functions in the variables $x^{(j)}_1,x^{(j)}_2,\cdots$, and throughout the paper $\zeta$ is a fixed primitive $m$-th root of unity. We define the power sum function $p_{\bm{\la}}$ associated with $m$-multipartition $\bm{\la}=(\la^{(1)},\la^{(2)},\cdots,\la^{(m)})$ by
\begin{align}\label{e:def-p}
p_{\bm{\la}}(x)=\prod_{i=1}^{m}\prod_{j=1}^{l(\la^{(i)})}p^{(i)}_{\la^{(i)}_j}(x).
\end{align}
Similarly, the Schur function $s_{\bm{\la}}(x)$ associated with $\bm{\la}$ is defined by
\begin{align}
s_{\bm{\la}}(x)=\prod_{i=1}^{m}s_{\la^{(i)}}(x^{(i)}),
\end{align}
where $s_{\la^{(i)}}(x^{(i)})$ denotes the ordinary Schur function associated with partition $\la^{(i)}$ in the variables $x^{(i)}_1,x^{(i)}_2,\cdots$. We refer to $s_{\bm{\la}}(x)$ as {\em multi-Schur functions}.

Clearly the tensor products of basis elements from a given basis of $\Lambda_{\mathbb{Q}}$ will form
a basis for the tensor product space $\Lambda_{\mathbb{Q}}^{\otimes m}$. In particular, $\{s_{\bm{\la}}\mid\bm{\la}\in \mathcal{P}_{n,m}\}$ forms a basis in $\Lambda_{\mathbb{Q}}^{\otimes m}$.
We canonically extend the inner product of $\Lambda_{\mathbb{Q}}$ to $\Lambda_{\mathbb{Q}}^{\otimes m}$. Namely,
\begin{align}\label{e:tensor}
\langle u_1\otimes\cdots\otimes u_m, v_1\otimes\cdots\otimes v_m \rangle_{\otimes}:=\prod_{i=1}^m\langle u_i,v_i \rangle
\end{align}
for any two sets of elements $\{u_i\}$, $\{v_i\}$ in $\Lambda_{\mathbb{Q}}$ and extend bilinearly.
Then we have
\begin{align}\label{e:orth}
\langle s_{\bm{\la}}, s_{\bm{\mu}} \rangle_{\otimes}=\delta_{\bm{\la}\bm{\mu}}.
\end{align}
This means $\{s_{\bm{\la}}\mid\bm{\la}\in \mathcal{P}_{n,m}\}$ forms an orthonormal basis in $\Lambda_{\mathbb{Q}}^{\otimes m}$.

\subsection{Frobenius-type formula for the cyclotomic Hecke algebra}\label{ss:Fro}
Let $\mathscr{W}_{m,n}$ be the complex reflection group of type $G(m, 1, n)$ (cf. \cite{ST}) (called complex reflection group for short), which is 
generated by $s_1,s_2,\cdots, s_n$ with the relations
\begin{align}
\begin{split}
&s_{1}^{m}=1, \quad s_{2}^2=\cdots=s^2_{n}=1,\\
&s_1s_2s_1s_2=s_2s_1s_2s_1,\\
&s_is_j=s_js_i, \quad\text{if $|i-j|>1$},\\
&s_is_{i+1}s_i=s_{i+1}s_is_{i+1}, \quad\text{for $2\leq i<n$}.
\end{split}
\end{align}
Denote by $\mathfrak{S}_n$ the symmetric group on $n$ letters. It is well-known that $\mathscr{W}_{m,n}=(\mathbb{Z}/m\mathbb{Z})^n \rtimes \mathfrak{S}_n$, where $s_2,\cdots, s_n$ are respectively the transpositions $(1,2), \cdots, (n-1,n)$ in $\mathfrak{S}_n$ and $s_1$ corresponds to the generator of the cyclic group. Several Murnaghan--Nakayama type rules for the complex reflection group $\mathscr{W}_{m,n}$ were
found in \cite{Osi, St, AK, GJ}.

Let $\mathbb{K}=\mathbb{C}(q,\bm{u})$ be the field of rational functions in indeterminates $q$ (called {\em Hecke parameter}) and $\bm{u}=(u_1,\cdots,u_m)$ (called {\em cyclotomic parameters}). The cyclotomic Hecke algebra $\mathscr{H}_{m,n}(q,\bm{u})$ \cite{AK,BM,C} (also called the Ariki-Koike algebra) associated to $\mathscr{W}_{m,n}$ is the unital associative $\mathbb{K}$-algebra generated by $g_1, g_2, \cdots, g_n$ and subject to relations\footnote{The presentation is slightly different from that in \cite{AK}. They are isomorphic by the substitutions: $g_i\rightarrow q^{-1}g_i$, $u_i\rightarrow q^{-1}u_i$ and $q^2\rightarrow q$.}
\begin{align}
\begin{split}
&(g_1-u_1)\cdots(g_1-u_m)=0,\\
&g_1g_2g_1g_2=g_2g_1g_2g_1,\\
&g_i^2=(q-q^{-1})g_i+1, \quad \text{for $2\leq i\leq n$},\\
&g_ig_j=g_jg_i, \quad \text{for $|i-j| > 1$},\\
&g_ig_{i+1}g_i = g_{i+1}g_ig_{i+1},\quad \text{for $2 \leq i < n$}.
\end{split}
\end{align}
It is known \cite{AK} that $\mathscr{H}_{m,n}(q,\bm{u})$ is a free $\mathbb{K}$-module of rank $n!m^n$. The algebra $\mathscr{H}_{m,n}(q,\bm{u})$ reduces to $\mathscr{W}_{m,n}$ by specializing $q=1$ and $\bm{u}=\bm{\zeta}=(1,\zeta,\cdots,\zeta^{m-1})$, i.e., $u_i=\zeta^{i-1}$, $i=1,2,\cdots,m$.

Let us recall the irreducible character table of $\mathscr{H}_{m,n}(q,\bm{u})$. It is well known that the set of conjugacy classes in $\mathscr{W}_{m,n}$ is parametrized by $\mathcal{P}_{n,m}$ \cite{AK}.
For a survey on the representation theory of $\mathscr{H}_{m,n}(q,\bm{u})$ and $q$-Schur algebras, see \cite{Mat}.

Let $A$ be the matrix of degree $m$ whose $(a,b)$-entry is $u^a_b$  for $0 \leq a\leq m-1$ and $1\leq b \leq m$. 
Its determinant $\Delta=\det A=\prod_{i>j}(u_i-u_j)$ is the Vandermonde determinant. Let $B=A^*$ be the adjoint matrix of $A$. We write $B=(v_{b,a}(\bm{u}))$, where $v_{b,a}(\bm{u})$ is a polynomial in $\bm{u}$. For $1 \leq a \leq m$, we define $F_b(X)$ as the generating polynomial of $v_{b, a}(\bm u)$ as follows:
\begin{align*}
F_b(X)=\sum_{0\leq a\leq m-1}v_{b,a}(\bm{u})X^a.
\end{align*}

Following \cite[Sec. 3.6]{Sho}, let $\mathscr{H}^{\natural}$ be the associative algebra over $\mathbb{K}$ generated by $g_2,\cdots, g_n$ and $\xi_{1},\cdots,\xi_{n}$ subject to the following relations:
\begin{align*}
&(g_i-q)(g_i + q^{-1})=0, \quad 2 \leq i \leq n;\quad (\xi_i-u_1)\cdots(\xi_i-u_m) = 0,\quad 1 \leq i \leq n;\\
&g_ig_{i+1}g_i = g_{i+1}g_ig_{i+1},\quad 2 \leq i < n;\quad g_ig_j = g_jg_i,\quad |i-j| \geq 2;\\
&\xi_i\xi_j = \xi_j\xi_i,\quad 1 \leq i, j \leq n;\quad g_j\xi_i = \xi_ig_j ,\quad i \neq j-1, j;\\
&g_j\xi_j = \xi_{j-1}g_j + \Delta^{-2}\sum_{a<b}(u_a-u_b)(q-q^{-1})F_a(\xi_{j-1})F_b(\xi_j),\quad 2\leq j\leq n;\\
&g_j\xi_{j-1} = \xi_{j}g_j-\Delta^{-2}\sum_{a<b}(u_a-u_b)(q-q^{-1})F_a(\xi_{j-1})F_b(\xi_j),\quad 2\leq j\leq n.
\end{align*}
Then $\mathscr{H}^{\natural}$ is isomorphic to $\mathscr{H}_{m,n}(q,\bm{u})$ owing to \cite[Theorem 3.7]{Sho}.

We define the {\it standard elements} $g(\bm{\mu})\in \mathscr{H}_{m,n}(q,\bm{u})$ as follows. Fix $m\ge 1$ and $n\ge 1$.  For any composition $\alpha=(\alpha_1,\ldots,\alpha_\ell)$ of $n$
(with $\alpha_k\ge 1$), set $b_0=0$ and $b_t=\alpha_1+\cdots+\alpha_t$ for $1\le t\le \ell$.
For each $t$, let $\mathscr{H}^{\natural}_{m,\alpha_t}$ be the cyclotomic Hecke algebra on
$\alpha_t$ letters (with generators $g_2,\ldots,g_{\alpha_t}$ and $\xi_1,\ldots,\xi_{\alpha_t}$),
and consider the block-embedding
\[
\iota_{\alpha}:\ 
\mathscr{H}^{\natural}_{m,\alpha_1}\otimes\cdots\otimes \mathscr{H}^{\natural}_{m,\alpha_\ell}
\ \hookrightarrow\ \mathscr{H}^{\natural}_{m,n}
\]
defined on generators by
\[
\iota_{\alpha}\big(g^{(t)}_{k}\big)=g_{b_{t-1}+k},\quad (2\le k\le \alpha_t),\qquad
\iota_{\alpha}\big(\xi^{(t)}_{k}\big)=\xi_{b_{t-1}+k},\quad (1\le k\le \alpha_t),
\]
where $g^{(t)}_{k},\xi^{(t)}_{k}$ denote the generators in the $t$-th tensor factor.
(Thus the images of different tensor factors commute, giving a single embedding.)

For $1\le i\le m$ and $r\ge 1$, define the $i$-cycle element in $\mathscr{H}^{\natural}_{m,r}$ by
\[
g_{r,i}:=\xi_r^{\,i-1}\,g_r g_{r-1}\cdots g_2,
\]
(with the convention $g_{1,i}:=\xi_1^{\,i-1}$).
Now for a partition $\mu=(\mu_1,\ldots,\mu_\ell)$ of $r$, define
\[
g(\mu,i):=\iota_{\mu}\big(g_{\mu_1,i}\otimes g_{\mu_2,i}\otimes\cdots\otimes g_{\mu_\ell,i}\big)
\ \in\ \mathscr{H}^{\natural}_{m,r},
\]
where $\iota_{\mu}$ stands for the above embedding with $\alpha=\mu$.

Finally, for a multipartition $\bm{\mu}=(\mu^{(1)},\ldots,\mu^{(m)})\in\mathcal{P}_{n,m}$, define
the standard element
\[
g(\bm{\mu}) := g(\mu^{(1)},1)\, g(\mu^{(2)},2)\cdots g(\mu^{(m)},m)\ \in\ \mathscr{H}^{\natural}_{m,n},
\]
and view it in $\mathscr{H}_{m,n}(q,\bm u)$ via the isomorphism
$\mathscr{H}^{\natural}\simeq \mathscr{H}_{m,n}(q,\bm u)$.

It was proved in \cite[Proposition 7.5]{Sho} that the characters $\chi^{\bm{\la}}$ of $\mathscr{H}_{m,n}(q,\bm{u})$ are determined completely by their values on $g(\bm{\mu})$ for all $\bm{\mu}\in \mathcal{P}_{n,m}$. For simplicity, we write $\chi^{\bm{\la}}_{\bm{\mu}}$ for $\chi^{\bm{\la}}(g(\bm{\mu}))$ in short. $\chi^{\bm{\la}}_{\bm{\mu}}$ forms the irreducible character table of $\mathscr{H}_{m,n}(q,\bm{u})$ when $\bm{\la}$ and $\bm{\mu}$ run over $\mathcal{P}_{n,m}$.

Recall the Hall-Littlewood function $q_{k}(x;t)$ defined by
\begin{align*}
q_0(x;t)&=1,\\
q_k(x;t)&=(1-t)\sum_{i\geq 1}x_i^k\prod_{j\neq i}\frac{x_i-tx_j}{x_i-x_j} \quad (k\geq 1).
\end{align*}
It is expedient to consider the generating function of $q_{k}(x;t)$:
\begin{align}\label{e:generating-q}
Q(z;t)=\sum_{k=0}^{\infty}q_{k}(x;t)z^k=\prod_{i}\frac{1-x_itz}{1-x_iz}.
\end{align}
Now we introduce a family of symmetric functions $q_{\bm{\mu}}(x;q,\bm{u})$ as follows. For each $1\leq i\leq m$, let
\begin{align*}
q^{(i)}_{n}(x;q,\bm{u})=\frac{q^n}{q-q^{-1}}\sum_{c\in C_{n,m}}u^i_c\prod_{j=1}^{m}q_{c_j}(x^{(j)};q^{-2}) \qquad\text{for $n\geq 1$,}
\end{align*}
where $c=(c_1,\cdots,c_m)$ runs through $C_{n, m}$, and $u_c$ denotes $u_k$ for the largest integer $k$ such that $c_k\neq 0$. For convention $q^{(i)}_{0}(x;q,\bm{u}):=\frac{1}{q-q^{-1}}$.
Then define $q_{\bm{\mu}}(x;q,\bm{u})$ by
\begin{align*}
q_{\bm{\mu}}(x;q,\bm{u})=\prod_{i=1}^m\prod_{j=1}^{l(\mu^{(i)})}q^{(i)}_{\mu_j^{(i)}}(x;q,\bm{u}).
\end{align*}

The Frobenius formula for the characters of the cyclotomic Hecke algebra can be stated as follows \cite[Theorem 6.14]{Sho}.
\begin{thm}\label{t:fro}
Let $\bm{\la},\bm{\mu}\in \mathcal{P}_{n,m}$, then the irreducible character $\chi^{\bm{\la}}_{\bm{\mu}}$ is determined by
\begin{align}\label{e:fro}
q_{\bm{\mu}}(x;q,\bm{u})=\sum_{\bm{\la}\in \mathcal{P}_{n,m}} \chi^{\bm{\la}}_{\bm{\mu}}s_{\bm{\la}}.
\end{align}
\end{thm}

\subsection{$\la$-ring theory} For the $\la$-ring theory, we mainly refer the reader to \cite{Knu}. A {\em pre-$\lambda$-ring} is a commutative ring $\R$ with identity $1$ endowed a series of operations $\lambda^i: \R\rightarrow \R$ for $i\in\mathbb{N}$, such that for all $r\in\R$, the formal power series
\begin{align}\label{E:c1}
\lambda_z(r)=\sum_{i=0}^{\infty}\lambda^i(r)z^i=\lambda^0(r)+\lambda^1(r)z+\lambda^2(r)z^2+\cdots
\end{align}
 satisfies 
\begin{align}\label{E:la1}
\lambda^0(r)=1,\quad \lambda^1(r)=r\quad \text{and}\quad \lambda_z(r+r')=\lambda_z(r)\lambda_z(r') \quad \text{for all $r,r'\in\R$}.
\end{align}
\eqref{E:la1} gives $\lambda_z(0)=1$ and $\lambda_z(r)=\lambda_z(-r)^{-1}$. A {\em $\la$-ring} is a pre-$\la$-ring attached more conditions involving expressions for $\la^i(xy)$ and $\la^i(\la^j(x))$. However, we will not delve further into the general theory of $\la$-ring here. See \cite[Chapter I]{Knu} for details.

Consider the ring of symmetric functions $\Lambda_{\mathbb{Q}}$, which is a free $\lambda$-ring on one generator, i.e., $\Lambda_{\mathbb{Q}}=\mathbb{Q}[e_1,e_2,\cdots]$ with $\lambda^i(e_1)=e_i$ where $e_i$ is the $i$-th elementary symmetric function. If we choose $x_1,x_2,\cdots$ as the variables. Then $\lambda^i(X)=e_i$. Here $X=e_1(x)=x_1+x_2+\cdots$. If we define $\sigma^i(r):=(-1)^i\lambda^i(-r)$ and
\begin{align*}
\sigma_z(r):=\sum_{i=0}^{\infty}\sigma^i(r)z^i=\sum_{i=0}^{\infty}
\lambda^i(-r)(-z)^i=\lambda_{-z}(-r).
\end{align*}
Then we have
\begin{align}
e_i(x)&=\lambda^i(X), \qquad\quad h_i(x)=\sigma^i(X),\\
E(x;z)&=\lambda_z(X),\qquad H(x;z)=\sigma_z(X),
\end{align}
where $E(x;z)$ and $H(x;z)$ are respectively the generating functions of $e_r(x)$ and $h_r(x)$ ($r$-th complete symmetric function) given by
\begin{align}\label{e:defe}
E(x;z):=\sum_{n\geq0}e_n(x)z^n&=\prod_{i\geq1}\left(\sum_{0\le j\le 1}(x_iz)^j\right)=\prod_{i\geq1}(1+x_iz),\\
\label{e:defh}
H(x;z):=\sum_{n\geq0}h_n(x)z^n&=\prod_{i\geq1}\left(\sum_{j\geq0}(x_iz)^j\right)=\prod_{i\geq1}\frac{1}{1-x_iz}.
\end{align}
Now we rewrite $q_{k}(x;t)$ by using the $\la$-ring notation. By \eqref{e:generating-q}, we have
\begin{align}
    Q(z;t)=\frac{\sigma_z(X)}{\sigma_z(tX)}=\sigma_z\big((1-t)X\big),
\end{align}
which further yields
\begin{align}\label{e:q-laring}
   q_{k}(x;t)=\sigma^k\big((1-t)X\big).
\end{align}

For any two partitions $\mu\subset\la$, we define the {\em skew Schur operation} $s^{\la/\mu}(X)$ in the $\la$-ring $\Lambda_{\mathbb{Q}}$ by
\begin{align}
    s^{\la/\mu}(X):=\det_{1\leq i,j\leq l(\la)}(\sigma^{\la_i-i-\mu_j+j}(X)).
\end{align}
We have the following properties of $s^{\la/\mu}(X)$:
\begin{align}\label{e:schurrule}
\begin{split}
s^{\la/\mu}(X+Y)&=\sum_{\mu\subset\nu\subset\la}s^{\la/\nu}(X)s^{\nu/\mu}(Y) \quad \text{(sum rule)}\\
s^{\la/\mu}(-X)&=(-1)^{|\la/\mu|}s^{\la^{t}/\mu^{t}}(X) \quad \text{(duality rule)}\\
s^{\la/\mu}(tX)&=t^{|\la/\mu|}s^{\la/\mu}(X) \quad \text{(homogeneity)}.
\end{split}
\end{align}
Here $\la^{t}$ is the {\it conjugate} of $\la$, obtained from $\la$ by transposing the diagram $\D(\la)$ along the main diagonal. These identities can be found in the literature, for instance \cite[p. 43, p. 72]{Mac}, \cite[(6.2)]{HR95}, etc.

\section{Murnaghan--Nakayama rule for the cyclotomic Hecke algebra}\label{s:M-N}
In this section, we will present our Murnaghan--Nakayama rule for the cyclotomic Hecke algebra.

\subsection{Murnaghan--Nakayama rule}
A skew diagram $\lambda/\mu$ is a {\it vertical} (resp. {\it horizontal}) {\it strip} if each row (resp. column) contains at most one box. We say a skew diagram of $\la/\mu$ is a {\it ribbon} (also known as {\it border strip} or {\it rim hook}) if it is connected and contains no $2\times 2$ blocks \begin{tikzpicture}
  \draw[step=0.2cm, color=black] (-0.2,-0.2) grid (0.2,0.2); 
\end{tikzpicture}. The {\it height} (denoted by $ht(\la/\mu)$) of a ribbon $\la/\mu$ is defined as the number of the rows in $\la/\mu$ minus $1$.  More generally, a {\it generalized ribbon} is a skew diagram without $2\times 2$ blocks (not necessarily connected). Any generalized ribbon is a union of its connected components, each of which is a ribbon. For a generalized ribbon $\la/\mu$ with connected components $\xi_1,\xi_2,\cdots,\xi_m$, define the height of $\la/\mu$ as the sum of the heights of connected components of $\la/\mu$, i.e., $ht(\la/\mu)=\sum_{i=1}^{m}ht(\xi_i)$, and use $cc(\la/\mu)$ to denote the number of connected components of $\la/\mu$. Define the {\it weight} of a generalized ribbon $\la/\mu$ by
\begin{align}
\wt(\la/\mu;t):=(1-t)^{cc(\la/\mu)}(-t)^{ht(\la/\mu)}.
\end{align}

The following lemma is essential for us to obtain our Murnaghan--Nakayama rule, with part of the proof techniques inspired by \cite[(6.7)]{HR95}.
\begin{lem}\label{t:qs}
Let $k$ be a non-negative integer, then
\begin{align}\label{e:qs}
q_k(x;t)s_{\mu}=\sum_{\la}\wt(\la/\mu;t)s_{\la}
\end{align}
summed over all partitions $\la\vdash |\mu|+k$ such that $\la/\mu$ is a generalized ribbon.
\end{lem}
\begin{proof}
Recall the identity in \cite[Corollary 3.2]{JL3} $\big(q=t=0,\bm{a}=(1,0,\cdots)~\text{and}~\bm{b}=(t,0,\cdots)\big)$ that
\begin{align}
\sigma^k\big((1-t)X\big)s_{\mu}(x)=\sum_{\mu\subset\la}s^{\la/\mu}(1-t)s_{\la}(x).
\end{align}
By \eqref{e:q-laring}, it suffices to prove $s^{\la/\mu}(1-t)=\wt(\la/\mu;t)$ if $\la/\mu$ is a generalized ribbon and $0$ otherwise. By \eqref{e:schurrule}, we have
\begin{align*}
s^{\la/\mu}(1-t)&=\sum_{\nu}s^{\nu/\mu}(1)s^{\la/\nu}(-t) \quad \text{(by sum rule)}\\
&=\sum_{\nu}(-1)^{|\la/\nu|}s^{\nu/\mu}(1)s^{\la^{t}/\nu^{t}}(t)\quad \text{(by duality rule)}\\
&=\sum_{\nu}(-t)^{|\la/\nu|}s^{\nu/\mu}(1)s^{\la^{t}/\nu^{t}}(1)\quad \text{(by homogeneity)}.
\end{align*}
It is well known that 
\begin{align*}
    s^{\nu/\mu}(1)=\begin{cases} 1& \text{if $\nu/\mu$ is a horizontal strip}\\ 0& \text{otherwise}\end{cases}\quad \text{and}\quad 
     s^{\la^{t}/\nu^{t}}(1)=\begin{cases} 1& \text{if $\la/\nu$ is a vertical strip}\\ 0& \text{otherwise}\end{cases},
\end{align*}
which yields that $s^{\la/\mu}(1-t)=0$ unless $\la$ is obtained from $\mu$ as follows: firstly add a horizontal strip to $\mu$ to obtain $\nu$ and then add a vertical strip to $\nu$ to get $\la$, i.e., $\la/\mu$ is a generalized ribbon. In this case, then each box in $\la/\mu$ satisfies one of the following (see Fig. \ref{F:hv}): 
\begin{enumerate}[label={(\roman*)}]
\item  There is a box of $\la/\mu$ immediately to the left (labeled by $h$);
\item There is a box of $\la/\mu$ immediately below it (labeled by $v$);
\item  Neither (i) nor (ii) holds (labeled by $*$).
\end{enumerate}
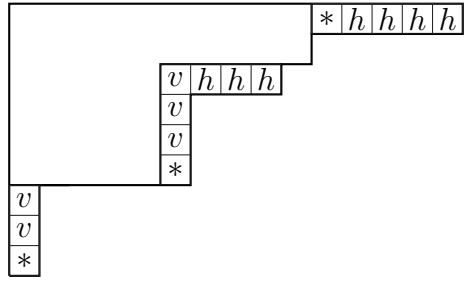
\begin{figure}
\begin{tikzpicture}[scale = 0.4]
      \draw [thick] (0,0) -| (1,3) -| (6,6) -| (9,7) -| (10,8) -| (15,9) -- (0,9) -- (0,0);
      \draw [thick] (0,3) -- (2,3);
      \draw [thick] (5,3) |- (9,7);
      \draw[thick] (10,8) -- (10,9);
      \begin{scope}
        \clip (0,0) -| (1,3) -- (0,3) -- (0,0);
        \draw [color=black!85] (0,0) grid (1,3);
      \end{scope}
      \begin{scope}
        \clip (5,3) -| (6,6) -| (9,7) -- (5,7) -- (5,3);
        \draw [color=black!85] (5,3) grid (9,7);
      \end{scope}
      \begin{scope}
        \clip (10,8) -| (15,9) -- (10,9) -- (10,8);
        \draw [color=black!85] (10,8) grid (15,9);
      \end{scope}
      \node [inner sep=2pt] at (0.5,0.5) { $*$}; 
       \node [inner sep=2pt] at (0.5,1.5) { $v$}; 
        \node [inner sep=2pt] at (0.5,2.5) { $v$}; 
         \node [inner sep=2pt] at (5.5,3.5) { $*$}; 
          \node [inner sep=2pt] at (5.5,4.5) { $v$}; 
           \node [inner sep=2pt] at (5.5,5.5) { $v$}; 
            \node [inner sep=2pt] at (5.5,6.5) { $v$}; 
             \node [inner sep=2pt] at (6.5,6.5) { $h$}; 
              \node [inner sep=2pt] at (7.5,6.5) { $h$}; 
               \node [inner sep=2pt] at (8.5,6.5) { $h$}; 
                \node [inner sep=2pt] at (10.5,8.5) { $*$}; 
                 \node [inner sep=2pt] at (11.5,8.5) { $h$}; 
                  \node [inner sep=2pt] at (12.5,8.5) { $h$}; 
                   \node [inner sep=2pt] at (13.5,8.5) { $h$}; 
                    \node [inner sep=2pt] at (14.5,8.5) { $h$}; 
    \end{tikzpicture}
    \caption{$\la$ is obtained from $\mu$ by adding a horizontal strip (to get $\nu$) and then adding a vertical strip (to $\nu$).}\label{F:hv}
    \end{figure}
Correspondingly, we have the following observations:
\begin{itemize}
    \item Each box in case (i) (resp. (ii)) must come from the application of the horizontal (resp. vertical) strip to $\mu$ and thus this box has weight $1$ (resp. $-t$);
    \item Each connected component (ribbon) $\xi_i$ in $\la/\mu$ has $c(\xi_i)-1$ (resp. $r(\xi_i)-1$) boxes in case (i) (resp. (ii));
\item Each connected component (ribbon) in $\la/\mu$ has exactly one box in case (iii), in which case this box has weight $1-t$.
\end{itemize}
Here $c(\xi_i)$ (resp. $r(\xi_i)$) presents the number of the columns (resp. rows) in $\xi_i$. Thus the product of these weights is exactly 
\begin{align*}
    (1-t)^{cc(\la/\mu)}\prod_{i}1^{c(\xi_i)-1}(-t)^{r(\xi_i)-1}=\wt(\la/\mu;t),
\end{align*}
which completes the proof.
\end{proof}

Suppose $\bm{\la}, \bm{\mu}\in \mathcal{P}_{n,m}$ and $\bm{\mu}\subset\bm{\la}$. We say $\bm{\la}/\bm{\mu}=(\la^{(1)}/\mu^{(1)},\cdots,\la^{(m)}/\mu^{(m)})$ is a {\it generalized $m$-multi-ribbon} if $\la^{(i)}/\mu^{(i)}$ is a generalized ribbon for all $1\leq i\leq m$. Analogously, we define
\begin{align}
cc(\bm{\la/\mu})&:=\sum_{i=1}^{m}cc(\la^{(i)}/\mu^{(i)}); \quad ht(\bm{\la/\mu}):=\sum_{i=1}^{m}ht(\la^{(i)}/\mu^{(i)}); \notag \\\label{e:defwt}
\wt(\bm{\la/\mu};t,\bm{u},r)&:=u^r_{\bm{\la/\mu}}\prod_{i=1}^{m}\wt(\la^{(i)}/\mu^{(i)};t)=u^r_{\bm{\la/\mu}}(1-t)^{cc(\bm{\la/\mu})}(-t)^{ht(\bm{\la/\mu})},
\end{align}
where $1\leq r\leq m$, $u_{\bm{\la/\mu}}:=u_s$ and $s$ presents the first non-empty entry in $\bm{\la/\mu}$ from right to left. For convention, we let $u_{\bla/\bla}:=1$. A generalized $m$-multi-ribbon $\bm{\la/\mu}$ is called a {\em $k$-generalized $m$-multi-ribbon} if $\sum_{i=1}^{m}|\la^{(i)}/\mu^{(i)}|=k$. Sometimes we will omit $m$, saying a generalized multi-ribbon for a generalized $m$-multi-ribbon in short.

\begin{exmp}
$\bm{\la/\mu}$ in \eqref{e:la/mu} is an $11$-generalized $5$-multi-ribbon. In this case, $cc(\bm{\la/\mu})=6$, $ht(\bm{\la/\mu})=1+1=2$ and $u_{\bla/\bmu}=u_5$. Thus, $\wt(\bm{\la/\mu};t,\bm{u},r)=u_5^rt^2(1-t)^6.$
\end{exmp}

With the above notations in hand, we are now in position to give the Murnaghan--Nakayama rule for the cyclotomic Hecke algebra, before which the following result is needed.
\begin{thm}\label{t:m-n}
Let $k$ be a non-negative integer and $1\leq r\leq m$, $\bm{\mu}\in \mathcal{P}_{n,m}$. Then
\begin{align}\label{e:m-n}
q^{(r)}_{k}(x;q,\bm{u})s_{\bm{\mu}}=\frac{q^k}{q-q^{-1}}\sum_{\bm{\la}}\wt(\bm{\la/\mu};q^{-2},\bm{u},r)s_{\bm{\la}}
\end{align}
summed over all multipartitions $\bm{\mu}\subset \bm{\la}\in\mathcal{P}_{n+k,m}$ such that $\bm{\la/\mu}$ is a $k$-generalized $m$-multi-ribbon.
\end{thm}
\begin{proof}
By definition, we have
\begin{align*}
q^{(r)}_{k}(x;q,\bm{u})s_{\bm{\mu}}=&\frac{q^k}{q-q^{-1}}\left(\sum_{c\in C_{k,m}}u^r_c\prod_{j=1}^{m}q_{c_j}(x^{(j)};q^{-2})\right)\prod_{i=1}^{m}s_{\mu^{(i)}}(x^{(i)})\\
=&\frac{q^k}{q-q^{-1}}\sum_{c\in C_{k,m}}u^r_c\prod_{j=1}^{m}q_{c_j}(x^{(j)};q^{-2})s_{\mu^{(j)}}(x^{(j)})\\
=&\frac{q^k}{q-q^{-1}}\sum_{c\in C_{k,m}}u^r_c\prod_{j=1}^{m}\left(\sum_{\la^{(j)}}\wt(\la^{(j)}/\mu^{(j)};q^{-2})s_{\la^{(j)}}(x^{(j)})\right)\quad \text{(by Lemma \ref{t:qs})}
\end{align*}
Swapping the inner summation and the product in the above equation implies that
\begin{align*}
q^{(r)}_{k}(x;q,\bm{u})s_{\bm{\mu}}=\frac{q^k}{q-q^{-1}}\sum_{c\in C_{k,m}}\sum_{\bm{\la}}\wt(\bm{\la/\mu};q^{-2},\bm{u},r)s_{\bm{\la}},
\end{align*}
where $\bm{\la}$ runs over all multipartitions $\bm{\mu}\subset \bm{\la}\in\mathcal{P}_{n+k,m}$ such that $\bm{\la/\mu}$ is a generalized multi-ribbon and $|\la^{(j)}/\mu^{(j)}|=c_j$, $j=1,2,\cdots,m$. This deduces \eqref{e:m-n}.
\end{proof}

As a straightforward result, the following corollary is an iterative formula for the irreducible characters, which is referred to as the Murnaghan--Nakayama rule.

\begin{cor}\label{t:m-n2}
Let $\bm{\la},\bm{\mu}\in \mathcal{P}_{n,m}$, $1\leq r\leq m$, $1\leq j\leq l(\mu^{(r)})$, then
\begin{align}\label{e:m-n2}
\chi_{\bm{\mu}}^{\bm{\la}}=\sum_{\bm{\nu}}\frac{q^{\mu^{(r)}_j}}{q-q^{-1}}\wt(\bm{\la/\nu};q^{-2},\bm{u},r)\chi^{\bm{\nu}}_{\bm{\mu}\setminus \mu^{(r)}_j}
\end{align}
summed over all multipartitions $\bm{\nu}\subset \bm{\la}$ such that $\bm{\la/\nu}$ is a $\mu^{(r)}_j$-generalized multi-ribbon, where $\bm{\mu}\setminus \mu^{(r)}_j$ is the multipartition obtained from $\bm{\mu}$ by removing the $j$-th part from the $r$-th partition $\mu^{(r)}$.
\end{cor}
\begin{proof}
By the Frobenius formula \eqref{e:fro}, we have
\begin{align}
q_{\bm{\mu}\setminus\mu^{(r)}_{j}}=\sum_{\bm{\nu}}\chi^{\bm{\nu}}_{\bm{\mu}\setminus\mu^{(r)}_{j}}s_{\bm{\nu}}.
\end{align}
Multiplying both sides by $q^{(r)}_{\mu^{(r)}_{j}}$ yields that
\begin{align}\label{e:I1}
\begin{split}
q_{\bm{\mu}}=&\sum_{\bm{\nu}}\chi^{\bm{\nu}}_{\bm{\mu}\setminus\mu^{(r)}_{j}}q^{(r)}_{\mu^{(r)}_{j}}s_{\bm{\nu}}\\
=&\sum_{\bm{\nu}}\chi^{\bm{\nu}}_{\bm{\mu}\setminus\mu^{(r)}_{j}}\frac{q^{\mu^{(r)}_{j}}}{q-q^{-1}}\sum_{\bm{\la}}\wt(\bm{\la/\nu};q^{-2},\bm{u},r)
s_{\bm{\la}} \quad\text{(by \eqref{e:m-n})}\\
=&\sum_{\bm{\la}}\sum_{\bm{\nu}}\chi^{\bm{\nu}}_{\bm{\mu}\setminus\mu^{(r)}_{j}}\frac{q^{\mu^{(r)}_{j}}}{q-q^{-1}}\wt(\bm{\la/\nu};q^{-2},\bm{u},r)
s_{\bm{\la}}
\end{split}
\end{align}
Note that $\{s_{\bm{\la}}\mid\bm{\la}\in \mathcal{P}_{n,m}\}$ forms an orthonormal basis in $\Lambda_{\mathbb{Q}}^{\otimes m}$. Comparing the coefficient of $s_{\bm{\la}}$ in \eqref{e:I1} and \eqref{e:fro} gives \eqref{e:m-n2}.
\end{proof}

\begin{exmp}
Suppose $\bm{\la}=((2,2),(4))$, $\bm{\mu}=((2,1),(3,2))$. Choose $r=j=2$. Then $\mu^{(r)}_{j}=2$ and there are $4$ $\bm{\nu}$'s such that $\bm{\la/\nu}$ is a $2$-generalized $2$-multi-ribbon. They are
\begin{align*}
\left(\begin{array}{cc}
\begin{tikzpicture}[scale = 0.33]
      \begin{scope}
        \clip (0,0) -| (2,2) -- (0,2) -- (0,0);
        \draw [color=black!25] (0,0) grid (2,2);
      \end{scope}
      \draw [thick] (0,0) -| (2,2) -- (0,2) -- (0,0);
      \draw [thick, rounded corners] (1.5,0.5) -- (1.5,1.5);
      \draw [color=black,fill=black,thick] (1.5,1.5) circle (.6ex); 
      \node [draw, circle, fill = white, inner sep = 1.5pt] at (1.5,0.5) { }; 
    \end{tikzpicture}~,&
    \begin{tikzpicture}[scale = 0.33]
    \begin{scope}
        \clip (0,0) -| (4,1) -- (0,1) -- (0,0);
        \draw [color=black!25] (0,0) grid (4,1);
      \end{scope}
      \draw [thick] (0,0) -| (4,1) -- (0,1) -- (0,0);
    \end{tikzpicture}
\end{array}\right)\quad
\left(\begin{array}{cc}
\begin{tikzpicture}[scale = 0.33]
      \begin{scope}
        \clip (0,0) -| (2,2) -- (0,2) -- (0,0);
        \draw [color=black!25] (0,0) grid (2,2);
      \end{scope}
      \draw [thick] (0,0) -| (2,2) -- (0,2) -- (0,0);
      \draw [thick, rounded corners] (0.5,0.5) -- (1.5,0.5);
      \draw [color=black,fill=black,thick] (1.5,0.5) circle (.6ex); 
      \node [draw, circle, fill = white, inner sep = 1.5pt] at (0.5,0.5) { }; 
    \end{tikzpicture}~,&
    \begin{tikzpicture}[scale = 0.33]
    \begin{scope}
        \clip (0,0) -| (4,1) -- (0,1) -- (0,0);
        \draw [color=black!25] (0,0) grid (4,1);
      \end{scope}
      \draw [thick] (0,0) -| (4,1) -- (0,1) -- (0,0);
    \end{tikzpicture}
\end{array}\right)\quad
\left(\begin{array}{cc}
\begin{tikzpicture}[scale = 0.33]
      \begin{scope}
        \clip (0,0) -| (2,2) -- (0,2) -- (0,0);
        \draw [color=black!25] (0,0) grid (2,2);
      \end{scope}
      \draw [thick] (0,0) -| (2,2) -- (0,2) -- (0,0);
    \end{tikzpicture}~,&
    \begin{tikzpicture}[scale = 0.33]
    \begin{scope}
        \clip (0,0) -| (4,1) -- (0,1) -- (0,0);
        \draw [color=black!25] (0,0) grid (4,1);
      \end{scope}
      \draw [thick] (0,0) -| (4,1) -- (0,1) -- (0,0);
      \draw [thick, rounded corners] (2.5,0.5) -- (3.5,0.5);
      \draw [color=black,fill=black,thick] (3.5,0.5) circle (.6ex); 
      \node [draw, circle, fill = white, inner sep = 1.5pt] at (2.5,0.5) { }; 
    \end{tikzpicture}
\end{array}\right)\quad
\left(\begin{array}{cc}
\begin{tikzpicture}[scale = 0.33]
      \begin{scope}
        \clip (0,0) -| (2,2) -- (0,2) -- (0,0);
        \draw [color=black!25] (0,0) grid (2,2);
      \end{scope}
      \draw [thick] (0,0) -| (2,2) -- (0,2) -- (0,0);
      \node [draw, circle, fill = white, inner sep = 1.5pt] at (1.5,0.5) { }; 
    \end{tikzpicture}~,&
    \begin{tikzpicture}[scale = 0.33]
    \begin{scope}
        \clip (0,0) -| (4,1) -- (0,1) -- (0,0);
        \draw [color=black!25] (0,0) grid (4,1);
      \end{scope}
      \draw [thick] (0,0) -| (4,1) -- (0,1) -- (0,0);
      \node [draw, circle, fill = white, inner sep = 1.5pt] at (3.5,0.5) { }; 
    \end{tikzpicture}
\end{array}\right).
\end{align*}
The weights of these $2$-generalized $2$-multi-ribbons are $u_1^2(1-t)(-t)$, $u_1^2(1-t)$, $u_2^2(1-t)$ and $u_2^2(1-t)^2$ respectively. Therefore,
\begin{align*}
\chi^{((2,2),(4))}_{((2,1),(3,2))}=&\frac{q^2}{q-q^{-1}}
\left(u_1^2(1-q^{-2})(-q^{-2})\chi^{((1,1),(4))}_{((2,1),(3))} +u_1^2(1-q^{-2})\chi^{((2),(4))}_{((2,1),(3))}\right.\\
&\left.+u_2^2(1-q^{-2})\chi^{((2,2),(2))}_{((2,1),(3))} +u_2^2(1-q^{-2})^2\chi^{((2,1),(3))}_{((2,1),(3))} \right).
\end{align*}
\end{exmp}

For given $\bm{\la}\in\mathcal{P}_{n,m}$, a {\it generalized $m$-multi-ribbon tableau} $\bm{\T}$ of shape $\bm{\la}$ is defined as a sequence
\begin{align}
\bm{\T}: \quad \bm{\varnothing}=\bm{\la^{(0)}}\subset\bm{\la^{(1)}}\subset\cdots\subset\bm{\la^{(l-1)}}\subset\bm{\la^{(l)}}=\bm{\la}
\end{align}
such that $\bm{\la^{(i)}/\la^{(i-1)}}$ $(1\leq i\leq l)$ is a generalized $m$-multi-ribbon.

Define an {\it order} on the set $\{i_j\mid i\geq 1, j\geq 1\}$ by
\begin{align}
a_b\leq c_d\quad \Longleftrightarrow \quad a< c~~\text{or } a=c \text{, } b\leq d.
\end{align}

A generalized $m$-multi-ribbon tableau $\bm{\T}$ can also be described by a tuple of tableaux, in which the cells of each partition are marked $i_j$ $(i\geq 1, j\geq 1)$, satisfying
\begin{enumerate}[label=(\roman*)]
\item the numbers $i_j$ inserted in the diagram of each partition of $\bm{\T}$ must increase weakly down each column and along each row from left to right;
\item the cells marked same number consist of a generalized $m$-multi-ribbon.
\end{enumerate}
The {\it content} of $\bm{\T}$ is defined as a multicomposition $\bm{c}=(c^{(1)},c^{(2)},\cdots)$ such that $c^{(i)}_j$ equals the number of cells marked $i_j$ in $\bm{\T}$. For a given $\bm{\T}$, introduce a map $f_{\bm{\T}}:~ \{i_j\mid i\geq1, j\geq1, i_j\in\bm{\T}\} \rightarrow \{1,2,\cdots,m\}$ by
\begin{align*}
f_{\bm{\T}}:~ i_j \longmapsto  \text{the rightmost entry in $\bm{\T}$ containing $i_j$}.
\end{align*}
For a generalized multi-ribbon tableau $\bm{\T}$ of shape $\bm{\la}$ with content $\bm{\mu}$, we define the {\it weight} of $\bm{\T}$ by
\begin{align}
\wt(\bm{\T};\bm{u},t):=(1-t)^{cc(\bm{\T})}(-t)^{ht(\bm{\T})}\prod_{i=1}^m\prod_{j=1}^{l(\mu^{(i)})}u^{i}_{f_{\bm{\T}}(i_j)}
\end{align}
where
\begin{align*}
cc(\bm{\T}):=\sum_{i=1}^m\sum_{j=1}^{l(\mu^{(i)})}cc(\bm{\xi}(i_j)),\quad ht(\bm{\T}):=\sum_{i=1}^m\sum_{j=1}^{l(\mu^{(i)})}ht(\bm{\xi}(i_j))
\end{align*}
and $\bm{\xi}(i_j)$ is the generalized multi-ribbon marked $i_j$.

\begin{exmp}
Below is a generalized multi-ribbon tableau $\bm{\T}$ of shape $\bm{\la}=((4,3),(2),(1,1))$ with content $\bm{\mu}=((3,2),(4),(2))$. In this case, $f_{\bm{\T}}(1_1)=1$, $f_{\bm{\T}}(1_2)=1$, $f_{\bm{\T}}(2_1)=3$, $f_{\bm{\T}}(3_1)=3$.
\begin{align*}
\bm{\T}=
\setlength{\arraycolsep}{15pt}
\left( \begin{array}{ccc}
\begin{tikzpicture}[scale=0.6]
   \begin{scope}
        \clip (0,0) -| (3,1) -| (4,2) -- (0,2) -- (0,0);
        \draw [color=black!25] (0,0) grid (4,2);
      \end{scope}
      \draw [thick] (0,0) -| (3,1) -| (4,2) -- (0,2) -- (0,0);
    \node at(0.5,0.5) {$1_1$};
    \node at(0.5,1.5) {$1_1$};
    \node at(1.5,1.5) {$1_1$};
    \node at(1.5,0.5) {$1_2$};
    \node at(2.5,1.5) {$1_2$};
    \node at(2.5,0.5) {$2_1$};
    \node at(3.5,1.5) {$2_1$};
    \end{tikzpicture},&
    \begin{tikzpicture}[scale=0.6]
    \begin{scope}
        \clip (0,0) -| (2,1) -- (0,1) -- (0,0);
        \draw [color=black!25] (0,0) grid (2,1);
      \end{scope}
      \draw [thick] (0,0) -| (2,1) -- (0,1) -- (0,0);
    \node at(0.5,0.5) {$2_1$};
\node at(1.5,0.5) {$3_1$};
    \end{tikzpicture}~,&
    \begin{tikzpicture}[scale=0.6]
    \begin{scope}
        \clip (0,0) -| (1,2) -- (0,2) -- (0,0);
        \draw [color=black!25] (0,0) grid (1,2);
      \end{scope}
      \draw [thick] (0,0) -| (1,2) -- (0,2) -- (0,0);
\node at(0.5,1.5) {$2_1$};
\node at(0.5,0.5) {$3_1$};
    \end{tikzpicture}
\end{array}  \right)
\end{align*}
$cc(\bm{\T})=1+2+4+2=9$, $ht(\bm{\T})=1$. So $\wt(\bm{\T};\bm{u},t)=u_1\cdot u_1\cdot u_3^2\cdot u_3^3(1-t)^{9}(-t)=u_1^2u_3^5(-t)(1-t)^9$.
\end{exmp}

\begin{cor}\label{t:m-n3}
Suppose $\bm{\la},\bm{\mu}\in\mathcal{P}_{n,m}$, we have the following general formula for $\chi^{\bm{\la}}_{\bm{\mu}}$ in terms of generalized $m$-multi-ribbon tableaux as follows:
\begin{align}
\chi^{\bm{\la}}_{\bm{\mu}}=\frac{q^n}{(q-q^{-1})^{l(\bm{\mu})}}\sum_{\bm{\T}}\wt(\bm{\T};\bm{u},q^{-2})
\end{align}
summed over all generalized $m$-multi-ribbon tableaux of shape $\bm{\la}$ with content $\bm{\mu}$.
\end{cor}
\begin{proof}
It is straightforward by using Corollary \ref{t:m-n2} repeatedly.
\end{proof}

\subsection{Examples and numerical relations} Now we apply our Murnaghan--Nakayama rule to compute some irreducible characters in special cases and establish an equality relation. 
We use the short-hand notation $\bvar^{1}_j$ (resp. $\bvar^{n}_j$) to present the multipartition with $(1^n)$ (resp. $(n)$) in the $j$-th entry and $\varnothing$ anywhere else.
\begin{exmp}
    We list the following simple results coming from our Murnaghan--Nakayama rule.
    \begin{itemize}
        \item If we denote by $f^{\la^{(i)}}$ the number of standard Young tableaux of shape $\la^{(i)}$, which is given by 
        \begin{align*}
            f^{\la^{(i)}}=\frac{|\la^{(i)}|!}{\prod_{x\in\la^{(i)}}h_x}
        \end{align*}
        where $h_x$ is the hook length at $x$. Then we have 
        \begin{align}
           \chi^{\bla}_{\bvar^{1}_j}=\prod_{i=1}^m f^{\la^{(i)}}u_i^{j|\la^{(i)}|}. 
        \end{align}
        \item Denote $L(\bla)=\#\{i\mid \la^{(i)}\neq\varnothing\}$ and $s=\max\{i\mid \la^{(i)}\neq \varnothing\}$. Then 
        \begin{align}
            \chi^{\bla}_{\bvar^{n}_j}=
        \begin{cases}
            q^{n-1}u^j_s(1-q^{-2})^{L(\bla)-1}(-q^{-2})^{\sum_{i}b_i}, &\text{if $\la^{(i)}=(a_i,1^{b_i})$, $i=1,\cdots,m$}\\
            0, &\text{otherwise}.
        \end{cases}
        \end{align}
        \item Denote $\mathcal{L}_{\bm{\mu}\bm{\mu}}:=\sum_{i=1}^m i l(\mu^{(i)})$.\footnote{We adopt this notation to be consistent with that in Theorem \ref{t:Sq}.} Then 
        \begin{align}
        \chi_{\bmu}^{\bvar^{1}_j}=&(-q)^{l(\bmu)-n}u_j^{\mathcal{L}_{\bmu\bmu}};\\
          \chi_{\bmu}^{\bvar^{n}_j}=&q^{n-l(\bmu)}u_j^{\mathcal{L}_{\bmu\bmu}}.  
        \end{align}
    \end{itemize}
    \end{exmp}
\begin{prop}\label{p:replace}
    Let $\bmu$ be an arbitrary $m$-multipartition of $n$. Suppose for $1\leq j\leq m$
    $$\bla=(\varnothing,\cdots,\varnothing,\la^{(j)},\cdots,\la^{(m)}) \quad \text{and}\quad \bla'=(\la^{(j)},\varnothing,\cdots,\varnothing,\la^{(j+1)},\cdots,\la^{(m)}).$$
    Then
    \begin{align}\label{e:replace}
\chi^{\bla}_{\bmu}=\chi^{\bla'}_{\bmu}\Big|_{u_1\rightarrow u_j}
    \end{align}
    where $f\Big|_{u_1\rightarrow u_j}$
 means replacing $u_1$ by $u_j$ in $f$. 
 \end{prop}
 \begin{proof}
    We observe that there exists a natural one-to-one correspondence between the set of generalized multi-ribbon tableaux of shape $\bla$ with content $\bmu$ and those of shape $\bla'$ with content $\bmu$. Applying Corollary \ref{t:m-n3} to both $\chi^{\bla}_{\bmu}$ and $\chi^{\bla'}_{\bmu}$, we find that the presence of $u_j$ in some monomial $M$ of $\chi^{\bla}_{\bmu}$ is equivalent to the presence of $i_j$ exclusively in $\la^{(j)}$. Under this one-to-one correspondence, this is further equivalent to the presence of $u_1$ in the corresponding monomial $M'$ of $\chi^{\bla'}_{\bmu}$. Moreover, all other factors remain identical. The only difference between the factors of $M$ and $M'$ lies in the replacement of $u_1$ with $u_j$. Thus, $\chi^{\bla}_{\bmu}$ can be obtained from $\chi^{\bla'}_{\bmu}$ by substituting $u_1$ with $u_j$.
 \end{proof}

 \section{Special cases of the Murnaghan--Nakayama rule}\label{s:special}
In this section, we will explain how our Murnaghan--Nakayama rule recovers the existing formulae as special cases. Three special cases of Corollary \ref{t:m-n2} will be considered case by case, which correspond to the complex reflection group of type $G(m, 1, n)$, the Iwahori--Hecke algebra in type $A$ and type $B$ respectively.
\subsection{Complex reflection group of type $G(m, 1, n)$}
As pointed out in the introduction, when we take the specialization $\bm{u}=\bm{\zeta}=(1,\zeta,\zeta^2,\cdots,\zeta^{m-1})$ and $q=1$, where $\zeta$ is the primitive
$m$-th root of unity, $\mathscr{H}_{m,n}(q,\bm{u})$ will reduce to $\mathscr{W}_{m,n}$. In this case, by definition, $\frac{q^{\mu_j^{(r)}}}{q-q^{-1}}\wt(\bm{\la/\mu};q^{-2},\bm{u},r)=0$ unless $cc(\bm{\la/\mu})=1$, i.e., $\bm{\la/\mu}$ has only one nonempty entry, in which the unique nonempty entry is a ribbon. We call such $\bm{\la/\mu}$ {\em $m$-multi-ribbon}. If we denote the unique nonempty entry by $s$, then $$\frac{q^{\mu_j^{(r)}}}{q-q^{-1}}\wt(\bm{\la/\mu};q^{-2},\bm{u},r)\Big|_{\substack{
\bm{u}=\bm{\zeta}\\ q=1}}=\zeta^{r(s-1)}(-1)^{ht(\bm{\la/\mu})}.$$
This explains why generalized multi-ribbons are not needed in the case of $\mathscr{W}_{m,n}$. Therefore, the Murnaghan--Nakayama rule for $\mathscr{W}_{m,n}$ (see also \cite[Prop. 2.2]{AK} and \cite{GJ}) can be stated as follows:
\begin{cor}
Let $\bm{\la},\bm{\mu}\in \mathcal{P}_{n,m}$, $1\leq r\leq m$, $1\leq j\leq l(\mu^{(r)})$. Denote $\phi^{\bm{\la}}$ as the irreducible character of $\mathscr{W}_{m,n}$ associated with the multipartition $\bm{\la}$. Then
\begin{align}\label{e:m-nW}
\phi_{\bm{\mu}}^{\bm{\la}}=\sum_{s=1}^{m}\sum_{\bm{\nu}}\zeta^{r(s-1)}(-1)^{ht(\bm{\la/\nu})}\phi^{\bm{\nu}}_{\bm{\mu}\setminus\mu^{(r)}_j}
\end{align}
summed over all $\bm{\nu}\in\mathcal{P}_{n-\mu^{(r)}_{j},m}$ such that the $s$-th skew diagram of $\bm{\la/\nu}$ is a $\mu^{(r)}_{j}$-ribbon, which is equivalent to
\begin{align}
\phi_{\bm{\mu}}^{\bm{\la}}=\sum_{\bm{\nu}}\zeta^{r(s-1)}(-1)^{ht(\bm{\la/\nu})}\phi^{\bm{\nu}}_{\bm{\mu}\setminus\mu^{(r)}_j}
\end{align}
summed over all $\bm{\nu}\in\mathcal{P}_{n-\mu^{(r)}_{j},m}$ such that $\bm{\la/\nu}$ is a $m$-multi-ribbon. Where $s$ is the unique nonempty entry in $\bm{\la/\nu}$.
\end{cor}

\subsection{Iwahori--Hecke algebra in type $A$}
The Iwahori--Hecke algebra in type $A$ (denoted by $\mathscr{H}_{n}(q)$) is the unital associative algebra
over $\mathbb{C}(q)$ generated by generators $T_{1}, T_{2},\ldots, T_{n-1}$ subject to the relations
\begin{align}\label{e:Hecke1}
T_{i}T_{j}&=T_{j}T_{i},  ~~~~\text{if}~ |i-j|>1,\\ \label{e:Hecke2}
T_{i}T_{i+1}T_{i}&=T_{i+1}T_{i}T_{i+1},\\ \label{e:Hecke3}
T_{i}^{2}&=(q-q^{-1})T_{i}+1.
\end{align}
It is clear that $\mathscr{H}_{m,n}(q,\bm{u})$ specializes to $\mathscr{H}_{n}(q)$ at $m=1$ and $\bm{u}=1$. In this case, a generalized multi-ribbon will reduce to a generalized ribbon, which is a disjoint union of ribbons as its connected components. Suppose there are $d$ ribbons in the generalized ribbons, then
$$\frac{q^{\mu_j^{(r)}}}{q-q^{-1}}\wt(\bm{\la/\mu};q^{-2},\bm{u},r)\Big|_{\substack{
m=1\\ \bm{u}=1}}=q^{\mu_j-1}(1-q^{-2})^{d-1}(-q^{-2})^{ht(\la/\mu)}.$$

As a specialization of Corollary \ref{t:m-n2}, the following is the Murnaghan--Nakayama rule for $\mathscr{H}_{n}(q)$ (see also \cite{JL1, Pfe1, Ram}).
\begin{cor}
Let $\la,\mu$ be two partitions of $n$, $1\leq j\leq l(\mu)$, and $\psi^{\la}$ be the irreducible character of $\mathscr{H}_{n}(q)$ associated with $\la$. Then
\begin{align}\label{e:m-nH}
\psi^{\la}_{\mu}=\sum_{\nu}q^{\mu_j-1}(1-q^{-2})^{d-1}(-q^{-2})^{ht(\la/\nu)}\psi^{\nu}_{\mu\setminus\mu_j}
\end{align}
summed over all partitions $\nu$ such that $\la/\nu$ is a $\mu_j$-generalized ribbon with $d$ ribbons as its connected components.
\end{cor}
\begin{rem}
We remark that \eqref{e:m-nH} is slightly different from those formulae in the literature due to the different presentation of $\mathscr{H}_{n}(q)$. One can obtain the previous one by noting that $\tilde{\psi}^{\la}_{\mu}(q^2)=q^{|\mu|-l(\mu)}\psi^{\la}_{\mu}(q)$ under the substitution $g_i\rightarrow q^{-1}T_{i-1}$ $(2\leq i\leq n)$, $u_1\rightarrow q^{-1}u_1$, followed by $q^2\rightarrow q$, where $\tilde{\psi}^{\la}_{\mu}$ is the previous irreducible character associated with $\la$ acting on the standard element $T_{\sigma_{\mu}}$ (see \cite[p. 29]{JL1} for the notation $T_{\sigma_{\mu}}$).
\end{rem}

\subsection{Iwahori--Hecke algebra in type $B$}
The Iwahori--Hecke algebra $\mathscr{B}_n(u, q^2)$ in type $B$ is the associative algebra over the field $\mathbb{C}(u,q)$ with generators
$g_1,g_2,\cdots,g_n$ and relations
\begin{align}
g_ig_j &= g_jg_i, \quad|i - j| > 1,\\
g_ig_{i+1}g_i &= g_{i+1}g_ig_{i+1},\quad 2 \leq i \leq n-1,\\
g_1g_2g_1g_2 &= g_2g_1g_2g_1,\\
g_1^2 &= (u - 1)g_1 + u,\\
g^2_i &= (q - q^{-1})g_i + 1, \quad2 \leq i \leq n.
\end{align}
It is clear that $\mathscr{H}_{m,n}(q,\bm{u})$ will reduce to $\mathscr{B}_{n}(u,q^2)$ when $m=2$ and $\bm{u}=(u_1,u_2)=(-1,u)$. We call $\bm{\la}\in\mathcal{P}_{n,2}$ a {\it double partition} of $n$. A $k$-generalized $2$-multi-ribbon is also called a {\em $k$-generalized double ribbon}. Suppose $\bm{\la},\bm{\mu}$ are two double partitions and $\bm{\mu}\subset\bm{\la}$, then, by definition, we have
\begin{align}
\wt(\bm{\la}/\bm{\mu};q,\bm{u},r)\Big|_{\substack{
m=2\\ \bm{u}=(-1,u)}}=
\begin{cases}
(-1)^r(1-q)^{cc(\la^{(1)}/\mu^{(1)})}(-q)^{ht(\la^{(1)}/\mu^{(1)})}, & \text{if $\mu^{(2)}=\la^{(2)}$},\\
u^r(1-q)^{cc(\bm{\la/\mu})}(-q)^{ht(\bm{\la/\mu})}, & \text{if $\mu^{(2)}\subsetneqq\la^{(2)}$}.
\end{cases}
\end{align}

One can easily obtain the Murnaghan--Nakayama rule for $\mathscr{B}_n(u,q^2)$ by specialization of Corollary \ref{t:m-n2}, which is stated below.
\begin{cor}
Let $\bm{\la}, \bm{\mu}$ be two double partitions of $n$ and $\mu^{(r)}_j$ ($r=1,2$) be an arbitrary nonzero part of $\bm{\mu}$. Denote by $\varphi^{\bm{\la}}$ the irreducible character of $\mathscr{B}_n(u,q^2)$ associated with $\bm{\la}$. Then for
fixed $r$
\begin{align}
\varphi^{\bm{\la}}_{\bm{\mu}}=\sum_{\bm{\nu}}\frac{q^{\mu^{(r)}_j}}{(q-q^{-1})}\epsilon^r(1-q^{-2})^{cc(\bm{\la/\nu})}(-q^{-2})^{ht(\bm{\la/\nu})}
\varphi^{\bm{\nu}}_{\bm{\mu}\setminus\mu^{(r)}_j}
\end{align}
summed over all double partitions $\bm{\nu}$ such that $\bm{\la/\nu}$ is a $\mu^{(r)}_j$-generalized double ribbon. Where $\epsilon=\begin{cases}
-1, & \text{if $|\nu^{(2)}|=|\la^{(2)}|$,}\\
u, & \text{if $|\nu^{(2)}|<|\la^{(2)}|$.}
\end{cases}$
\end{cor}
\begin{rem}
Pfeiffer's recursion in \cite[Sec.~12]{Pfe2} computes irreducible character values on the \emph{minimal-length} representatives of conjugacy classes in the type~$B$ Weyl group, and its coefficients in the hook-case may depend on additional \emph{content}-type data of the embedding. In contrast, our specialization is obtained from the Frobenius--type formula for the Ariki--Koike algebra and yields a uniform recurrence on our standard elements, with coefficients determined solely by the generalized (double) ribbon statistics $cc(\,\cdot\,)$ and $ht(\,\cdot\,)$ (together with the parameter choice $\epsilon$).
\end{rem}

\section{A dual Murnaghan--Nakayama rule via vertex operators}\label{s:OV}
In the previous section, we establish a Murnaghan--Nakayama rule, which can be regarded as an iterative formula on the lower multipartitions of the irreducible characters. In this section, as duality, an iteration of the characters on the upper multipartitions will be derived with the help of the vertex operator realization of Schur functions.

\subsection{Vertex operator realization of Schur functions}
Let $\mathfrak{h}$ be the infinite dimensional Heisenberg algebra over $\mathbb{C}$ with generators $h_n$, $n\in\mathbb{Z}\setminus\{0\}$ and the central element $c$ subject to the relation:
\begin{align}\label{e:hrel}
[h_m, h_n]=\delta_{m,-n}m\cdot c.
\end{align}

As is well known, there is a basic representation realized on the space $V=Sym(\mathfrak{h}^-)$ for the algebra $\mathfrak{h}$, where $V$ denotes the symmetric algebra generated by the elements $h_{-n}$, $n\in\mathbb{N}$. Explicitly, for any $n\in\mathbb{N}$, the element $h_{-n}$ acts as the multiplication by $h_{-n}$, while $h_n$ acts as the differentiation operator $n\frac{\partial}{\partial h_{-n}}$, which then satisfy the relation \eqref{e:hrel} with $c=1$.

Define the canonical inner product on $V$ by
\begin{align}\label{e:inner1}
\langle h_{-\la}, h_{-\mu} \rangle_{V}=\delta_{\la\mu}z_{\la},
\end{align}
where $h_{-\la}:=h_{-\la_1}h_{-\la_2}\cdots$ and $z_{\la}:=\prod_{i\geq 1}i^{m_i(\la)}m_i(\lambda)!$.

Therefore, $h_{-n}$ and $h_{n}$ are dual to each other with respect to \eqref{e:inner1}. Namely,
\begin{align}\label{e:dual1}
\langle h_{-n}f, g \rangle_{V}=\langle f, h_ng \rangle_{V} \quad\text{for any $f,g\in V$}.
\end{align}

Define the linear maps $X_n$ and $X^*_n$: $V\longrightarrow V[[z, z^{-1}]]$ by
\begin{align}\label{e:Schurop1}
X(z)&=\mbox{exp} \left( \sum\limits_{n\geq 1} \dfrac{1}{n}h_{-n}z^{n} \right) \mbox{exp} \left( -\sum \limits_{n\geq 1} \dfrac{1}{n}h_{n}z^{-n} \right)=\sum_{n\in\mathbb{Z}}X_nz^n,\\ \label{e:Schurop*1}
X^*(z)&=\mbox{exp} \left(-\sum\limits_{n\geq 1} \dfrac{1}{n}h_{-n}z^{n} \right) \mbox{exp} \left(\sum \limits_{n\geq 1} \dfrac{1}{n}h_{n}z^{-n} \right)=\sum_{n\in\mathbb{Z}}X^*_nz^{-n}.
\end{align}
Recall that $\Lambda_{\mathbb{Q}}$ denotes the space of symmetric functions. The {\it characteristic mapping} $\iota$ from $V$ to $\Lambda_{\mathbb{Q}}$ is defined by
\begin{align}
\iota(h_{-\la}):=p_{\la_1}\cdots p_{\la_l}=p_{\la} \quad \text{for any partition $\la$}.
\end{align}
By definition, $\iota(h_n)=\iota(n\frac{\partial}{\partial h_{-n}})=n\frac{\partial}{\partial p_n}:=p^*_n$. It is shown \cite[Theorem 3.6]{J2} that $\iota$ is an isometric isomorphism between $V$ and $\Lambda_{\mathbb{Q}}$ with respect to the inner product on $\Lambda_{\mathbb{Q}}$:
\begin{align}
\langle p_{\la}, p_{\mu} \rangle=\delta_{\la\mu}z_{\la}.
\end{align}
It follows from \eqref{e:dual1} that
\begin{align}
\langle p_{n}f, g \rangle=\langle f, p_n^*g \rangle \quad\text{for any $f,g\in \Lambda_{\mathbb{Q}}$}.
\end{align}
Note that $*$ is $\mathbb Q(t)$-linear and  anti-involutive.

Under the characteristic mapping, the images of $X(z)$ and $X^*(z)$ (denoted by $S(z)$ and $S^*(z)$ respectively) are given by
\begin{align}\label{e:Schurop2}
\iota(X(z))=S(z)&=\mbox{exp} \left( \sum\limits_{n\geq 1} \dfrac{1}{n}p_{n}z^{n} \right) \mbox{exp} \left( -\sum \limits_{n\geq 1} \dfrac{\partial}{\partial p_n}z^{-n} \right)=\sum_{n\in\mathbb{Z}}S_nz^n,\\ \label{e:Schurop*2}
\iota(X^*(z))=S^*(z)&=\mbox{exp} \left(-\sum\limits_{n\geq 1} \dfrac{1}{n}p_{n}z^{n} \right) \mbox{exp} \left(\sum \limits_{n\geq 1} \dfrac{\partial}{\partial p_n}z^{-n} \right)=\sum_{n\in\mathbb{Z}}S^*_nz^{-n}.
\end{align}

\begin{prop}\cite{J1} Let us recall the vertex operator realization of Schur functions and the relations between $S_{n}$ and $S^*_{m}$.
\begin{enumerate}
\item  For any composition $\mu=(\mu_{1},\ldots,\mu_{k})$, the
product $S_{\mu}.1:=S_{\mu_{1}}\cdots S_{\mu_{k}}.1=s_{\mu}$ is the
Schur function labeled by $\mu$. In general, $s_{\mu}=0$ or $\pm s_{\lambda}$ for a partition $\lambda$ such that $\lambda\in \mathfrak{S}_{l}(\mu+\delta)-\delta.$ Here $\delta=(l-1,l-2,\ldots,1,0),$ where $l = l(\mu)$.
Moreover, $S_{-n}.1=\delta_{n,0}, S^{*}_{n}.1=\delta_{n,0}, (n\geq0)$.

\item  The components of $S(z)$ and $S^{*}(z)$ obey the following commutation relations:
\begin{align}
S_{m}S_{n}+S_{n-1}S_{m+1}&=0,\\
S^{*}_{m}S^{*}_{n}+S^{*}_{n+1}S^{*}_{m-1}&=0,\\
S_{m}S^{*}_{n}+S^{*}_{n-1}S_{m-1}&=\delta_{m,n}.
\end{align}
\end{enumerate}
\end{prop}

Recall that $\{s_{\bla}\mid \bla\in\mathcal{P}_{n,m}\}$ forms an orthonormal basis of $\Lambda_{\mathbb{Q}}^{\otimes m}$. Therefore
\begin{align}\label{e:cha-inner}
\begin{split}
\chi^{\bm{\la}}_{\bm{\mu}}=&\langle q_{\bm{\mu}}(x;q,\bm{u}), s_{\bm{\la}} \rangle_{\otimes} \quad\text{(by \eqref{e:fro} and \eqref{e:orth})}\\
=&\langle q_{\bm{\mu}}(x;q,\bm{u}), s_{\la^{(1)}}\cdots S_{\la^{(j)}}.1\cdots s_{\la^{(m)}} \rangle_{\otimes}\quad\text{(by vertex operator realization)}\\
=&\langle S^*_{\la^{(j)}_1}q_{\bm{\mu}}(x;q,\bm{u}), s_{\la^{(1)}}\cdots s_{\la^{(j)}\setminus\la^{(j)}_1}\cdots s_{\la^{(m)}} \rangle_{\otimes}\\
=&\langle S^*_{\la^{(j)}_1}q_{\bm{\mu}}(x;q,\bm{u}), s_{\bm{\la}\setminus\la^{(j)}_1} \rangle_{\otimes}.
\end{split}
\end{align}
\subsection{An iterative formula on the upper multipartitions}
Next, we will find a formula for $S^{*}_{k}q_{\bm{\mu}}(x;q,\bm{u})$ in terms of $q_{\bm{\rho}}(x;q,\bm{u})S^{*}_{m}$, which further yields an iteration of $\chi^{\bm{\la}}_{\bm{\mu}}$.

\begin{lem}\label{l:sq}
\cite[Prop. 2.3]{JL1}. Let $k\in\mathbb{Z}_{+}$ and $h\in\mathbb{Z}_{\geq 0}$, then
\begin{align}
S^*_{k}q_{h}(x;q^{-2})=q_h(x;q^{-2})S^*_k+(1-q^{-2})\sum_{i=1}^h q_{h-i}(x;q^{-2})S^*_{k-i}.
\end{align}
\end{lem}
The following is the multi-partition generalization of Lemma \ref{l:sq}.
\begin{thm}\label{t:Sq}
Let $k\in\mathbb{Z}_{+}$ and $\bm{\mu}\in \mathcal{P}_{n,m}$, then we have
\begin{align}\label{e:Sq}
S^*_{k}(x^{(1)})q_{\bm{\mu}}(x;q,\bm{u})=\sum_{\bm{\tau}}(1-q^{-2})^{l(\bm{\tau})}q^{|\bm{\tau}|}u_{1}^{\mathcal{L}_{\bm{\mu}\bm{\tau}}}
q_{\bm{\mu}-\bm{\tau}}(x;q,\bm{u})S^*_{k-|\bm{\tau}|}(x^{(1)}).
\end{align}
summed over all $m$-multicompositions $\bm{\tau}$ such that $\bm{\tau}\subset\bm{\mu}$. Here $\mathcal{L}_{\bm{\mu}\bm{\tau}}$ is the non-negative integer associated with multipartition $\bm{\mu}$ and multicomposition $\bm{\tau}$ defined by
\begin{align*}
\mathcal{L}_{\bm{\mu}\bm{\tau}}:=\sum_{i=1}^m i\left(l(\mu^{(i)})-l(\mu^{(i)}-\tau^{(i)})\right).
\end{align*}
\end{thm}
\begin{proof}
Note that $S^*_{k}(x^{(1)})$ is commutative with $q_{t}(x^{(i)};q^{-2})$ $(i>1)$.
It follows that
\begin{align*}
&S^*_{k}(x^{(1)})q_{n}^{(i)}(x;q,\bm{u})\\
=&S^*_{k}(x^{(1)})\left(\frac{q^n}{q-q^{-1}}\sum_{c\in C_{n,m}}u_{c}^iq_{c_1}(x^{(1)};q^{-2})q_{c_2}(x^{(2)};q^{-2})\cdots q_{c_m}(x^{(m)};q^{-2})\right)\qquad\text{(by definition)}\\
=&\frac{q^n}{q-q^{-1}}\sum_{c\in C_{n,m}}u_{c}^i \left(q_{c_1}(x^{(1)};q^{-2})S^*_{k}(x^{(1)})+\sum_{a=1}^{c_1}(1-q^{-2})q_{c_1-a}(x^{(1)};q^{-2})S^*_{k-a}(x^{(1)})\right)q_{c_2}(x^{(2)};q^{-2})\\
&\hspace{6.8em}\cdots q_{c_m}(x^{(m)};q^{-2})\qquad\text{(by Lemma \ref{l:sq})}\\
=&q^{(i)}_n(x;q,\bm{u})S^*_{k}(x^{(1)})+\sum_{a=1}^{n-1}(1-q^{-2})q^aq^{(i)}_{n-a}(x;q,\bm{u})S^*_{k-a}(x^{(1)})+q^{n-1}u_1^iS^*_{k-n}(x^{(1)}).
\end{align*}
The last identity holds by swapping the summations.
Recall that $q^{(i)}_0(x;q,\bm{u})=\frac{1}{q-q^{-1}}$. A direct iteration gives that
\begin{align*}
S^*_{k}(x^{(1)})q^{(i)}_{\mu}(x;q,\bm{u})=\sum_{\tau}(1-q^{-2})^{l(\tau)}q^{|\tau|}u_1^{i(l(\mu)-l(\mu-\tau))}
q^{(i)}_{\mu-\tau}(x;q,\bm{u})S^*_{k-|\tau|}(x^{(1)})
\end{align*}
summed over all compositions $\tau\subset\mu$.

By a further iteration, we have
\begin{align*}
S^*_{k}(x^{(1)})q_{\mu}(x;q,\bm{u})=\sum_{\bm{\tau}\subset\bm{\mu}}(1-q^{-2})^{l(\bm{\tau})}q^{|\bm{\tau}|}u_1^{\sum_{i}i(l(\mu^{(i)})-l(\mu^{(i)}-\tau^{(i)}))}
q_{\bm{\mu}-\bm{\tau}}(x;q,\bm{u})S^*_{k-|\bm{\tau}|}(x^{(1)}),
\end{align*}
completing the proof.
\end{proof}

We remark that if we let the RHS of \eqref{e:Sq} act on the vacuum vector $1$, then by \eqref{e:Schurop*2} we will find $S^*_{k-|\bm{\tau}|}(x^{(1)}).1$ is exactly the elementary symmetric function (up to a sign). Explicitly, 
\begin{align}\label{e:S.1}
    S^*_{k-|\bm{\tau}|}(x^{(1)}).1=
    \begin{cases}
        (-1)^{\bm{|\tau|}-k}e_{\bm{|\tau|}-k}(x^{(1)}) & \text{if $k\leq |\bm{\tau}|$},\\
        0, & \text{otherwise}.
    \end{cases}
\end{align}

\begin{lem}[Multi-partition version of Pieri rule]\label{l:MPS}
 For $k\in \mathbb{N}$, $\bla\in\mathcal{P}_{n,m}$, then
 \begin{align}
     e_k(x^{(1)})s_{\bla}=\sum_{\brho}s_{\brho}
 \end{align}
 summed over all multipartitions $\brho\in\mathcal{P}_{n-k,m}$ such that $\rho^{(i)}=\la^{(i)}$, $i> 1$ and $\la^{(1)}/\rho^{(1)}$ is a $k$-vertical strip.
\end{lem}
\begin{proof}
    This can be derived by the classical Pieri rule for Schur functions, which can be found in \cite[p. 331]{Mac}.
\end{proof}

\begin{thm}\label{t:daulMN}
    Let $(\bm{\mu-\tau})^*$ denote the multipartition obtained by rearranging each composition in the multicomposition $\bm{\mu-\tau}$. We have the following iterative formula for $\chi^{\bm{\la}}_{\bm{\mu}}$.
    \begin{align}\label{e:daulMN}
      \chi^{\bm{\la}}_{\bm{\mu}}=\sum_{\bm{\tau,\rho}}(-1)^{|\bm{\tau}|-\la^{(1)}_1}(1-q^{-2})^{l(\bm{\tau})}q^{|\bm{\tau}|}u_1^{\mathcal{L}_{\bm{\mu\tau}}}\chi^{\bm{\rho}}_{(\bm{\mu-\tau})^*}  
    \end{align}
    summed over all multicompositions $\bm{\tau}\subset\bm{\mu}$ and multipartitions $\bm{\rho}$ such that $\rho^{(k)}=\la^{(k)}$, $k> 1$ and $(\la^{(1)}\setminus\la^{(1)}_1)/\rho^{(1)}$ is a $(|\bm{\tau}|-\la^{(1)}_1)$-vertical strip.
\end{thm}
\begin{proof}
    We have 
    \begin{align*}
        \chi^{\bm{\la}}_{\bm{\mu}}=&\langle S^*_{\la^{(1)}_1}q_{\bm{\mu}}(x;q,\bm{u}), s_{\bm{\la}\setminus\la^{(1)}_1} \rangle_{\otimes}\quad\text{(by \eqref{e:cha-inner})}\\
        =&\sum_{\bm{\tau}}(1-q^{-2})^{l(\bm{\tau})}q^{|\bm{\tau}|}u_{1}^{\mathcal{L}_{\bm{\mu}\bm{\tau}}}
\langle q_{\bm{\mu}-\bm{\tau}}(x;q,\bm{u})S^*_{\la^{(1)}_1-|\bm{\tau}|}(x^{(1)}).1, s_{\bm{\la}\setminus\la^{(1)}_1} \rangle_{\otimes}\quad\text{(by \eqref{e:Sq})}\\
=&\sum_{\bm{\tau}}(1-q^{-2})^{l(\bm{\tau})}q^{|\bm{\tau}|}u_{1}^{\mathcal{L}_{\bm{\mu}\bm{\tau}}}
\langle q_{\bm{\mu}-\bm{\tau}}(x;q,\bm{u})(-1)^{\bm{|\tau|}-\la^{(1)}_1}e_{\bm{|\tau|}-\la^{(1)}_1}(x^{(1)}), s_{\bm{\la}\setminus\la^{(1)}_1} \rangle_{\otimes}\quad\text{(by \eqref{e:S.1})}\\
=&\sum_{\bm{\tau}}(-1)^{\bm{|\tau|}-\la^{(1)}_1}(1-q^{-2})^{l(\bm{\tau})}q^{|\bm{\tau}|}u_{1}^{\mathcal{L}_{\bm{\mu}\bm{\tau}}}
\langle q_{\bm{\mu}-\bm{\tau}}(x;q,\bm{u}), e^*_{\bm{|\tau|}-\la^{(1)}_1}(x^{(1)})s_{\bm{\la}\setminus\la^{(1)}_1} \rangle_{\otimes}\\
=&\sum_{\bm{\tau}}(-1)^{\bm{|\tau|}-\la^{(1)}_1}(1-q^{-2})^{l(\bm{\tau})}q^{|\bm{\tau}|}u_{1}^{\mathcal{L}_{\bm{\mu}\bm{\tau}}}
\left\langle q_{\bm{\mu}-\bm{\tau}}(x;q,\bm{u}), \sum_{\bm{\rho}}s_{\bm{\rho}} \right\rangle_{\otimes}\quad\text{(by Lemma \ref{l:MPS})}
    \end{align*}
    summed over all multipartitions $\bm{\rho}$ such that $\rho^{(k)}=\la^{(k)}$, $k\neq 1$ and $(\la^{(1)}\setminus\la^{(1)}_1)/\rho^{(1)}$ is a $(|\bm{\tau}|-\la^{(1)}_1)$-vertical strip.
\end{proof}

\begin{cor}\label{c:dualMN2}
    Let $1\leq j\leq m$, we have
    \begin{align}\label{e:dualMN2}
\chi^{(\varnothing,\cdots,\varnothing,\la^{(j)},\cdots,\la^{(m)})}_{\bm{\mu}}=\sum_{\bm{\tau,\rho}}(-1)^{|\bm{\tau}|-\la^{(j)}_1}(1-q^{-2})^{l(\bm{\tau})}q^{|\bm{\tau}|}u_j^{\mathcal{L}_{\bm{\mu\tau}}}\chi^{\bm{\rho}}_{(\bm{\mu-\tau})^*}  
    \end{align}
    summed over all multicompositions $\bm{\tau}\subset\bm{\mu}$ and multipartitions $\bm{\rho}$ such that
    \begin{itemize}
        \item $\rho^{(i)}=\varnothing$ for $i=1,\cdots,j-1$;
        \item $(\la^{(j)}\setminus\la^{(j)}_1)/\rho^{(j)}$ is a $(|\bm{\tau}|-\la^{(j)}_1)$-vertical strip;
        \item $\rho^{(i)}=\la^{(i)}$ for $i=j+1,\cdots,m$.
    \end{itemize}
\end{cor}
\begin{proof}
    This can be derived as follows:
    \begin{align*}    \chi^{(\varnothing,\cdots,\varnothing,\la^{(j)},\cdots,\la^{(m)})}_{\bm{\mu}}=&\chi^{(\la^{(j)},\varnothing,\cdots,\varnothing,\la^{(j+1)},\cdots,\la^{(m)})}_{\bm{\mu}}\big|_{u_1\rightarrow u_j} \quad\text{(by Proposition \ref{p:replace})} \\
      =&\sum_{\bm{\tau,\tilde{\rho}}}(-1)^{|\bm{\tau}|-\la^{(j)}_1}(1-q^{-2})^{l(\bm{\tau})}q^{|\bm{\tau}|}u_1^{\mathcal{L}_{\bm{\mu\tau}}}\chi^{\bm{\tilde{\rho}}}_{(\bm{\mu-\tau})^*}\big|_{u_1\rightarrow u_j} \quad\text{(by Theorem \ref{t:daulMN})}
    \end{align*}
    summed over all multicompositions $\bm{\tau}\subset\bm{\mu}$ and multipartitions $\bm{\tilde{\rho}}$ such that
    \begin{itemize}
        \item $(\la^{(j)}\setminus\la^{(j)}_1)/\tilde{\rho}^{(1)}$ is a $(|\bm{\tau}|-\la^{(j)}_1)$-vertical strip;
        \item $\tilde{\rho}^{(i)}=\varnothing$ for $i=2,\cdots,j$ and $\tilde{\rho}^{(i)}=\la^{(i)}$ for $i=j+1,\cdots,m$.
    \end{itemize}
    Using Proposition \ref{p:replace} again completes the proof.
\end{proof}

\begin{rem}
    Corollary \ref{c:dualMN2} ensures that the iteration can continue even when the parts of $\la^{(1)}$ are exhausted. The Murnaghan--Nakayama rule expresses $\chi^{\bm{\la}}_{\bm{\mu}}$ in terms of irreducible characters where one part is removed from the lower multipartition $\bm{\mu}$. In contrast, Corollary \ref{c:dualMN2} computes $\chi^{\bm{\la}}_{\bm{\mu}}$ using irreducible characters in which at least one part is removed from the upper multipartition $\bm{\la}$. For this reason, we refer to equation \eqref{e:dualMN2} as the dual Murnaghan--Nakayama rule for $\chi^{\bm{\la}}_{\bm{\mu}}$.
\end{rem}

\begin{rem}
It is worth pointing out that there is a slight difference in how we apply formulas \eqref{e:m-n2} and \eqref{e:dualMN2}. When using formula \eqref{e:m-n2}, we remove one part from a partition of $\bmu$ at each step. In contrast, when using formula \eqref{e:dualMN2}, we subtract some parts from all partitions of $\bmu$ simultaneously at each step. As a result, the final lower multipartition obtained through formula \eqref{e:m-n2} is the empty multipartition, while the final lower multipartition obtained through formula \eqref{e:dualMN2} is a multipartition where every part in each partition is zero. The corresponding irreducible character values for these two cases are different. For the former, the value is $1$, while for the latter, it is $\left(\frac{1}{q - q^{-1}}\right)^{l(\bmu)}$ for the reason that $q^{(i)}_0(x;q,\bm{u})=\frac{1}{q-q^{-1}}$. 
Here, as an example, we use \eqref{e:m-n2} and \eqref{e:dualMN2} respectively to compute $\chi^{((n),\varnothing,\cdots,\varnothing)}_{\bmu}$. 
\begin{enumerate}
    \item use \eqref{e:m-n2}: 
    \begin{align*}
\chi^{((n),\varnothing,\cdots,\varnothing)}_{\bmu}&=\frac{q^n}{(q-q^{-1})^{l(\bmu)}}u_1^{\mathcal{L}_{\bmu\bmu}}(1-q^{-2})^{l(\bmu)}(-q^{-2})^0\chi^{((0),\varnothing,\cdots,\varnothing)}_{\bm{\varnothing}} \\
&=q^{n-l(\bmu)}u_1^{\mathcal{L}_{\bmu\bmu}}
    \end{align*}
    \item use \eqref{e:dualMN2}: 
    \begin{align*}
\chi^{((n),\varnothing,\cdots,\varnothing)}_{\bmu}&=(1-q^{-2})^{l(\bmu)}q^n u_1^{\mathcal{L}_{\bmu\bmu}}\chi^{((0),\varnothing,\cdots,\varnothing)}_{((0,\cdots,0),\cdots,(0,\cdots,0))} \\
&=(1-q^{-2})^{l(\bmu)}q^n u_1^{\mathcal{L}_{\bmu\bmu}}\frac{1}{(q-q^{-1})^{l(\bmu)}}\\
&=q^{n-l(\bmu)}u_1^{\mathcal{L}_{\bmu\bmu}}.
    \end{align*}
\end{enumerate}
\end{rem}

\section{Application I: Regev-type formula}\label{s:Regev}
In this section, as the first application of our Murnaghan--Nakayama rule, we are devoted to giving the Regev formula for the cyclotomic Hecke algebra $\mathscr{H}_{m,n}(q,\bm{u})$. Once the Murnaghan--Nakayama rule is set up right,
 the proof of the Regev formula, which is based on the classical deep work on symmetric groups and Hecke algebras in type $A$, is remarkably natural. To proceed with it, we initially, in accordance with \cite{Zhao1}, review the permutation super representation of the cyclotomic Hecke algebra $\mathscr{H}_{m,n}(q,\bm{u})$, as well as the Schur–Sergeev duality that exists between a certain quantum superalgebra $\mathscr{U}_q(\bm{g})$ and $\mathscr{H}_{m,n}(q,\bm{u})$.

We refer to a superspace $W$ as a $\mathbb{Z}_{2}$-graded vector space with decomposition $W=W_{\Bar{0}}\oplus W_{\Bar{1}}$. An element $w\in W$ is homogeneous with degree $\Bar{i}$ if $w\in W_{\Bar{i}}$, $i=0,1$. We also define the dimension of $W$ by $\rm{dim}(W):={\rm dim}(W_{\Bar{0}})| {\rm dim}(W_{\Bar{1}})$.

Let $V$ represent a direct sum of $m$ superspaces, denoted as $V=\bigoplus_{i=1}^m V^{(i)}$. In this context, we assume ${\rm dim}(V^{(i)})=k_i| l_i$ and $\{v^{(i)}_1,v^{(i)}_2,\cdots, v^{(i)}_{k_i}, v^{(i)}_{k_i+1},\cdots, v^{(i)}_{k_i+l_i}\}$ forms a homogeneous basis of $V^{(i)}$. Evidently, ${\rm dim}(V)=k| l$ with $k=\sum k_i$ and $l=\sum l_i$. Additionally, $\mathfrak{B}:=\{v_1,v_2,\cdots,v_{k+l}\}$ constitutes a homogeneous basis for $V$ with $v_d=v^{(s)}_t$ if $d=\sum_{i=1}^{s}(k_i+l_i)+t$. The entity $c(v_i)=r$ is defined contingent upon the condition $v_i\in V^{(r)}$. The term $c(\cdot)$ is occasionally referred to as the {\em color function}.

Introduce the linear map $\omega:~V\longrightarrow V$ defined by $\omega(v_i)=u_{c(v_i)}v_i$. Define the linear operators $\pi, \sigma:~ V\otimes V\longrightarrow V\otimes V$ by 
\begin{align}
    \pi(v_i\otimes v_j)&=
    \begin{cases}
        (q-q^{-1})v_i\otimes v_j+(-1)^{\Bar{i}\Bar{j}}v_j\otimes v_i, & \text{if $i<j$};\\
        \frac{(q-q^{-1})+(-1)^{\Bar{i}}(q+q^{-1})}{2}v_i\otimes v_i, & \text{if $i=j$};\\
        (-1)^{\Bar{i}\Bar{j}}v_j\otimes v_i, & \text{if $i>j$}.
        \end{cases}\\
        \sigma(v_i\otimes v_j)&=
    \begin{cases}
        \pi(v_i\otimes v_j), & \text{if $c(v_i)=c(v_j)$};\\
        (-1)^{\Bar{i}\Bar{j}}v_j\otimes v_i, & \text{if $c(v_i)\neq c(v_j)$}.
    \end{cases}
\end{align}
We can demonstrate $\pi^2=(q-q^{-1})\pi+1$, which implies $\pi$ is invertible and $\pi^{-1}=\pi-(q-q^{-1}){\rm Id}$.
We canonically extend $\pi$ and $\sigma$ to linear maps $\pi_i,\sigma_i (i=1,2,\cdots,n-1):~V^{\otimes n}\longrightarrow V^{\otimes n}$ by
\begin{align}
    \pi^{\pm 1}_{i}&:={\rm Id}^{\otimes (i-1)}\otimes \pi^{\pm 1}\otimes {\rm Id}^{\otimes(n-i-1)}\\
    \sigma_{i}&:={\rm Id}^{\otimes (i-1)}\otimes \sigma\otimes {\rm Id}^{\otimes(n-i-1)}\\
    \omega_{i}&:={\rm Id}^{\otimes (i-1)}\otimes \omega\otimes {\rm Id}^{\otimes(n-i)}.
\end{align}
Furthermore, we define $\pi_0:=\pi_{1}^{-1}\pi_2^{-1}\cdots\pi_{n-1}^{-1}\sigma_{n-1}\cdots\sigma_{1}\omega_1$. As is shown in \cite[Theorem 3.6]{Zhao1}, the map
\begin{align*}
    \Xi^{\bm{u},q}_{\bm{k}\mid\bm{l};n}: ~\mathscr{H}_{m,n}(q,\bm{u})&\longrightarrow {\rm End}_{\mathbb{K}}V^{\otimes n}\\
    g_{i}~ &\longmapsto \pi_{i-1}, \quad i=1,2,\cdots,n.
\end{align*}
gives a (super) representation of $\mathscr{H}_{m,n}(q,\bm{u})$, called {\em permutation super representation} of $\mathscr{H}_{m,n}(q,\bm{u})$. We denote by $\chi^{\bm{u},q}_{\bm{k}\mid\bm{l};n}(g(\bm{\mu}))$ the value of the character $\chi^{\bm{u},q}_{\bm{k}\mid\bm{l};n}$ of $\Xi^{\bm{u},q}_{\bm{k}\mid\bm{l};n}$ on the standard element $g(\bm{\mu})$.

Let $\bm{k}=(k_1,k_2,\cdots,k_m)$ and $\bm{l}=(l_1,l_2,\cdots,l_m)$ be two compositions. For a given nonnegative integer $n$, we define the set $\mathfrak{C}(n;\bm{k}\mid\bm{l})$ of pairs of $m$-multicompositions associated with $\bm{k}$ and $\bm{l}$ by
\begin{align*}
 \mathfrak{C}(n;\bm{k}|\bm{l}):=\left\{(\bm{\alpha},\bm{\beta})\left| \, 
\begin{array}{l}
 \bm{\alpha}=(\alpha^{(1)},\cdots,\alpha^{(m)}) \text{ with $l(\alpha^{(i)})\leq k_i$ } \\
\bm{\beta}=(\beta^{(1)},\cdots,\beta^{(m)})\text{ with $l(\beta^{(i)})\leq l_i$ }
\end{array}
\right.\text{and $|\bm{\alpha}|+|\bm{\beta}|=n$} \right\}.   
\end{align*}
The main result in this section is the following formula for $\chi^{\bm{u},q}_{\bm{k}\mid\bm{l};n}(g(\bm{\mu}))$, which is referred to as Regev's formula.
\begin{thm}\label{t:Regev}
Let $\bm{\mu}$ be a $m$-multipartition of $n$, then the value of the character of the permutation super representation on the standard element $g(\bm{\mu})$ is given by
\begin{align}\label{e:Regev}
\begin{split}
    &\chi^{\bm{u},q}_{\bm{k}|\bm{l};n}(g(\bm{\mu}))\\
=&q^{n-l(\bm{\mu})}\prod_{r=1}^m\prod_{j=1}^{l(\mu^{(r)})}\sum_{(\bm{\alpha};\bm{\beta})\in \mathfrak{C}(\mu^{(r)}_j;\bm{k}\mid\bm{l})} \bm{u}^r_{\overrightarrow{l}(\bm{\alpha};\bm{\beta})}(1-q^{-2})^{l(\balp;\bbeta)-1}(-q^{-2})^{|\bbeta|-l(\bbeta)}\prod_{i=1}^m\binom{k_i}{l(\alpha^{(i)})}\binom{l_i}{l(\beta^{(i)})}
\end{split}
\end{align}
where $\overrightarrow{l}(\bm{\alpha};\bm{\beta}):=\max\{1\leq i\leq m\mid l(\alpha^{(i)})+l(\beta^{(i)})>0\}$ and $l(\balp;\bbeta):=l(\balp)+l(\bbeta)$.
\end{thm}

\begin{rem}
    In a recent development, Zhao \cite{Zhao2} has proposed an alternative method to derive Theorem \ref{t:Regev}. It should be noted that his formulation is slightly different from ours, due to a distinct choice of presentation of $\mathscr{H}_{m,n}(q,\bm{u})$, as well as the standard elements $g(\bm{\mu})$.
\end{rem}

Next we will introduce a series of notations. Let $\bm{\la}$ be a $m$-multipartition with diagram $\D(\bm{\la})$. A {\em $(\bm{k},\bm{l})$-semistandard tableau} $\bm{\T}_{\bm{\la}}$ of shape $\bm{\la}$ is a $m$-tuple of $(k_i,l_i)$-semistandard tableaux (cf. \cite{BR}) obtained by filling each component $\D(\bm{\la}^{(i)})$ of $\D(\bm{\la})$ with variables $\{a^{(i)}_1<a^{(i)}_2<\cdots<a^{(i)}_{k_i}<b^{(i)}_1<b^{(i)}_2<\cdots<b^{(i)}_{l_i}, i=1,2,\cdots,m\}$, satisfying
\begin{enumerate}
    \item the cells filled with $a^{(i)}_{j}$'s form a tableau;
    \item $a^{(i)}_{j}$'s are nondecreasing in rows, strictly increasing in columns;
    \item $b^{(i)}_{j}$'s are nondecreasing in columns, strictly increasing in rows.
\end{enumerate}
Denote $s_{\bm{k},\bm{l}}(\bm{\la}):=\#\{(\bm{k},\bm{l})\text{-semistandard tableaux of shape $\la$}\}$. For a given $(\bm{\alpha},\bm{\beta})\in \mathfrak{C}(n;\bm{k}|\bm{l})$, we say that a $(\bm{k},\bm{l})$-semistandard tableau $\bm{\T}_{\bm{\la}}$ is of weight $(\bm{\alpha},\bm{\beta})$ if the number of $a^{(i)}_j$ (resp. $b^{(i)}_j$) in $\bm{\T}_{\bm{\la}}$ equals $\alpha^{(i)}_j$ (resp. $\beta^{(i)}_j$). Furthermore, denote by $s_{(\bm{\alpha},\bm{\beta})}(\bm{\la})$ the number of all $(\bm{k},\bm{l})$-semistandard tableaux of shape $\bm{\la}$ with weight $(\bm{\alpha},\bm{\beta})$. Clearly we have
\begin{align}\label{e:semi-tab}
    s_{\bm{k},\bm{l}}(\bm{\la})=\sum_{(\bm{\alpha},\bm{\beta})\in\mathfrak{C}(n;\bm{k}|\bm{l})}s_{(\bm{\alpha},\bm{\beta})}(\bm{\la})
\end{align} 
for any $m$-multipartition $\bm{\la}$ of $n$.

Define $H(k| l;n):=\{\la\vdash n\mid \la_{k+1}\leq l\}$. An element in $H(k| l;n)$ is called a {\em $(k,l)$-hook partition}. More generally, a multipartition is called a {\em $(\bm{k},\bm{l})$-hook multipartition} if each $\la^{(i)}$ is a $(k_i,l_i)$-hook partition for $i=1,2,\cdots,m$. Denote by $H(\bm{k}|\bm{l};n)$ the set of all $(\bm{k},\bm{l})$-hook multipartitions of $n$. Clearly, $H((1^m)| (1^m);n)$ consists of all multipartitions $\bm{\la}$ of $n$ such that each component of $\bm{\la}$ is a hook partition. Note that the assertion in \cite[Definition 2.4]{BR} implies that $s_{\bm{k},\bm{l}}(\bm{\la})\neq 0$ if and only if $\bm{\la}\in H(\bm{k}| \bm{l};n)$ for a multipartition $\bm{\la}$ of $n$.

Now we will express $\chi^{\bm{u},q}_{\bm{k}\mid\bm{l};n}$ as a linear combination of irreducible characters $\chi^{\bm{\la}}$ using the above notations. To continue with this development, the Schur–Sergeev duality between $\mathscr{H}_{m,n}(q,\bm{u})$ and the quantum superalgebra $\mathscr{U}_{q}(\bm{g})$ is needed. Where $\mathscr{U}_{q}(\bm{g}):=\mathscr{U}_{q}(\bm{gl}(k_1| l_1))\otimes\cdots\otimes \mathscr{U}_{q}(\bm{gl}(k_m| l_m))$, among which $\mathscr{U}_{q}(\bm{gl}(k_i| l_i))$, introduced in \cite{BKK}, is the quantum enveloping superalgebra of $\bm{gl}(k_i| l_i)$. $V$ can be a representation of $\mathscr{U}_{q}(\bm{g})$, referred to as {\em vector presentation}. The Hopf structure of $\mathscr{U}_{q}(\bm{g})$ induces an action on $V^{\otimes n}$, denoted by $\Theta_n$. As is demonstrated in \cite{BKK}, the irreducible summands of the representation $V^{\otimes n}$ of $\mathscr{U}_q(gl(k,l))$ are indexed by the $(k,l)$-hook partitions of $n$, with each irreducible summand possessing the dimension $s_{k,l}(\la)$. Consequently, the irreducible summands of the representation $(\Theta_n, V^{\otimes n})$ of $\mathscr{U}_q(\bm{g})$ are characterized by the $(\bm{k},\bm{l})$-hook partitions of $n$, with each irreducible summand affording the dimension $s_{\bm{k},\bm{l}}(\bm{\la})$.

It is confirmed in \cite[Theorem 4.9]{Zhao1} that $(\Xi^{\bm{u},q}_{\bm{k}\mid\bm{l};n}, V^{\otimes n})$ and $(\Theta_n, V^{\otimes n})$ satisfy the Schur–Sergeev duality, i.e.,
\begin{align}
    \Xi^{\bm{u},q}_{\bm{k}\mid\bm{l};n}(\mathscr{H}_{m,n}(q,\bm{u}))={\rm End}_{\mathscr{U}_{q}(\bm{g})}V^{\otimes n}, \quad \Theta_n(\mathscr{U}_{q}(\bm{g}))={\rm End}_{\mathscr{H}_{m,n}(q,\bm{u})}V^{\otimes n}
\end{align}
Moreover, as a $(\mathscr{U}_q(\bm{g}), \mathscr{H}_{m,n}(q,\bm{u})$-bimodule,
\begin{align}
    V^{\otimes n}\cong\bigoplus_{\bm{\la}\in H(\bm{k}\mid\bm{l};n)}W_{\bm{\la}}\otimes W^{\bm{\la}},
\end{align}
where $W_{\bm{\la}}$ (resp. $W^{\bm{\la}}$) is the irreducible module of $\mathscr{U}_q(\bm{g})$ (resp. $\mathscr{H}_{m,n}(q,\bm{u})$) indexed by $\bm{\la}$. This gives the following decomposition:
\begin{align}
    V^{\otimes n}\cong\bigoplus_{\bm{\la}\in H(\bm{k}\mid\bm{l};n)}s_{\bm{k},\bm{l}}(\bm{\la})W^{\bm{\la}},
\end{align}
which is equivalent to
\begin{align}\label{e:expression}
\chi^{\bm{u},q}_{\bm{k}\mid\bm{l};n}=\sum_{\bm{\la}\in H(\bm{k}\mid\bm{l};n)} s_{\bm{k},\bm{l}}(\bm{\la})\chi^{\bm{\la}}.
\end{align}

Now we introduce the {\em skew characters} of $\mathscr{H}_{m,n}(q,\bm{u})$ via the Jacobi--Trudi rule.
Let $\bmu\subset\bla$ be two $m$-multipartitions (with $|\bmu|\le |\bla|$). For each $1\le r\le m$ and $a\in\mathbb Z$, set
\[
h_a^{(r)} :=
\begin{cases}
\chi^{\bvar_r^a}, & a\ge 1,\\
1, & a=0,\\
0, & a<0,
\end{cases}
\]
where $1$ denotes the unit character (for the degree $0$ algebra). We define
\begin{align}\label{e:def-skew-character}
\chi^{\bla/\bmu}:=\widehat{\bigotimes_{r=1}^{m}}
\det_{1\le i,j\le l(\la^{(r)})}\big(h^{(r)}_{\la^{(r)}_i-i-\mu^{(r)}_j+j}\big),
\end{align}
where the determinant is taken in the character ring with multiplication given by the outer product.


Next, we will restate our Murnaghan--Nakayama rule in terms of the skew characters. Before that, we need some terminologies about skew diagrams. Let $\bla/\bmu$ be a $m$-multi skew diagram, denote by ${\rm SE}(\bla/\bmu)$ the subset of $\bla/\bmu$ consisting of the southeast ribbon which starts from the northeastmost box and ends at the southeastmost box in each connected component in $\D(\bla/\bmu)$.
\begin{figure}
\begin{align*}
\left( \begin{array}{ccccc}
\begin{tikzpicture}[scale = 0.5]
      \begin{scope}
        \clip (0,0) -| (3,2) -| (4,3) -| (6,4) -- (2,4) |- (1,3) |- (0,1) -- (0,0);
        \draw [color=black!25] (0,0) grid (6,4);
      \end{scope}
      \draw [thick] (0,0) -| (3,2) -| (4,3) -| (6,4) -- (2,4) |- (1,3) |- (0,1) -- (0,0);
      \filldraw[fill = green]
    (0,0) rectangle (1,1) (1,0) rectangle (2,1) (2,0)rectangle (3,1) (2,1)rectangle(3,2) (2,2)rectangle(3,3) (3,2)rectangle(4,3) (3,3)rectangle(4,4) (4,3)rectangle(5,4) (5,3)rectangle(6,4);
      \node at (2.5,0.5) {$\bullet$};
      \node at (2.5,1.5) {$\bullet$};
    \end{tikzpicture},~&
    \begin{tikzpicture}[scale = 0.5]
      \begin{scope}
        \clip (0,0) -| (3,3) -| (6,4) -- (2,4) |- (1,2) |- (0,1) -- (0,0);
        \draw [color=black!25] (0,0) grid (6,4);
      \end{scope}
      \draw [thick] (0,0) -| (3,3) -| (6,4) -- (2,4) |- (1,2) |- (0,1) -- (0,0);
      \filldraw[fill = green]
    (0,0) rectangle (1,1) (1,0) rectangle (2,1) (2,0)rectangle (3,1) (2,1)rectangle(3,2) (2,2)rectangle(3,3) (2,3)rectangle(3,4) (3,3)rectangle(4,4) (4,3)rectangle(5,4) (5,3)rectangle(6,4);
      \node at (2.5,0.5) {$\bullet$};
      \node at (4.5,3.5) {$\bullet$};
      \node at (5.5,3.5) {$\bullet$};
    \end{tikzpicture},~&
    \begin{tikzpicture}[scale = 0.5]
      \begin{scope}
        \clip (0,0) -| (2,2) -| (5,3) -| (9,4) -- (3,4) |- (2,3) |- (0,2) -- (0,0);
        \draw [color=black!25] (0,0) grid (9,4);
      \end{scope}
      \draw [thick] (0,0) -| (2,2) -| (5,3) -| (9,4) -- (3,4) |- (2,3) |- (0,2) -- (0,0);
      \filldraw[fill = green]
    (0,0) rectangle (1,1) (1,0) rectangle (2,1) (1,1)rectangle (2,2) (2,2)rectangle(3,3) (3,2)rectangle(4,3) (4,2)rectangle(5,3) (4,3)rectangle(5,4) (5,3)rectangle(6,4) (6,3)rectangle(7,4) (7,3)rectangle(8,4) (8,3)rectangle(9,4);
      \node at (1.5,0.5) {$\bullet$};
      \node at (4.5,2.5) {$\bullet$};
      \node at (7.5,3.5) {$\bullet$};
      \node at (8.5,3.5) {$\bullet$};
    \end{tikzpicture}
\end{array}  \right)
\end{align*}
\caption{${\rm SE}(\bla/\bmu)$ consists of the green boxes. The subset formed by the boxes with bullets is an element of $\R_9(\bla/\bmu)$.}\label{F:SE}
    \end{figure}
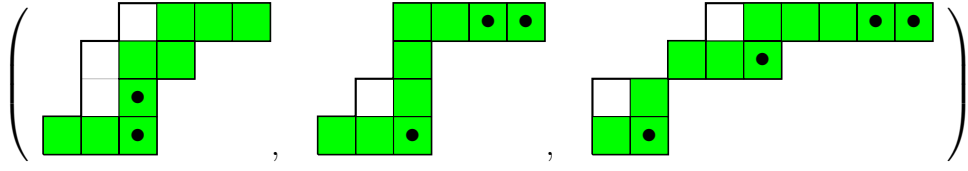
We identify a box in the $a$-th row and the $b$-th column of $\D(\bla^{(i)}/\bmu^{(i)})$ with the coordinate $(a^{(i)},b^{(i)})$, $i=1,2,\cdots,m$. Let $\R_{t}(\bla/\bmu)$ be the set of all subsets $R$ of ${\rm SE}(\bla/\bmu)$ satisfying
\begin{itemize}
    \item $R$ has size $t$;
    \item For any $(a^{(i)},b^{(i)})\in R$ and $(x^{(i)},y^{(i)})\in {\rm SE}(\bla/\bmu)\setminus R$, we have either $x^{(i)}<a^{(i)}$ or $y^{(i)}<b^{(i)}$.
\end{itemize}
Figure \ref{F:SE} presents an example for ${\rm SE}(\bla/\bmu)$ and $\R_{t}(\bla/\bmu)$. It is not hard to see that for any $m$-multi skew diagram $\D(\bla/\bmu)$ and $R\in\R_t(\bla/\bmu)$, $\D(\bla/\bmu)\setminus R$ is still a $m$-multi skew diagram. For $R\in\R_{t}(\bla/\bmu)$, denote by $cc(R)$ the number of connected components of $R$. If $R$ has $d$ connected components $R_1,R_2,\cdots,R_d$, define $ht(R):=\sum_{i=1}^d ht(R_i)$.


Let $\bm\nu\in\mathcal P_{|\bla|-|\bmu|,m}$, and fix $1\le r\le m$, $1\le j\le l(\nu^{(r)})$. Then we have the following skew-diagram version of the Murnaghan--Nakayama rule:
\begin{align}\label{e:skew-m-n}
    \chi^{\bla/\bmu}_{\bm{\nu}}
    =q^{\nu^{(r)}_j-1}\sum_{R\in\R_{\nu^{(r)}_j}(\bla/\bmu)}
    u^{r}_{\overrightarrow{R}}(1-q^{-2})^{cc(R)-1}(-q^{-2})^{ht(R)}
    \chi^{(\bla/\bmu)\setminus R}_{\bm{\nu}\setminus \nu^{(r)}_j}.
\end{align}
Here
\[
\overrightarrow{R}:=\max\{1\leq i\leq m\mid \text{there exists at least one box of $R$ in }\bla^{(i)}/\bmu^{(i)}\}.
\]
Equation \eqref{e:skew-m-n} can be deduced from Corollary \ref{t:m-n2} using arguments similar to those in \cite[Chap.~2.4]{JK} (see also \cite[Theorem~5.2]{Zhao3}). More precisely, one applies Corollary \ref{t:m-n2} to the terms arising from the Jacobi--Trudi determinant expansion defining $\chi^{\bla/\bmu}$, and then regroups the resulting terms according to the subsets $R\in \R_{\nu^{(r)}_j}(\bla/\bmu)$ of ${\rm SE}(\bla/\bmu)$. The regrouped terms are precisely the skew characters $\chi^{(\bla/\bmu)\setminus R}$.

For $(\balp,\bbeta)\in\mathfrak{C}(n;\bm{k}|\bm{l})$, denote by $R_{\balp,\bbeta}$ some (any) $m$-multi skew diagrams $\D(\bla/\bmu)$ (This is not unique!) such that the connected components of $D(\la^{(i)}/\mu^{(i)})$ are horizontal strips $H^{(i)}_1,H^{(i)}_2,\cdots,H^{(i)}_{k_i}$ and vertical strips $V^{(i)}_1,V^{(i)}_2,\cdots,V^{(i)}_{l_i}$ $(i=1,2,\cdots,m)$ satisfying 
\begin{align}
   (H^{(i)}_1,H^{(i)}_2,\cdots,H^{(i)}_{k_i})=(\alpha^{(i)}_1,\alpha^{(i)}_2,\cdots,\alpha^{(i)}_{k_i}),\quad (V^{(i)}_1,V^{(i)}_2,\cdots,V^{(i)}_{l_i})=(\beta^{(i)}_1,\beta^{(i)}_2,\cdots,\beta^{(i)}_{l_i}). 
\end{align}

\begin{figure}
\begin{align*}
\left( \begin{array}{ccccc}
 \begin{tikzpicture}[scale = 0.4]
      \begin{scope}
        \clip (0,0) -| (1,3) -| (4,4) -| (5,6) -| (7,7) -- (0,7) -- (0,0);
        \draw [color=black!25] (0,0) grid (7,7);
      \end{scope}
      \draw [thick] (0,0) -| (1,3) -| (4,4) -| (5,6) -| (7,7) -- (0,7) -- (0,0);
      \filldraw[fill = green]
    (0,0) rectangle (1,1) (0,1) rectangle (1,2) (1,3)rectangle (2,4) (2,3)rectangle(3,4) (3,3)rectangle(4,4) (4,4)rectangle(5,5) (4,5)rectangle(5,6) (5,6)rectangle(6,7) (6,6)rectangle(7,7);
    \end{tikzpicture},\qquad&
    \begin{tikzpicture}[scale = 0.4]
      \begin{scope}
        \clip (0,0) -| (2,1) -| (3,4) -| (5,6) -| (8,7) -- (0,7) -- (0,0);
        \draw [color=black!25] (0,0) grid (8,7);
      \end{scope}
      \draw [thick] (0,0) -| (2,1) -| (3,4) -| (5,6) -| (8,7) -- (0,7) -- (0,0);
      \filldraw[fill = green]
    (0,0) rectangle (1,1) (1,0) rectangle (2,1) (2,1)rectangle (3,2) (2,2)rectangle(3,3) (2,3)rectangle(3,4) (3,4)rectangle(4,5) (4,4)rectangle(5,5) (5,6)rectangle(6,7) (6,6)rectangle(7,7) (7,6)rectangle(8,7);
    \end{tikzpicture}
\end{array}  \right)
\end{align*}
\caption{An example for $R_{\balp,\bbeta}$ where $\balp=((2,3), (3,2,2))$ and $\bbeta=((2,2),(3))$. The corresponding skew character can be written as $\chi^{(\balp,\bbeta)}=\chi^{(2)}\hat{\otimes}\chi^{(3)}\hat{\otimes}\chi^{(1^2)}\hat{\otimes}\chi^{(1^2)}\hat{\otimes}\chi^{(3)}\hat{\otimes}\chi^{(2)}\hat{\otimes}\chi^{(2)}\hat{\otimes}\chi^{(1^3)}.$}\label{F:HV}
    \end{figure}
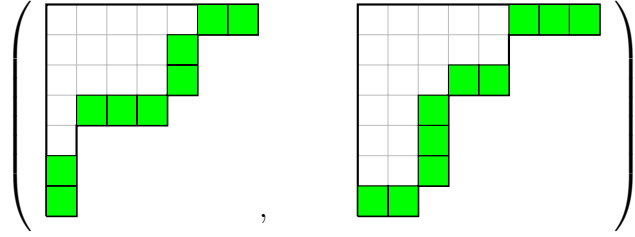

For $R_{\balp,\bbeta}$, define the character related to $R_{\balp,\bbeta}$ by 
\begin{align}
    \chi^{(\balp;\bbeta)}:=\widehat{\bigotimes\limits_{1\leq i\leq m}}\left(\chi^{(\alpha^{(i)}_1)}\hat{\otimes}\cdots\hat{\otimes}\chi^{(\alpha^{(i)}_{k_i})}\hat{\otimes}\chi^{(1^{\beta^{(i)}_{1}})}\hat{\otimes}\cdots\hat{\otimes}\chi^{(1^{\beta^{(i)}_{l_i}})}\right).
\end{align}
See Figure \ref{F:HV} for an example for $R_{\balp,\bbeta}$. 

Due to \cite[Lemma 3.23]{BR} (see also \cite[Lemma 3.2]{Tay}, \cite[Corollary 5.5]{Zhao3}), we have the following decomposition for $\chi^{(\balp;\bbeta)}$ in terms of irreducible characters.
\begin{align}\label{e:alp-beta}
    \chi^{(\balp;\bbeta)}=\sum_{\bla\in H(\bm{k}|\bm{l};n)}s_{(\balp,\bbeta)}(\bla)\chi^{\bla}.
\end{align}

Now we restate our skew version of the Murnaghan--Nakayama rule in a special case. $R_{\balp,\bbeta}$ is a $m$-multi skew diagram as mentioned above. Then 
\begin{align}
    \R_t(R_{\balp,\bbeta})=\{R_{\bgam,\bdel}\mid (\bgam,\bdel)\subset(\balp,\bbeta), |\bgam|+|\bdel|=t\}.
\end{align}
Furthermore, if $R_{\bgam,\bdel}\in \R_t(R_{\balp,\bbeta})$, then we have $cc\left(R_{\bgam,\bdel}\right)=l(\bgam)+l(\bdel)=l(\bgam;\bdel)$, $ht(R_{\bgam,\bdel})=|\bdel|-l(\bdel)$ and $\overrightarrow{R_{\bgam,\bdel}}=\overrightarrow{l}(\bgam;\bdel)$. We can restate \eqref{e:skew-m-n} in the following special setting.
\begin{align}\label{e:skew-m-n2}
    \chi^{(\balp;\bbeta)}_{\bm{\nu}}=q^{\nu^{(r)}_j-1}\sum_{\substack{
(\bgam,\bdel)\subset(\balp,\bbeta)\\ |\bgam|+|\bdel|=\nu^{(r)}_j}} u^{r}_{\overrightarrow{l}(\bgam,\bdel)}(1-q^{-2})^{l(\bgam;\bdel)-1}(-q^{-2})^{|\bdel|-l(\bdel)}\chi^{(\balp;\bbeta)\setminus(\bgam;\bdel)}_{\bm{\nu}\setminus\nu^{(r)}_j}.
\end{align}

We are ready to prove Theorem \ref{t:Regev}.\\
{\em Proof of Theorem \ref{t:Regev}.} 
\begin{align*}
    &\chi^{\bm{u},q}_{\bm{k}|\bm{l};n}(g(\bm{\mu}))\\
    =&\sum_{\bm{\la}\in H(\bm{k}\mid\bm{l};n)} s_{\bm{k},\bm{l}}(\bm{\la})\chi^{\bm{\la}} \quad \text{(by \eqref{e:expression})}\\
    =&\sum_{\bm{\la}\in H(\bm{k}\mid\bm{l};n)} 
    \sum_{(\bm{\alpha},\bm{\beta})\in\mathfrak{C}(n;\bm{k}|\bm{l})}s_{(\bm{\alpha},\bm{\beta})}(\bm{\la})\chi^{\bm{\la}} \quad \text{(by \eqref{e:semi-tab})}\\
    =&\sum_{(\bm{\alpha},\bm{\beta})\in\mathfrak{C}(n;\bm{k}|\bm{l})} \chi^{(\balp;\bbeta)} \quad \text{(by \eqref{e:alp-beta})}\\
    =&\sum_{(\bm{\alpha},\bm{\beta})\in\mathfrak{C}(n;\bm{k}|\bm{l})}q^{\mu^{(r)}_j-1}\sum_{\substack{
(\bgam,\bdel)\subset(\balp,\bbeta)\\ |\bgam|+|\bdel|=\mu^{(r)}_j}} u^{r}_{\overrightarrow{l}(\bgam,\bdel)}(1-q^{-2})^{l(\bgam;\bdel)-1}(-q^{-2})^{|\bdel|-l(\bdel)}\chi^{(\balp;\bbeta)\setminus(\bgam;\bdel)}_{\bm{\mu}\setminus\mu^{(r)}_j}\quad \text{(by \eqref{e:skew-m-n2})}.
\end{align*}
For a given $(\bm{\alpha}^{'},\bm{\beta}^{'})\in\mathfrak{C}(n-\mu^{(r)}_j;\bm{k}|\bm{l})$, if $(\bm{\alpha},\bm{\beta})\setminus(\bgam,\bdel)=(\bm{\alpha}^{'},\bm{\beta}^{'})$, then the number of $(\bgam,\bdel)$'s contributing the coefficient
\begin{align*}
    q^{\mu^{(r)}_j-1}u^{r}_{\overrightarrow{l}(\bgam,\bdel)}(1-q^{-2})^{l(\bgam;\bdel)-1}(-q^{-2})^{|\bdel|-l(\bdel)}
\end{align*}
is exactly
\begin{align*}
    \prod_{i=1}^m\binom{k_i}{l(\gamma^{(i)})}\binom{l_i}{l(\delta^{(i)})}.
\end{align*}
Therefore,
\begin{align*}
    &\chi^{\bm{u},q}_{\bm{k}|\bm{l};n}(g(\bm{\mu}))\\
    =&\left(\sum_{(\bm{\gamma},\bm{\delta})\in\mathfrak{C}(\mu^{(r)}_j;\bm{k}|\bm{l})}q^{\mu^{(r)}_j-1} u^{r}_{\overrightarrow{l}(\bgam,\bdel)}(1-q^{-2})^{l(\bgam;\bdel)-1}(-q^{-2})^{|\bdel|-l(\bdel)}\prod_{i=1}^m\binom{k_i}{l(\gamma^{(i)})}\binom{l_i}{l(\delta^{(i)})}\right)\\
    &\times \sum_{(\bm{\alpha}^{'},\bm{\beta}^{'})\in\mathfrak{C}(n-\mu^{(r)}_j;\bm{k}|\bm{l})}\chi^{(\balp^{'};\bbeta^{'})}_{\bm{\mu}\setminus\mu^{(r)}_j}.
\end{align*}
An induction argument completes the proof. $\hspace{15em}\Box$

Next, we consider the special case $(\bm{k}|\bm{l})=((1^m)|(1^m))$. It follows from the proof of Theorem \ref{t:Regev} that
\begin{align}\label{e:equality}
    \chi^{\bm{u},q}_{(1^m)|(1^m);n}=\sum_{(\bm{\alpha},\bm{\beta})\in\mathfrak{C}(n;(1^m)|(1^m))} \chi^{(\balp;\bbeta)}. 
\end{align}
Note that every $(\balp,\bbeta)$ in $\mathfrak{C}(n;(1^m)|(1^m))$ has the form
$$(\balp,\bbeta)=((\alpha_1),(\alpha_2),\cdots,(\alpha_m),(\beta_1),(\beta_2),\cdots,(\beta_m))$$
with $\sum_{i=1}^m(\alpha_i+\beta_i)=n$. So
\begin{align}\label{e:q2}
    \text{RHS of \eqref{e:equality} }=\sum_{\substack{
c_1+c_2+\cdots+c_m=n\\ c_i\geq 0}}\widehat{\bigotimes_{1\leq i\leq m}}\sum_{\alpha_i=0}^{c_i}\chi^{(\alpha_i;c_i-\alpha_i)}.
\end{align}
By \cite[Remark 3.3]{Tay}, 
\begin{align}
    \sum_{\alpha_i=0}^{c_i}\chi^{(\alpha_i;c_i-\alpha_i)}=
    \begin{cases}
        1 &\text{if $c_i=0$};\\
        2\sum_{\alpha_i=0}^{c_i-1}\chi^{(c_i-\alpha_i,1^{\alpha_i})} &\text{if $c_i\geq 1$}.
    \end{cases}
\end{align}
Substituting this into \eqref{e:q2} gives that
\begin{align}\label{e:q3}
 \text{RHS of \eqref{e:equality} }=\sum_{\bla\in H((1^m)| (1^m);n)}   2^{L(\bla)}\chi^{\bla}.
\end{align}
Recall that $H((1^m)| (1^m);n)$ consists of all $m$-multipartitions $\bm{\la}$ of $n$ such that each component $\la^{(i)}$ is a hook partition. Here $L(\bla):=\#\{1\leq i\leq m\mid \la^{(i)}\neq \varnothing\}$.

On the other hand, by \eqref{e:Regev},
\begin{align}\label{e:q4}
 \chi^{\bm{u},q}_{(1^m)|(1^m);n}(g(\bm{\mu}))
=q^{n-l(\bm{\mu})}\prod_{r=1}^m\prod_{j=1}^{l(\mu^{(r)})}\sum_{(\alpha;\beta)} \bm{u}^r_{\overrightarrow{l}(\alpha;\beta)}(1-q^{-2})^{l(\alpha)+l(\beta)-1}(-q^{-2})^{|\beta|-l(\beta)}
\end{align}
the inner sum runs over all pairs of $m$-compositions $(\alpha;\beta)$ with $|\alpha|+|\beta|=\mu^{(r)}_j$.

Note that
\begin{align}
\begin{split}\label{e:q5}
 &\sum_{\substack{(\alpha;\beta)\\
 |\alpha|+|\beta|=\mu^{(r)}_j}} \bm{u}^r_{\overrightarrow{l}(\alpha;\beta)}(1-q^{-2})^{l(\alpha)+l(\beta)-1}(-q^{-2})^{|\beta|-l(\beta)}\\
 =&\frac{1}{1-q^{-2}}\sum_{\substack{
c_1+c_2+\cdots+c_m=\mu^{(r)}_j\\ c_i\geq 0}}\bm{u}^r_{\overrightarrow{c}}\sum_{\substack{(\alpha;\beta)\\
 |\alpha|+|\beta|=c_i}}(1-q^{-2})^{\sum_{i=1}^m(2-\delta_{\beta_i,c_i}-\delta_{\beta_i,0})}(-q^{-2})^{\sum_{i=1}^m(\beta_i+\delta_{\beta_i,0}-1)}\\
 =&\frac{1}{1-q^{-2}}\sum_{\substack{
c_1+c_2+\cdots+c_m=\mu^{(r)}_j\\ c_i\geq 0}}\bm{u}^r_{\overrightarrow{c}}\sum_{\substack{(\alpha;\beta)\\
 |\alpha|+|\beta|=c_i}}\prod_{i=1}^m(1-q^{-2})^{2-\delta_{\beta_i,c_i}-\delta_{\beta_i,0}}(-q^{-2})^{\beta_i+\delta_{\beta_i,0}-1}\\
 =&\frac{1}{1-q^{-2}}\sum_{\substack{
c_1+c_2+\cdots+c_m=\mu^{(r)}_j\\ c_i\geq 0}}\bm{u}^r_{\overrightarrow{c}}\prod_{i=1}^m\sum_{\beta_i=0}^{c_i}(1-q^{-2})^{2-\delta_{\beta_i,c_i}-\delta_{\beta_i,0}}(-q^{-2})^{\beta_i+\delta_{\beta_i,0}-1}\\
=&\sum_{\substack{
c_1+c_2+\cdots+c_m=\mu^{(r)}_j\\ c_i\geq 0}}\bm{u}^r_{\overrightarrow{c}}2^{\#c}(1-q^{-2})^{\#c-1}\prod_{i=1}^m[c_i]_{-q^{-2}},
\end{split}
\end{align}
where $\overrightarrow{c}:=\max\{1\leq i\leq m\mid c_i\neq 0\}$, $c:=\{1\leq i\leq m\mid c_i\neq 0\}$ and $[k]_q:=\frac{1-q^k}{1-q}$ (Here $k\geq 0$ and $[0]_q:=1$) is the standard $q$-integer. The last equation holds for the following identity:
\begin{align*}
   \sum_{\beta_i=0}^{c_i}(1-q^{-2})^{2-\delta_{\beta_i,c_i}-\delta_{\beta_i,0}}(-q^{-2})^{\beta_i+\delta_{\beta_i,0}-1}=
   \begin{cases}
       2(1-q^{-2})[c_i]_{-q^{-2}} & \text{if $c_i>0$}\\
       1 & \text{if $c_i=0$}.
   \end{cases}
\end{align*}

We say that a $m$-tuple of matrices $(M^{(1)},M^{(2)},\cdots,M^{(m)})$ is a {\em $m$-multimatrix} if each $M^{(i)}$ has $m$ columns with entries nonnegative integers.  Denote by $\mathcal{M}_{m}$ the set of all $m$-multimatrices. For a $m$-multimatrix $\bm{M}=(M^{(1)},M^{(2)},\cdots,M^{(m)})$, define
\begin{align*}
    \bm{r}(\bm{M}):=\sum_{i=1}^m r(M^{(i)}),\quad   \bm{u}_{\bm{M}}:=\prod_{i=1}^{m}\prod_{j=1}^{r(M^{(i)})}u^{i}_{\overrightarrow{r}_j(M^{(i)})}.
\end{align*}
Here we use the notations  $r(M^{(i)})$ and $\overrightarrow{r}_j(M^{(i)})$ to present the number of rows in matrix $M^{(i)}$ and the $j$-th row vector of matrix $M^{(i)}$, respectively.
Thus $u_{\overrightarrow{r}_j(M^{(i)})}=u_{f_{ij}}$ with $f_{ij}$ the column containing the first non-zero element from right to left in the $j$th row of $M^{(i)}$. For instance, if
\begin{align}\label{e:example}
\bm{M}=\left( \begin{array}{ccc}
\begin{pmatrix}3 &  \textcolor{orange}{2} & 0 \\
 0 & 1 & \textcolor{orange}{6}
\end{pmatrix},~&
    \begin{pmatrix}\textcolor{orange}{2} & 0 & 0 \\
 0 & 4 & \textcolor{orange}{2}\\
 \textcolor{orange}{8}& 0 & 0
\end{pmatrix},~&
\begin{pmatrix}2, \textcolor{orange}{2}, 0 
\end{pmatrix} 
\end{array}  \right).
\end{align}
We have marked the last non-zero element in each row in \textcolor{orange}{orange}. Thus,
\begin{align*}
(f_{ij})=\left( \begin{array}{ccc}
\begin{pmatrix}  2  \\
 3
\end{pmatrix},~&
    \begin{pmatrix}1 \\
 3\\
 1
\end{pmatrix},~&
 \begin{pmatrix}2
\end{pmatrix}
\end{array}  \right),\quad \bm{u}_{\bm{M}}=u^4_1u^4_2u_3^3.
\end{align*}

Substituting \eqref{e:q5} into \eqref{e:q4} gives that
\begin{align}\label{e:q6}
 \chi^{\bm{u},q}_{(1^m)|(1^l);n}(g(\bm{\mu}))=\sum_{\bm{M}}u_{\bm{M}}2^{\#\bm{M}}(q-q^{-1})^{\#\bm{M}-l(\bmu)}\prod_{m^{(r)}_{ij}\in\bm{M}}[m^{(r)}_{ij}]_{-q^{-2}},
 \end{align}
 summed over all $m$-multimatrices $(M^{(1)},\cdots,M^{(m)})\in\mathcal{M}_m$ such that the $j$-th row sum of $M^{(r)}$ equals $\mu^{(r)}_j$. Where $\#\bm{M}$ is the number of nonzero entries in $\bm{M}$.

 Now we have obtained the following result.
 \begin{cor}\label{t:special-hook}
    Let $\bmu$ be a multipartition of $n$. Then we have the following combinatorial formula, which is also referred as the Regev formula. 
     \begin{align}\label{e:special-hook}
 \sum_{\substack{\bla\vdash n\\ {\rm all}~\la^{(i)}~{\rm hook}}}   2^{L(\bla)}\chi^{\bla}_{\bmu}=\sum_{\bm{M}}u_{\bm{M}}2^{\#\bm{M}}(q-q^{-1})^{\#\bm{M}-l(\bmu)}\prod_{m^{(r)}_{ij}\in\bm{M}}[m^{(r)}_{ij}]_{-q^{-2}},
 \end{align}
 summed over all $m$-multimatrices $(M^{(1)},\cdots,M^{(m)})$ such that the $j$th row sum of $M^{(r)}$ equals $\mu^{(r)}_j$.
     \begin{proof}
         It follows from \eqref{e:equality}, \eqref{e:q3} and \eqref{e:q6} by noticing that $H((1^m)| (1^m);n)$ consists of all multipartitions $\bla$ of $n$ with each component $\la^{(i)}$ being hook shape.
     \end{proof}
 \end{cor}

\begin{rem}
We remark that Zhao \cite{Zhao2} obtained an asymptotic formula for the LHS of \eqref{e:special-hook}.
It is pertinent to point out that an alternative proof for Corollary \ref{t:special-hook} is presented by establishing a connection between the irreducible characters of $\mathscr{H}_{m,n}(q,\bm{u})$ and those of Hecke algebra in type $A$. This proof relies on the Regev formula for Hecke algebra in type $A$. To prevent redundancy, we leave this proof in the appendix \ref{app}.
   \end{rem}
   
 Recall that $\mathscr{H}_{m,n}(q,\bm{u})$ will reduce to Hecke algebra $\mathscr{H}_{n}(q)$ (resp. complex reflection group $\mathscr{W}_{m,n}$) when $m=1$, $u_1=1$ (resp. $q=1$, $\bm{u}=\{1,\zeta,\cdots,\zeta^{m-1}\}$). Here $\zeta$ is the primitive $m$-th root of unity. Before closing this section, we want to examine how the Regev formula adapts in these two special cases. To do so, we will analyze each case individually.  For the first case, i.e., $m=1$, $u_1=1$, it is easy to see \eqref{e:special-hook} is degraded to
 \begin{align}
     \sum_{\la~\text{hook}}\psi^{\la}_{\mu}=2^{l(\mu)-1}\prod_{i=1}^{l(\mu)}[\mu_i]_{-q^{-2}}.
 \end{align}
Where $\psi^{\la}$ is the irreducible character of $\mathscr{H}_{n}(q)$ indexed by $\la$. This is exactly the Regev formula for the Hecke algebra in type $A$ established in \cite{JL1,Zhao3}. For the second case, i.e., $q=1$, $\bm{u}=\{1,\zeta,\cdots,\zeta^{m-1}\}$, we have
\begin{align*}
   &\sum_{\bm{M}}u_{\bm{M}}2^{\#\bm{M}}(q-q^{-1})^{\#\bm{M}-l(\bmu)}\prod_{m^{(r)}_{ij}\in\bm{M}}[m^{(r)}_{ij}]_{-q^{-2}}=
   \begin{cases}
       \sum\limits_{\bm{M}} \bm{u}_{\bm{M}} 2^{l(\bm{\mu})} & \text{if all $\mu^{(r)}_j$'s are odd}\\
       0 & \text{otherwise}
   \end{cases}
\end{align*}
summed over all $m$-multimatrices $\bm{M}$ such that $j$-th row of $M^{(r)}$ has only one nonzore entry $\mu^{(r)}_{j}$. Moreover, if all $\mu^{(r)}_j$'s are odd, then
\begin{align*}
  \sum\limits_{\bm{M}} \bm{u}_{\bm{M}} 2^{l(\bm{\mu})}=\prod_{r=1}^m\prod_{j=1}^{l(\mu^{(r)})}2\sum_{i=1}^{m}\zeta^{(r-1)i}=
  \begin{cases}
      (2m)^{l(\mu^{(1)})}, & \text{if $\bmu=(\mu^{(1)},0,\cdots,0)$};\\
      0, &\text{otherwise}.
  \end{cases}
\end{align*}
The last equation holds for the following identity:
\begin{align*}
    \sum_{i=1}^m\zeta^{(r-1)i}=
    \begin{cases}
        m, &\text{if $r=1$};\\
        0, &\text{if $1<r\leq m$}.
    \end{cases}
\end{align*}

We have reached the Regev formula for $\mathscr{W}_{m,n}$.
\begin{cor}
    Suppose $\bmu$ is a $m$-multipartition of $n$, then
    \begin{align}
        \sum_{\substack{\bla\vdash n\\ \la^{(i)}~{\rm hook}}}\phi^{\bla}_{\bmu}=
        \begin{cases}
            (2m)^{l(\mu^{(1)})}, &\text{if $\bmu=(\mu^{(1)},0,\cdots,0)$ with all parts odd;}\\
            0, &\text{otherwise}.
        \end{cases}
    \end{align}
    Where $\phi^{\bla}$ is the irreducible character of $\mathscr{W}_{m,n}$ indexed by multipartition $\bla$.
\end{cor}

\section{Application II:  L\"ubeck--Prasad--Adin--Roichman-type formula}\label{s:LPAR}
In this section, as the second application of our Murnaghan--Nakayama rule, we will establish a  L\"ubeck--Prasad--Adin--Roichman formula for the cyclotomic Hecke algebra $\mathscr{H}_{m,n}(q,\bm{u})$. To proceed with this, We will present an algorithm for $\chi^{(\la^{(1)},\la^{(2)},\cdots,\la^{(m)})}_{\bmu}$ using the corresponding single partition $\la$ with empty $m$-core and $(\la^{(1)},\la^{(2)},\cdots,\la^{(m)})$ as its $m$-quotient. Finally, we will explain how our formula recovers the Adin-Roichman formula for the complex reflection group $\mathscr{W}_{m,n}$ \cite{AR}, and particularly the L\"ubeck-Prasad formula for the hyperoctahedral group \cite{LP}.

\subsection{A restatement of the Murnaghan--Nakayama rule}
\begin{figure}
	\centering
    \begin{tikzpicture}[>=Stealth, line width=.9pt, every node/.style={inner sep=0pt}]

\node (root) at (0,0)   {$\left(\YDthreeone\qquad\YDtwoonone\qquad\YDtwtwo\right)$};

\node (L)    at (-4,-3) {$\left(\YDtwo\quad\YDtwoonone\quad\YDtwtwo\right)$};
\node (R)    at ( 4,-3) {$\left(\YDonesone\quad\YDtwoonone\quad\YDtwtwo\right)$};

\node (LL)   at (-6,-6) {$\left(\YDtwo\quad\YDfourcol\quad\YDtwtwo\right)$};
\node (LR)   at (-2,-6) {$\left(\YDtwo\quad\YDtweone\quad\YDtwtwo\right)$};

\node (RL)   at ( 2,-6) {$\left(\YDonesone\quad\YDfourcol\quad\YDtwtwo\right)$};
\node (RR)   at ( 6,-6) {$\left(\YDonesone\quad\YDtweone\quad\YDtwtwo\right)$};

\draw[->,shorten <=6pt,shorten >=6pt] (root) -- (L);
\draw[->,shorten <=6pt,shorten >=6pt] (root) -- (R);

\draw[->,shorten <=6pt,shorten >=6pt] (L) -- (LL);
\draw[->,shorten <=6pt,shorten >=6pt] (L) -- (LR);

\draw[->,shorten <=6pt,shorten >=6pt] (R) -- (RL);
\draw[->,shorten <=6pt,shorten >=6pt] (R) -- (RR);

\end{tikzpicture}
	\caption{Four possible removal processes for $\eta^{(3,1),(2,1,1),(2,2)}_{(2,1);r}$. Here $t=1$ and we choose $\sigma=id$. The sum terms they contribute are respectively $u^r_2(1-q^{-2})^3(-q^{-2})^0$, $u^r_2(1-q^{-2})^3(-q^{-2})^0$, $u^r_2(1-q^{-2})^2(-q^{-2})^0$ and $u^r_2(1-q^{-2})^2(-q^{-2})^0$.}\label{F:remove}
\end{figure}
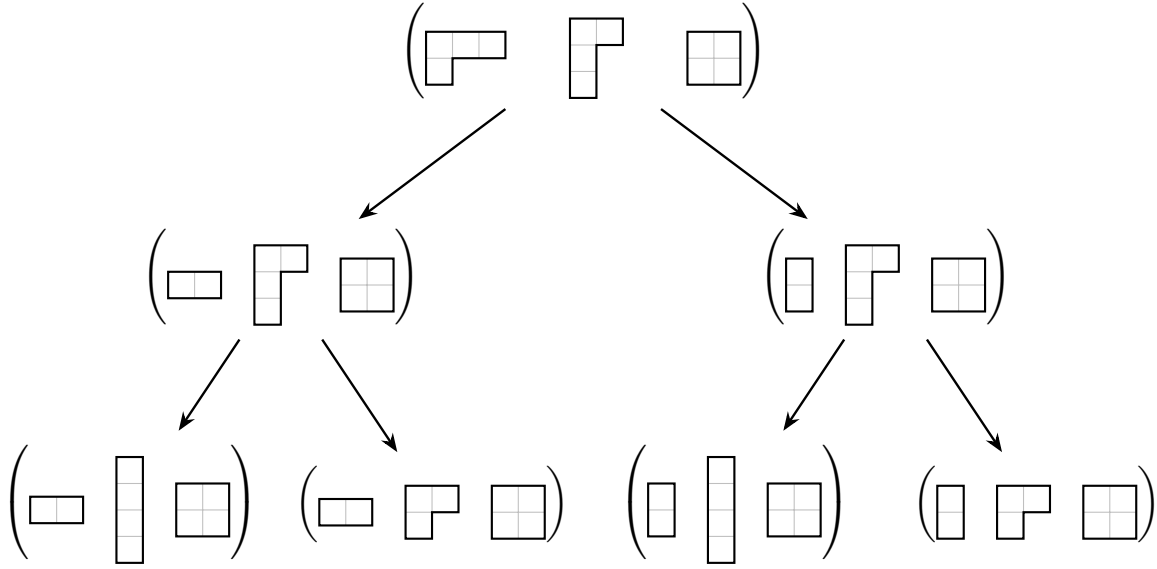

In this subsection, we restate the Murnaghan--Nakayama rule established in Section \ref{s:M-N} in terms of $\eta^{\bla}_{\btau;r}$ defined as below. 
\begin{defn}\label{d:eta}
    Let $\bla=(\la^{(1)},\la^{(2)},\cdots,\la^{(m)})\in \mathcal{P}_{n,m}$, $\btau=(\tau^{(1)},\tau^{(2)},\cdots,\tau^{(t)})$ be a $t$-tuple of partitions with each partition having at most $m$ parts and $r$ be a positive integer. Define $\eta^{\bla}_{\btau;r}$ as the sum of values obtained by all possible applications of the following algorithm:\\
    {\em Initialization}: $\varepsilon:=1$.
    \begin{enumerate}
        \item If $t=0$, end the program and output $\varepsilon$.
        \item If $l(\tau^{(t)})<m$, add $0$'s to the tail to make sure $\tau^{(t)}$ has $m$ components.
        \item Choose a permutation $\sigma\in \mathfrak{S}_m$. Remove a $\tau^{(t)}_{\sigma(j)}$-generalized ribbon $\theta_j$ from $\la^{(j)}$ for all $j=1,2,\cdots,m$. If there is no such $\sigma$ to do the removal process, set $\varepsilon=0$ and end the program.
        \item Multiply $\varepsilon$ by $u^r_{\m}(1-q^{-2})^{-1+\sum_{j=1}^{m}cc(\theta_j)}(-q^{-2})^{\sum_{j=1}^{m}ht(\theta_j)}$, where $cc(\theta_j)$ (resp. $ht(\theta_j)$) is the number of connected components (resp. the height) of $\theta_j$. And $\m:=\max\{1\leq j\leq m\mid \tau^{(t)}_{\sigma(j)}\neq 0\}$.
        \item Redefine $\btau:=(\tau^{(1)},\tau^{(2)},\cdots,\tau^{(t-1)})$ and $t:=t-1$.
        \item Return to step (1). 
    \end{enumerate}
\end{defn}
Figure \ref{F:remove} presents four possible removal processes. 

For a partition $\mu=(\mu_1,\mu_2,\cdots,\mu_t)$ $(\mu_t>0)$, define $\mathcal{P}_m(\mu):=\{(\tau^{(1)},\tau^{(2)},\cdots,\tau^{(t)})\mid \tau^{(j)}\vdash \mu_j ~\text{and}~ l(\tau^{(j)})\leq m, 1\leq j\leq t\}$. Let $\bla,\bmu\in\mathcal{P}_{n,m}$. Then we have the following restatement of the Murnaghan--Nakayama rule:
\begin{align}\label{e:restate-MN}
    \chi^{\bla}_{\bmu}=q^{n-l(\bmu)}\sum_{r=1}^m\sum_{\btau\in\mathcal{P}_m(\mu^{(r)})}\eta^{\bla}_{\btau;r}.
\end{align}

\subsection{A  L\"ubeck--Prasad--Adin--Roichman-type formula}
\begin{figure}
\begin{tikzpicture}[scale=0.5]
    \coordinate (Origin)   at (0,0);
    \coordinate (XAxisMin) at (0,0);
    \coordinate (XAxisMax) at (14,0);
    \coordinate (YAxisMin) at (0,0);
    \coordinate (YAxisMax) at (0,-14);
    \draw [thin, black] (0,0) -- (7,0);
    \draw [thin, black] (0,-1) --(7,-1);
    \draw [thin, black] (0,-2) -- (7,-2);
    \draw [thin, black] (0,-3) -- (7,-3);
    \draw [thin, black] (0,-4) -- (6,-4);
    \draw [thin, black] (0,-5) -- (3,-5);
    \draw [thin, black] (0,-6) -- (3,-6);
    \draw [thin, black] (0,-7) -- (2,-7);
     \draw [thin, black] (0,-8) -- (1,-8);
     \draw [thin, black] (0,0) -- (0,-8);
     \draw [thin, black] (1,0) -- (1,-8);
      \draw [thin, black] (2,0) -- (2,-7);
     \draw [thin, black] (3,0) -- (3,-6);
     \draw [thin, black] (4,0) -- (4,-4);
     \draw [thin, black] (5,0) -- (5,-4);
      \draw [thin, black] (6,0) -- (6,-4);
      \draw [thin, black] (7,0) -- (7,-3);
      \draw[postaction={decorate}, decoration={markings, 
        mark=at position 0.6 with {\arrow{stealth}}}] (0,-8) -- (1,-8);
        \draw[postaction={decorate}, decoration={markings, 
        mark=at position 0.6 with {\arrow{stealth}}}] (1,-8) -- (1,-7);
        \draw[postaction={decorate}, decoration={markings, 
        mark=at position 0.6 with {\arrow{stealth}}}] (1,-7) -- (2,-7);
        \draw[postaction={decorate}, decoration={markings, 
        mark=at position 0.6 with {\arrow{stealth}}}] (2,-7) -- (2,-6);
        \draw[postaction={decorate}, decoration={markings, 
        mark=at position 0.6 with {\arrow{stealth}}}] (2,-6) -- (3,-6);
        \draw[postaction={decorate}, decoration={markings, 
        mark=at position 0.6 with {\arrow{stealth}}}] (3,-6) -- (3,-5);
        \draw[postaction={decorate}, decoration={markings, 
        mark=at position 0.6 with {\arrow{stealth}}}] (3,-5) -- (3,-4);
        \draw[postaction={decorate}, decoration={markings, 
        mark=at position 0.6 with {\arrow{stealth}}}] (3,-4) -- (4,-4);
        \draw[postaction={decorate}, decoration={markings, 
        mark=at position 0.6 with {\arrow{stealth}}}] (4,-4) -- (5,-4);
        \draw[postaction={decorate}, decoration={markings, 
        mark=at position 0.6 with {\arrow{stealth}}}] (5,-4) -- (6,-4);
        \draw[postaction={decorate}, decoration={markings, 
        mark=at position 0.6 with {\arrow{stealth}}}] (6,-4) -- (6,-3);
        \draw[postaction={decorate}, decoration={markings, 
        mark=at position 0.6 with {\arrow{stealth}}}] (6,-3) -- (7,-3);
        \draw[postaction={decorate}, decoration={markings, 
        mark=at position 0.6 with {\arrow{stealth}}}] (7,-3) -- (7,-2);
        \draw[postaction={decorate}, decoration={markings, 
        mark=at position 0.6 with {\arrow{stealth}}}] (7,-2) -- (7,-1);
        \draw[postaction={decorate}, decoration={markings, 
        mark=at position 0.6 with {\arrow{stealth}}}] (7,-1) -- (7,0);
     \end{tikzpicture}
    \caption{$\la=(7,7,7,6,3,3,2,1)$, the boundary of $\la$ is marked by arrows. Thus, $\partial(\la)=10101001 \mid 1101000$ and the anchor of $\partial(\la)$ is $8$.}\label{F:partial}
    \end{figure}
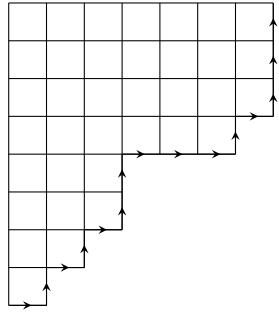

Now, we will identify a partition $\la$ with its {\em boundary sequence} $\partial(\la)$, which is a finite $0/1$ sequence constructed as follows: commence at the southwestern vertex of the diagram of $\la$, and traverse along the edges of the southeastern boundary proceeding towards the northeastern vertex; represent each horizontal movement by a $1$, and each vertical movement by a $0$. Clearly, $\la$ is determined by its boundary sequence $\partial(\la)$ and invariant under adding $0$'s (resp. $1$'s) to the head (resp. tail) of $\partial(\la)$. We define the {\em anchor} of $\partial(\la)=(s_1,s_2,\cdots,s_t)$ is the position $i$ such that $\#\{1\leq j\leq i\mid s_j=1\}=\#\{i+1\leq j\leq t\mid s_j=0\}$. We usually insert a vertical line ``$\mid$" into $\partial(\la)$ to indicate the anchor. See Figure \ref{F:partial} for an example. 

Now we describe the bijection $\Psi_m$ between the set of $m$-multipartitions of $n$ (denoted by $\mathcal{P}_{n,m}$) and the set of partitions of $mn$ with empty $m$-core (denoted by $\mathcal{P}^{0}_{mn}$). Here, the $m$-core of a partition $\la$ is obtained by removing as many $m$-ribbons from $\la$ as possible.

\begin{defn}\label{d:quotient}
   For a given $\bla=(\la^{(1)},\la^{(2)},\cdots,\la^{(m)})\in \mathcal{P}_{n,m}$, suppose $\partial(\la^{(i)})=(s^{(i)}_{1},\cdots,s^{(i)}_{d})$, $i=1,2,\cdots,m$. Here we add the leading $0$'s and trailing $1$'s to make $\partial(\la^{(i)})$ have the same length $d$ and the same anchor for all $1\leq i\leq m$. We let $\Psi_m(\bla)$ be the unique partition such that $\partial(\Psi_m(\bla))$ is
   \begin{align*}       s^{(1)}_{1},s^{(2)}_1,\cdots,s^{(m)}_1,s^{(1)}_{2},s^{(2)}_2,\cdots,s^{(m)}_2,\cdots,s^{(1)}_{d},s^{(2)}_d,\cdots,s^{(m)}_d.
   \end{align*}
We claim $\Psi_m$ is a bijection. Indeed, $\bla=(\la^{(1)},\la^{(2)},\cdots,\la^{(m)})$ is the $m$-quotient of $\Psi_m(\bla)$. Moreover, we have $|\Psi_m(\bla)|=m|\bla|$ (see also \cite[Remark 2.9]{AR}).    
\end{defn}

\begin{exmp}
    We choose $\bla$ as follows:
    \begin{align*}
\setlength{\arraycolsep}{15pt}
\left( \begin{array}{ccc}
\begin{tikzpicture}[scale=0.6]
   \coordinate (Origin)   at (0,0);
    \coordinate (XAxisMin) at (0,0);
    \coordinate (XAxisMax) at (4,0);
    \coordinate (YAxisMin) at (0,-3);
    \coordinate (YAxisMax) at (0,0);
\draw [thin, black] (0,0) -- (3,0);
    \draw [thin, black] (0,-1) -- (3,-1);
    \draw [thin, black] (0,-2) -- (1,-2);
    \draw [thin, black] (0,0) -- (0,-2);
    \draw [thin, black] (1,0) -- (1,-2);
    \draw [thin, black] (2,0) -- (2,-1);
    \draw [thin, black] (3,0) -- (3,-1);
    \end{tikzpicture},&
    \begin{tikzpicture}[scale=0.6]
   \coordinate (Origin)   at (0,0);
    \coordinate (XAxisMin) at (0,0);
    \coordinate (XAxisMax) at (4,0);
    \coordinate (YAxisMin) at (0,-3);
    \coordinate (YAxisMax) at (0,0);
\draw [thin, black] (0,0) -- (2,0);
    \draw [thin, black] (0,-1) -- (2,-1);
    \draw [thin, black] (0,-2) -- (1,-2);
    \draw [thin, black] (0,-3) -- (1,-3);
    \draw [thin, black] (0,0) -- (0,-3);
    \draw [thin, black] (1,0) -- (1,-3);
    \draw [thin, black] (2,0) -- (2,-1);
    \end{tikzpicture}~,&
    \begin{tikzpicture}[scale=0.6]
   \coordinate (Origin)   at (0,0);
    \coordinate (XAxisMin) at (0,0);
    \coordinate (XAxisMax) at (4,0);
    \coordinate (YAxisMin) at (0,-3);
    \coordinate (YAxisMax) at (0,0);
\draw [thin, black] (0,0) -- (2,0);
    \draw [thin, black] (0,-1) -- (2,-1);
    \draw [thin, black] (0,-2) -- (2,-2);
    \draw [thin, black] (0,0) -- (0,-2);
    \draw [thin, black] (1,0) -- (1,-2);
     \draw [thin, black] (2,0) -- (2,-2);
    \end{tikzpicture}
\end{array}  \right).
\end{align*}
The boundary sequence of each component of $\bla$ is $10\mid110$, $100\mid10$ and $11\mid00$, respectively. Add the leading $0$'s and tailing $1$'s to make sure they have the same length and the same anchor. So their boundary sequences can be presented as follows:
\begin{align*}
    010 \mid 110 \\ 
100 \mid 101\\ 
011 \mid 001 \\ 
\end{align*}
Therefore, the boundary sequence of $\Psi_{m}(\bla)$ is $010101001\mid110100011$. This is exactly the partition in Figure \ref{F:partial} when remove the leading $0$'s and the tailing $1$'s.
\end{exmp}

\begin{figure}
\begin{tikzpicture}[scale=0.5]
    \coordinate (Origin)   at (0,0);
    \coordinate (XAxisMin) at (0,0);
    \coordinate (XAxisMax) at (14,0);
    \coordinate (YAxisMin) at (0,0);
    \coordinate (YAxisMax) at (0,-14);
    \draw [thin, black] (0,0) -- (7,0);
    \draw [thin, black] (0,-1) --(7,-1);
    \draw [thin, black] (0,-2) -- (7,-2);
    \draw [thin, black] (0,-3) -- (7,-3);
    \draw [thin, black] (0,-4) -- (6,-4);
    \draw [thin, black] (0,-5) -- (3,-5);
    \draw [thin, black] (0,-6) -- (3,-6);
    \draw [thin, black] (0,-7) -- (2,-7);
     \draw [thin, black] (0,-8) -- (1,-8);
     \draw [thin, black] (0,0) -- (0,-8);
     \draw [thin, black] (1,0) -- (1,-8);
      \draw [thin, black] (2,0) -- (2,-7);
     \draw [thin, black] (3,0) -- (3,-6);
     \draw [thin, black] (4,0) -- (4,-4);
     \draw [thin, black] (5,0) -- (5,-4);
      \draw [thin, black] (6,0) -- (6,-4);
      \draw [thin, black] (7,0) -- (7,-3);
       \node at(-0.5,-0.5) {$1$};
       \node at(-0.5,-1.5) {$3$};
       \node at(-0.5,-2.5) {$2$};
       \node at(-0.5,-3.5) {$3$};
       \node at(-0.5,-4.5) {$2$};
       \node at(-0.5,-5.5) {$1$};
       \node at(-0.5,-6.5) {$2$};
       \node at(-0.5,-7.5) {$3$};
        \filldraw[fill = green]
     (4,-3) rectangle (5,-4)
     (5,-3) rectangle (6,-4)
     (5,-2) rectangle (6,-3) 
     (6,-2) rectangle (7,-3) 
     (6,-1) rectangle (7,-2)
     (6,0) rectangle (7,-1) 
     (1,-6) rectangle (2,-7)
     (1,-5) rectangle (2,-6) 
     (2,-5) rectangle (3,-6); 
     \end{tikzpicture}
    \caption{$\tilde{\partial}(\la)=010101001 \mid 110100011$, $k(\tilde{\partial}(\la))=9$ and the row-color sequence is marked in the leftmost column. Thus, the green generalized ribbon is colored by $1$.}\label{F:row-color}
    \end{figure}
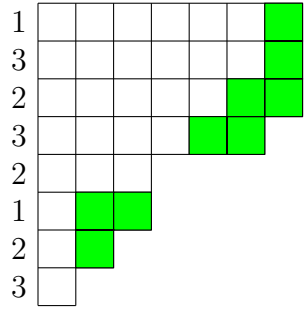

\begin{lem}\label{l:fact}
    Let $\la=(\la_1,\cdots,\la_l)$ $(\la_l\geq 0)$ be a partition with its boundary sequence $\partial(\la)=(s_1,s_2,\cdots)$. We collect the following facts from \cite{AR}:
    \begin{enumerate}
\item Removing a $k$-ribbon from $\la$ is equivalent to exchanging $s_{i}=1$ with $s_{i+k}=0$, for some $i$.
   \item The position in $\partial(\la)$ of the zero corresponding to $\la_i$ is equal to $\la_i+l-i+1$. 
    \end{enumerate}
\end{lem}

For a given partition $\la$ with empty $m$-core and length $l$, denote by $\tilde{\partial}(\la)$ the boundary sequence given by Definition \ref{d:quotient} (before removing the leading $0$'s and the tailing $1$'s). Note that $\tilde{\partial}(\la)$ relies on the choice of $\partial(\la^{(1)}), \partial(\la^{(2)}),\cdots,\partial(\la^{(m)})$. Define $k(\tilde{\partial}(\la)):=\#\{j\mid \tilde{\partial}(\la)_j=0\}$. We define the {\em row-color sequence} $a(\la)$ of $\la$ by $a(\la):=(a_1,\cdots,a_l)$, where $1\leq a_i\leq m$ and $a_i\equiv \la_i+k(\tilde{\partial}(\la))-i+1 ~({\rm mod}~ m)$, $i=1,2,\cdots,l$.\footnote{We remark that the definition of the row-color sequence is slightly different from that in \cite[Definition 3.3]{AR}.} We claim that $a(\la)$ is independent of the choice of $\tilde{\partial}(\la)$. Indeed, if we have another boundary sequence $\hat{\partial}(\la)$ obtained from Definition \ref{d:quotient} via distinct choice of $\partial(\la^{(1)}), \partial(\la^{(2)}),\cdots,\partial(\la^{(m)})$, then $m\mid(k(\tilde{\partial}(\la))-k(\hat{\partial}(\la)))$, which implies they gives the same $a(\la)$ by definition. 

We say the {\em starting point} (resp. {\em ending point}) of a ribbon is the northeastmost (resp. southwestmost) box of the ribbon. A ribbon is {\em colored by $i$} if the row-color of the row containing the starting point of the ribbon equals $i$. That is, if the starting point is in the $j$-th row then $a_j=i$. More generally, We say a generalized ribbon is colored by $i$ if all of its connected components are colored by $i$. See Figure \ref{F:row-color} for an example.

\begin{prop}\label{p:equiv}
    Let $\la$ be a partition with empty $m$-core and $m$-quotient $(\la^{(1)},\la^{(2)},\cdots,\la^{(m)})$. We have that removing a $k$-generalized ribbon with $p$ connected components from $\la^{(i)}$ is equivalent to removing a $km$-generalized ribbon with $p$ connected components colored by $i$ from $\la$. 
\end{prop}
\begin{proof}  
Suppose $\partial(\la^{(i)})=(s^{(i)}_1,s^{(i)}_2,\cdots)$ is the boundary sequence of $\la^{(i)}$ and $\tilde{\partial}(\la)=(s_1,s_2,\cdots)$ is the boundary sequence given by Definition \ref{d:quotient}. By Lemma \ref{l:fact} (1), removing a $k$-generalized ribbon with $p$ connected components from $\la^{(i)}$ is equivalent to exchanging $s^{(i)}_{i_1}=1$ with $s^{(i)}_{j_1}=0$, $s^{(i)}_{i_2}=1$ with $s^{(i)}_{j_2}=0$, $\cdots$, $s^{(i)}_{i_p}=1$ with $s^{(i)}_{j_p}=0$ in $\partial(\la^{(i)})$ for some $i_1<j_1<i_2<j_2<\cdots<i_p<j_p$ such that $\sum_{a=1}^p(j_a-i_a)=k$. It follows from the definition of $\Psi_m$ that this is further equivalent to exchanging $s_{u_1}=1$ with $s_{v_1}=0$, $s_{u_2}=1$ with $s_{v_2}=0$, $\cdots$, $s_{u_p}=1$ with $s_{v_p}=0$ such that $\sum_{a=1}^p(v_a-u_a)=km$ and $u_a\equiv v_a\equiv i~({\rm mod}~m)$, $1\leq a\leq p$. Lemma \ref{l:fact} (1) implies that this is equivalent to removing a $km$-generalized ribbon with $p$ connected components from $\la$. If we assume the starting points of the connected components of this $km$-generalized ribbon are respectively in the row $t_1, t_2,\cdots, t_p$ in $\la$, then Lemma \ref{l:fact} (2) yields that $\la_{t_a}+k(\tilde{\partial}(\la))-t_a+1=v_a$, $1\leq a\leq p$. So $\la_{t_a}+k(\tilde{\partial}(\la))-t_a+1\equiv i~({\rm mod}~m)$, $1\leq a\leq p$, i.e., this $km$-generalized ribbon is colored by $i$, which completes the proof.
\end{proof}

\begin{exmp}
 If we remove the generalized ribbon (its two components are in color yellow and green) from $\la^{(2)}=(2,1,1)$ as follows:   
\begin{align*}
\setlength{\arraycolsep}{15pt}
\left( \begin{array}{ccc}
\begin{tikzpicture}[scale=0.6]
   \coordinate (Origin)   at (0,0);
    \coordinate (XAxisMin) at (0,0);
    \coordinate (XAxisMax) at (4,0);
    \coordinate (YAxisMin) at (0,-3);
    \coordinate (YAxisMax) at (0,0);
\draw [thin, black] (0,0) -- (3,0);
    \draw [thin, black] (0,-1) -- (3,-1);
    \draw [thin, black] (0,-2) -- (1,-2);
    \draw [thin, black] (0,0) -- (0,-2);
    \draw [thin, black] (1,0) -- (1,-2);
    \draw [thin, black] (2,0) -- (2,-1);
    \draw [thin, black] (3,0) -- (3,-1);
    \end{tikzpicture},&
    \begin{tikzpicture}[scale=0.6]
   \coordinate (Origin)   at (0,0);
    \coordinate (XAxisMin) at (0,0);
    \coordinate (XAxisMax) at (4,0);
    \coordinate (YAxisMin) at (0,-3);
    \coordinate (YAxisMax) at (0,0);
\draw [thin, black] (0,0) -- (2,0);
    \draw [thin, black] (0,-1) -- (2,-1);
    \draw [thin, black] (0,-2) -- (1,-2);
    \draw [thin, black] (0,-3) -- (1,-3);
    \draw [thin, black] (0,0) -- (0,-3);
    \draw [thin, black] (1,0) -- (1,-3);
    \draw [thin, black] (2,0) -- (2,-1);
     \filldraw[fill = green]
     (0,-1) rectangle (1,-2) (0,-2) rectangle (1,-3); 
     \filldraw[fill = yellow]
     (1,0) rectangle (2,-1);
    \end{tikzpicture}~,&
    \begin{tikzpicture}[scale=0.6]
   \coordinate (Origin)   at (0,0);
    \coordinate (XAxisMin) at (0,0);
    \coordinate (XAxisMax) at (4,0);
    \coordinate (YAxisMin) at (0,-3);
    \coordinate (YAxisMax) at (0,0);
\draw [thin, black] (0,0) -- (2,0);
    \draw [thin, black] (0,-1) -- (2,-1);
    \draw [thin, black] (0,-2) -- (2,-2);
    \draw [thin, black] (0,0) -- (0,-2);
    \draw [thin, black] (1,0) -- (1,-2);
     \draw [thin, black] (2,0) -- (2,-2);
    \end{tikzpicture}
\end{array}  \right).
\end{align*}
This is equivalent to removing the following generalized ribbon colored by $2$ (its two components are in color yellow and green) from $(7,7,7,6,3,3,2,1)$:
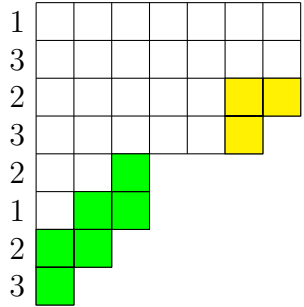
\begin{gather*}
\begin{tikzpicture}[scale=0.5]
    \coordinate (Origin)   at (0,0);
    \coordinate (XAxisMin) at (0,0);
    \coordinate (XAxisMax) at (14,0);
    \coordinate (YAxisMin) at (0,0);
    \coordinate (YAxisMax) at (0,-14);
    \draw [thin, black] (0,0) -- (7,0);
    \draw [thin, black] (0,-1) --(7,-1);
    \draw [thin, black] (0,-2) -- (7,-2);
    \draw [thin, black] (0,-3) -- (7,-3);
    \draw [thin, black] (0,-4) -- (6,-4);
    \draw [thin, black] (0,-5) -- (3,-5);
    \draw [thin, black] (0,-6) -- (3,-6);
    \draw [thin, black] (0,-7) -- (2,-7);
     \draw [thin, black] (0,-8) -- (1,-8);
     \draw [thin, black] (0,0) -- (0,-8);
     \draw [thin, black] (1,0) -- (1,-8);
      \draw [thin, black] (2,0) -- (2,-7);
     \draw [thin, black] (3,0) -- (3,-6);
     \draw [thin, black] (4,0) -- (4,-4);
     \draw [thin, black] (5,0) -- (5,-4);
      \draw [thin, black] (6,0) -- (6,-4);
      \draw [thin, black] (7,0) -- (7,-3);
       \node at(-0.5,-0.5) {$1$};
       \node at(-0.5,-1.5) {$3$};
       \node at(-0.5,-2.5) {$2$};
       \node at(-0.5,-3.5) {$3$};
       \node at(-0.5,-4.5) {$2$};
       \node at(-0.5,-5.5) {$1$};
       \node at(-0.5,-6.5) {$2$};
       \node at(-0.5,-7.5) {$3$};
        \filldraw[fill = green]
     (0,-7) rectangle (1,-8) (0,-6) rectangle (1,-7) (1,-6) rectangle (2,-7) (1,-5) rectangle (2,-6) (2,-5) rectangle (3,-6) (2,-4) rectangle (3,-5); 
     \filldraw[fill = yellow]
     (5,-3) rectangle (6,-4) (5,-2) rectangle (6,-3) (6,-2) rectangle (7,-3);
     \end{tikzpicture}
  \end{gather*}
\end{exmp}

For any finite sequence $\bm{a}=(a_1,a_2,\cdots)$, we define its {\em inversion number} by  \begin{align*}
        {\rm inv}(\bm{a}):=\# \{(i,j)\mid i<j, a_i>a_j\}.
    \end{align*}
Let $a(\la)=(a_1,a_2,\cdots,a_l)$ be the row-color sequence of $\la$. Define its {\em $m$-inversion number} by
    \begin{align*}
        {\rm inv}_m(\la):={\rm inv}(a(\la)).
    \end{align*}
For any finite sequence $\bm{a}=(a_1,\cdots,a_d)$, we assume that $\bm{b}=(b_1,\cdots,b_d)$ is a permutation of $\bm{a}$. Define $\I(\bm{a};\bm{b})$ by
\begin{align*}
    \I(\bm{a};\bm{b}):=\min\{{\rm inv}(\sigma)\mid \sigma\in\mathfrak{S}_d, ~\text{such that}~a_i=b_{\sigma(i)}~\text{for all}~ i=1,2,\cdots,d\}.
\end{align*}
Here the permutation $\sigma$ is viewed as a word (or sequence). For instance, 
$\sigma=(13254)$ is identified with the word $35214$ for the reason that
\begin{align*}
\sigma=
   \begin{pmatrix}1 &  2 & 3 & 4 & 5 \\
 3 & 5 & 2 & 1 &4
\end{pmatrix}. 
\end{align*}
For example, we choose $\bm{a}=(1,2,5,1,3,2,7)$, $\bm{b}=(2,5,1,3,2,7,1)$. Then the permutation $\sigma$ with the minimal inversion number satisfying $a_i=b_{\sigma(i)}$ is $3127456$. So $\I(\bm{a};\bm{b})=5$.

\begin{rem}\label{r:graph}
    $\I(\bm{a};\bm{b})$ can be computed graphically. Indeed, $\I(\bm{a};\bm{b})$ is equal to the minimal value of all possible applications of the following process:
    \begin{enumerate}
        \item Put sequences $\bm{a}$ and $\bm{b}$ in two rows and align them;
        \item  Connect the same number in the generic manner: no three lines intersect at one point;
         \item  Count the intersection number of the lines.
    \end{enumerate}
    In the above example, we can draw the lines, which produce the minimal intersection number as follows:
    
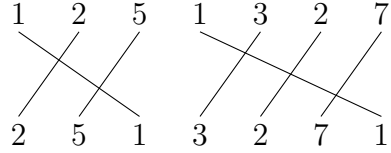
\begin{gather*}
     \begin{tikzpicture}[scale=0.8]
   \coordinate (Origin)   at (0,0);
    \coordinate (XAxisMin) at (0,0);
    \coordinate (XAxisMax) at (4,0);
    \coordinate (YAxisMin) at (0,-3);
    \coordinate (YAxisMax) at (0,0);
\draw [thin, black] (0,-0.3) -- (2,-1.7);
    \draw [thin, black] (1,-0.3) -- (0,-1.7);
    \draw [thin, black] (2,-0.3) -- (1,-1.7);
    \draw [thin, black] (3,-0.3) -- (6,-1.7);
    \draw [thin, black] (4,-0.3) -- (3,-1.7);
    \draw [thin, black] (5,-0.3) -- (4,-1.7);
    \draw [thin, black] (6,-0.3) -- (5,-1.7);
    \node at(0,0) {$1$}; \node at(1,0) {$2$}; \node at(2,0) {$5$}; \node at(3,0) {$1$}; \node at(4,0) {$3$}; \node at(5,0) {$2$};\node at(6,0) {$7$};
    \node at(0,-2) {$2$}; \node at(1,-2) {$5$}; \node at(2,-2) {$1$}; \node at(3,-2) {$3$}; \node at(4,-2) {$2$}; \node at(5,-2) {$7$};\node at(6,-2) {$1$};
    \end{tikzpicture} 
    \end{gather*}
    So $\I(\bm{a};\bm{b})=5$.
\end{rem}

\begin{lem}\label{l:inv}
    Let $\la^{'}$ be a partition obtained from $\la$ by removing a generalized ribbon with connected components $\xi_1,\xi_2,\cdots,\xi_p$. Suppose $u_i$ and $v_i$ are respectively the indices of the rows in $\la$ containing the starting point and ending point of $\xi_i$ $(1\leq i\leq p)$ such that $u_1<v_1<u_2<v_2<\cdots<u_p<v_p$. If we assume the row-color sequence of $\la$ is
    \begin{align*}
        a(\la)=(\cdots,\underbrace{a_{u_1},\cdots,a_{v_1}}_{v_1-u_1+1},\cdots,\underbrace{a_{u_2},\cdots,a_{v_2}}_{v_2-u_2+1},\cdots,\underbrace{a_{u_p},\cdots,a_{v_p}}_{v_p-u_p+1},\cdots).
    \end{align*}
  Then we have
\begin{align}
     a(\la^{'})=(\cdots,\underbrace{a_{u_{1}+1},\cdots,a_{v_1},a_{u_1}}_{v_1-u_1+1},\cdots,&\underbrace{a_{u_{2}+1},\cdots,a_{v_2},a_{u_2}}_{v_2-u_2+1},\cdots,\underbrace{a_{u_{p}+1},\cdots,a_{v_p},a_{u_p}}_{v_p-u_p+1},\cdots);\label{e:ala}\\
     \I(\la/\la^{'}):=\I(a(\la);a(\la^{'}))&=\sum_{i=1}^p\#\{j\mid u_i<j\leq v_i, a_j\neq a_{u_i}\}; \label{e:I}\\
     (-1)^{{\rm inv}_m(\la)-{\rm inv}_m(\la^{'})}&=(-1)^{\I(\la/\la^{'})}.\label{e:iI}
\end{align}
\end{lem}
\begin{proof}
    \eqref{e:ala} can be obtained by using the same arguments as those in \cite[Lemma 3.5]{AR}. \eqref{e:I} follows from \eqref{e:ala} by definition (or Remark \ref{r:graph}). Now we show \eqref{e:iI}. It follows from \eqref{e:ala} that 
    \begin{align*}
        {\rm inv}_m(\la)-{\rm inv}_m(\la^{'})=\sum_{i=1}^p\big(\#\{j\mid u_i<j\leq v_i, a_j< a_{u_i}\}-\#\{j\mid u_i<j\leq v_i, a_j> a_{u_i}\}\big),
    \end{align*}
 which has the same parity as $\I(\la/\la^{'})$ by \eqref{e:I}. This completes the proof.   
\end{proof}

\begin{defn}\label{d:eta'}
 Let $\la$ be a partition of $mn$, $\btau=(\tau^{(1)},\tau^{(2)},\cdots,\tau^{(t)})$ be a $t$-tuple of partitions with each partition having at most $m$ parts and $r$ be a positive integer. Define $\eta^{\la}_{\btau;r}$ as the sum of values obtained by all possible applications of the following algorithm:\\
 {\em Initialization}: $\varepsilon:=1$.
    \begin{enumerate}
        \item If $t=0$, end the program and output $\varepsilon$.
        \item If $l(\tau^{(t)})<m$, add $0$'s to the tail to make sure $\tau^{(t)}$ has $m$ components.
        \item Choose a permutation $\sigma\in \mathfrak{S}_m$. Remove a $m\tau^{(t)}_{\sigma(1)}$-generalized ribbon $\theta^{'}_1$ colored by $1$ from $\la$, and then remove a $m\tau^{(t)}_{\sigma(2)}$-generalized ribbon $\theta^{'}_2$ colored by $2$ from $\la/\theta^{'}_1$. Continue this process until we remove a $m\tau^{(t)}_{\sigma(m)}$-generalized ribbon $\theta^{'}_{m}$ colored by $m$ from $\la/(\theta^{'}_1\cup\cdots\cup\theta^{'}_{m-1})$. If there is no such $\sigma$ to do the removal process, set $\varepsilon=0$ and end the program.
        \item Multiply $\varepsilon$ by $u^r_{\m}(1-q^{-2})^{-1+\sum_{j=1}^{m}cc(\theta^{'}_j)}(-q^{-2})^{\sum_{j=1}^{m}(ht(\theta^{'}_j)-\I(\theta^{'}_j))}$, where $cc(\theta^{'}_j)$ (resp. $ht(\theta^{'}_j)$) is the number of connected components (resp. the height) of $\theta^{'}_j$. And $\m:=\max\{1\leq j\leq m\mid \tau^{(t)}_{\sigma(j)}\neq 0\}$.
        \item Redefine $\btau:=(\tau^{(1)},\tau^{(2)},\cdots,\tau^{(t-1)})$ and $t:=t-1$.
        \item Return to step (1). 
    \end{enumerate}
\end{defn}

Now we are ready to give a L\"ubeck--Prasad--Adin--Roichman formula for the cyclotomic Hecke algebra, computing the irreducible characters by using partitions instead of $m$-multipartitions.
\begin{thm}\label{t:AR}
Let $\bla=(\la^{(1)},\la^{(2)},\cdots,\la^{(m)}),\bmu\in \mathcal{P}_{n,m}$. We assume that $\la$ is the unique partition with empty $m$-core and $(\la^{(1)},\la^{(2)},\cdots,\la^{(m)})$ as its $m$-quotient. Then we have
\begin{align}\label{e:AR}
    \chi^{\bla}_{\bmu}=q^{n-l(\bmu)}\sum_{r=1}^m\sum_{\btau\in\mathcal{P}_m(\mu^{(r)})}\eta^{\la}_{\btau;r},
\end{align}
\end{thm}
\begin{proof}
    By \eqref{e:restate-MN}, Definition \ref{d:eta} and Proposition \ref{p:equiv}, it suffices to show that $\eta^{\bla}_{\btau;r}=\eta^{\la}_{\btau;r}$, for which it is sufficient to verify that the factors of corresponding removal process under the bijection established in Proposition \ref{p:equiv} differ by some power of $-q^{-2}$. To be precise, if we remove a $k$-generalized ribbon $\theta_i$ with height $ht(\theta_{i})$ from $\la^{(i)}$, which, due to Proposition \ref{p:equiv}, is equivalent to removing a $km$-generalized ribbon $\theta^{'}_i$ colored by $i$ with height $ht(\theta^{'}_{i})$ from $\la$, then we need to verify $ht(\theta_i)=ht(\theta^{'}_i)-\I(\theta^{'}_i)$. For the sake of argument, we make the following assumptions:
    \begin{itemize}
        \item $\theta_i$ (resp. $\theta^{'}_i$) has $p$ connected components $\xi_1,\xi_2,\cdots,\xi_p$ (resp. $\xi^{'}_1,\xi^{'}_2,\cdots,\xi^{'}_p$).
        \item Let $u_1<v_1<u_2<v_2<\cdots<u_p<v_p$ (resp. $u^{(i)}_1<v^{(i)}_1<u^{(i)}_2<v^{(i)}_2<\cdots<u^{(i)}_p<v^{(i)}_p$) be the indices of rows in $\la$ (resp. $\la^{(i)}$) such that the $u_q$-th (resp. $u^{(i)}_q$-th) row and the $v_q$-th (resp. $v^{(i)}_q$-th) row contain respectively the starting point and ending point of $\xi^{'}_q$ (resp. $\xi_{q}$), $q=1,2,\cdots,p.$
        \item Let $b_1<d_1<b_2<d_2<\cdots<b_p<d_p$ (resp. $b^{(i)}_1<d^{(i)}_1<b^{(i)}_2<d^{(i)}_2<\cdots<b^{(i)}_p<d^{(i)}_p$) be the indices of positions in $\partial(\la)$ (resp. $\partial(\la^{(i)})$) such that we swap $s_{b_q}=1$ with $s_{d_q}=0$ (resp. $s^{i}_{b^{(i)}_q}=1$ with $s^{(i)}_{d^{(i)}_q}=0$) when remove $\xi^{'}_q$ (resp. $\xi_q$) from $\la$ (resp. $\la^{(i)}$).
    \end{itemize}
    Thus, with the above notations in hand, we have 
    $$ ht(\xi^{'}_q)=\#\{t\mid s_t=0, b_q\leq t<d_q\}, \quad ht(\xi_q)=\#\{t\mid s^{(i)}_t=0, b^{(i)}_q\leq t<d^{(i)}_q\}.$$
    By the definition of $\Psi_m$, we hence have 
    $$\#\{t\mid s^{(i)}_t=0, b^{(i)}_q\leq t<d^{(i)}_q\}=\#\{t\mid s_t=0, b_q\leq t<d_q, t\equiv i~({\rm mod}~ m)\}.$$
    Therefore, 
    \begin{align*}
        ht(\theta^{'}_i)-ht(\theta_i)=&\sum_{q=1}^p \big(ht(\xi_q^{'})-ht(\xi_q)\big)\\
        =&\sum_{q=1}^p \#\{t\mid s_t=0, b_q\leq t<d_q, t\not\equiv i~({\rm mod}~ m)\}\\
        =&\sum_{q=1}^p \#\{t\mid u_q< t\leq v_q, a_t\neq a_{u_q}\}\quad \text{(by the definition of $a(\la)$)}\\
        =&\I(\theta^{'}_i)\quad \text{(by Lemma \ref{l:inv})}.
    \end{align*}

    This finishes the proof.
\end{proof}


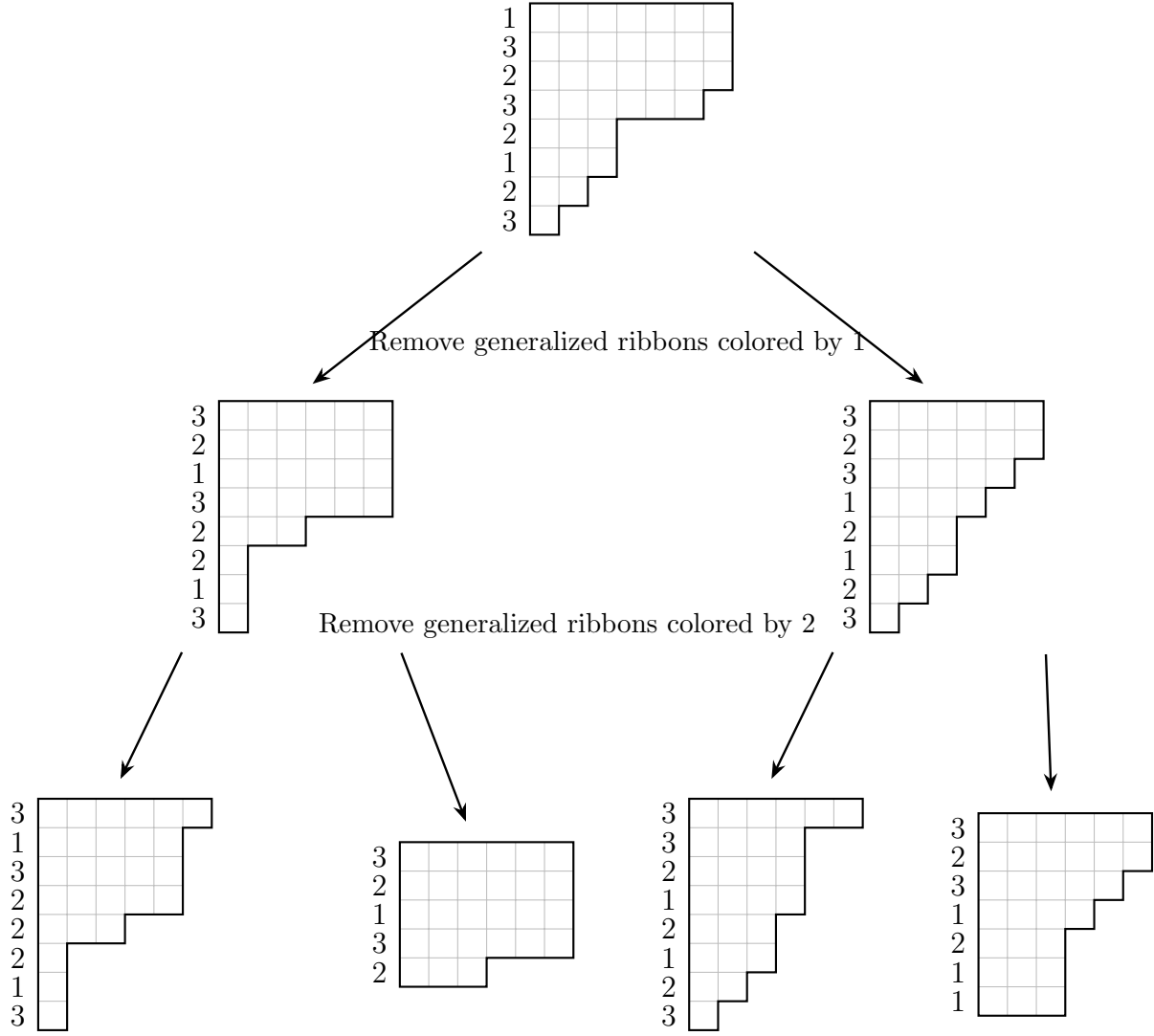
\begin{figure}[ht]
\centering
\begin{tikzpicture}[>=Stealth, line width=.9pt, every node/.style={inner sep=0pt}]

\node (T)  at (0,0)        {\YDTop};
\node (ML) at (-4.5,-5.5)  {\YDMidL};
\node (MR) at ( 4.5,-5.5)  {\YDMidR};

\node (B1) at (-7,-11)  {\YDBotI};
\node (B2) at (-2,-11)  {\YDBotII};
\node (B3) at ( 2,-11)  {\YDBotIII};
\node (B4) at ( 6,-11)  {\YDBotIV};

\node[font=\small] at (0,-3.1) {Remove generalized ribbons colored by $1$};
\node[font=\small] at (-0.7,-7) {Remove generalized ribbons colored by $2$};

\draw[->,shorten <=10pt,shorten >=10pt] (T.south west) -- (ML.north);
\draw[->,shorten <=10pt,shorten >=10pt] (T.south east) -- (MR.north);

\draw[->,shorten <=8pt,shorten >=8pt] (ML.south west) -- (B1.north);
\draw[->,shorten <=8pt,shorten >=8pt] (ML.south east) -- (B2.north);

\draw[->,shorten <=8pt,shorten >=8pt] (MR.south west) -- (B3.north);
\draw[->,shorten <=8pt,shorten >=8pt] (MR.south east) -- (B4.north);

\end{tikzpicture}
\caption{Four removal processes corresponding to those in Figure \ref{F:remove}. The sum terms they contribute are same as those in Figure \ref{F:remove}.}\label{F:remove2}
\end{figure}

\begin{exmp}
    Figure \ref{F:remove2} presents four removal processes corresponding to those of Figure \ref{F:remove}. Let us take the leftmost one as an example. We firstly remove a $6$-generalized ribbon $\theta^{'}_1$ colored by $1$, which has two connected components. So $ht(\theta^{'}_1)=2+1=3$. As we can see, $a(\la)=(1,3,2,3,2,1,2,3)$ and $a(\la\setminus\theta^{'}_1)=(3,2,1,3,2,2,1,3)$. So $\I(\theta^{'}_1)=3$. Similarly, we have $ht(\theta^{'}_2)=2$ and $\I(\theta^{'}_2)=2$. Thus, this process contribute the sum term $u^r_2(1-q^{-2})^3(-q^{-2})^0$.
\end{exmp}

\subsection{Specialization to $\mathscr{W}_{m,n}$ case} This subsection aims to explain how our result (Theorem \ref{t:AR}) reduces to the special case for $\mathscr{W}_{m,n}$, which was considered in \cite{AR}. As we mentioned in section \ref{s:special}, the irreducible character $\phi^{\bm{\la}}$ can be obtained from $\chi^{\bla}$ by specializing $\bm{u}=\bm{\zeta}=(1,\zeta,\zeta^2,\cdots,\zeta^{m-1})$ and $q=1$, where $\zeta$ is the primitive $m$-th root of unity. In this case, $\eta^{\la}_{\btau;m}=0$ unless each partition $\tau^{(i)}$ ($1\leq i\leq l(\mu^{(m)})$) has only one row, i.e., $\tau^{(i)}=(\mu^{(m)}_i)$. Moreover, only those processes of removing ribbons (not generalized ribbons) from $\la$ can contribute the summation. Therefore, if $\btau\in\mathcal{P}_m(\mu^{(m)})$, then $\eta^{\la}_{\btau;m}\Big|_{\substack{
\bm{u}=\bm{\zeta}\\ q=1}}$ equals the sum of values obtained by all possible applications of the following algorithm:\\
    {\em Initialization}: $\varepsilon:=1$.
    \begin{enumerate}
        \item If $t=0$, end the program and output $\varepsilon$.
\item Remove a $m\mu^{(m)}_t$-ribbon $\theta$ from $\la$. If we can not find such a ribbon, set $\varepsilon=0$ and end the program.
     \item Multiply $\varepsilon$ by $(-1)^{ht(\theta)-\I(\theta)}$. 
     \item Redefine $\btau:=\big((\mu^{(m)}_1),(\mu^{(m)}_2),\cdots,(\mu^{(m)}_{t-1})\big)$ and $t:=t-1$.
        \item Return to step (1). 
      \end{enumerate}
       Let $\varnothing=\la^{(0)}\subset \la^{(1)}\subset\cdots\subset\la^{(r)}=\la$ be a chain such that $\la^{(i+1)}/\la^{(i)}$ is a ribbon. Then it follows from equation \eqref{e:iI} that $(-1)^{\sum_{i=0}^{r-1}\I(\la^{(i+1)}/\la^{(i)})}=(-1)^{\sum_{i=0}^{r-1}\big({\rm inv}_m(\la^{(i+1)})-{\rm inv}_m(\la^{(i)})\big)}=(-1)^{{\rm inv}_m(\la)-{\rm inv}_m(\varnothing)}$. This, together with the above algorithm, implies that 
       \begin{align}\label{e:special-eta}
         \eta^{\la}_{\btau;m}\Big|_{\substack{
\bm{u}=\bm{\zeta}\\ q=1}}=(-1)^{{\rm inv}_m(\la)-{\rm inv}_m(\varnothing)} \gamma^{\la}_{m\mu^{(m)}}.
       \end{align}
Here $\gamma^{\la}$ is the irreducible character of the symmetric group indexed by the partition $\la$ with $\bla=(\la^{(1)},\la^{(2)},\cdots,\la^{(m)})$ as its $m$-quotient, $m\mu^{(m)}=(m\mu^{(m)}_1, m\mu^{(m)}_2,\cdots)$ and $\btau\in\mathcal{P}_m(\mu^{(m)})$.

Summarizing above, we get the special case of Theorem \ref{t:AR} as below, which firstly appears in \cite[Theorem 3.1]{AR}.
\begin{cor}
    Let $\la\vdash mn, \mu\vdash n$. Suppose the $m$-quotient of $\la$ is $(\la^{(1)},\la^{(2)},\cdots,\la^{(m)})$. Then we have
    \begin{align}\label{e:special-AR}
    \phi^{(\la^{(1)},\la^{(2)},\cdots,\la^{(m)})}_{(\varnothing,\cdots,\varnothing,\mu)}=(-1)^{{\rm inv}_m(\la)-{\rm inv}_m(\varnothing)}\gamma^{\la}_{m\mu},
    \end{align}
    where $\phi^{(\la^{(1)},\la^{(2)},\cdots,\la^{(m)})}$ (resp. $\gamma^{\la}$) is the irreducible character of the complex reflection group (resp. the symmetric group) indexed by $(\la^{(1)},\la^{(2)},\cdots,\la^{(m)})$ (resp. $\la$).
\end{cor}
\begin{proof}
    Taking specialization $\bm{u}=\bm{\zeta}$, $q=1$ and $\bmu=(\varnothing,\cdots,\varnothing,\mu)$ on both sides of \eqref{e:AR} gives that
    \begin{align}
    \phi^{(\la^{(1)},\la^{(2)},\cdots,\la^{(m)})}_{(\varnothing,\cdots,\varnothing,\mu)}=\eta^{\la}_{((\mu^{(m)}_1),(\mu^{(m)}_2),\cdots,);m}\Big|_{\substack{
\bm{u}=\bm{\zeta}\\ q=1}},
    \end{align}
    which, together with \eqref{e:special-eta}, yields \eqref{e:special-AR}. As desired. 
\end{proof}

\begin{rem}
It would be interesting to find a representation-theoretic interpretation of
\(
\sum_{\btau\in \mathcal{P}_m(\mu^{(m)})}\eta^{\la}_{\btau;m}
\)
(up to an explicit normalization factor, e.g. involving $q^{n-l(\mu^{(m)})}$ and possibly a power of $-q^{-2}$),
for $\la$ with empty $m$-core and $m$-quotient
$(\la^{(1)},\la^{(2)},\cdots,\la^{(m)})$.
More precisely, one may ask whether this quantity arises as the value of a natural character (or virtual character) associated with a suitable algebra related to $\mathscr{H}_{m,n}(q,\bm{u})$, for example a subalgebra, a corner algebra, or a quotient/subquotient, evaluated at an element corresponding to $g((\varnothing,\cdots,\varnothing,\mu^{(m)}))$.
In the $\mathscr{W}_{m,n}$ case, this interpretation is given by the symmetric group $\mathfrak{S}_n$.
\end{rem}


\section{Multiple bitrace for the irreducible characters}\label{s:mbtr}
In this section, we firstly introduce multiple bitrace for the irreducible characters of the cyclotomic Hecke algebra. Subsequently, a combinatorial formula for it will be given. Finally, we will consider the special case for $\mathscr{W}_{m,n}$, which gives the second orthogonality relation for the irreducible
characters of $\mathscr{W}_{m,n}$. For this reason, the multiple bitrace can be also viewed as the second orthogonality relation for $\chi^{\bla}_{\bmu}$.

To proceed with this, let us recall the basis of this algebra. Consider the subalgebra $\mathscr{W}^{0}_{m,n}$ of $\mathscr{W}_{m,n}$ generated by $s_2,s_3,\cdots,s_n$. Then $\mathscr{W}^0_{m,n}:=<s_2,s_3,\cdots,s_n>\cong \mathfrak{S}_n$. Similarly, the subalgebra $\mathscr{H}^0_{m,n}:=<g_2,g_3,\cdots,g_{n}>$ is isomorphic to the Hecke algebra in type $A$.   For any $\sigma\in \mathfrak{S}_n$, if $\sigma=s_{i_1}s_{i_2}\cdots s_{i_r}$ is any reduced decomposition of $\sigma$ in $\mathscr{W}^0_{m,n}$, then we define $T_{\sigma}:=g_{i_1}g_{i_2}\cdots g_{i_r}$. The isomorphism above guarantees this definition is well-defined. We introduce the {\em Jucys–Murphy elements} by 
\begin{align}
    J_1:=g_1, \quad J_i:=g_i\cdots g_2g_1g_2\cdots g_i, \quad i=2,3,\cdots,n.
\end{align}
More generally, $J^{\alpha}:=J_1^{\alpha_1}\cdots J_n^{\alpha_n}$ for a composition $\alpha$. It is shown in \cite{AK} that $\{J^{\alpha}T_{\sigma}|\alpha\in \mathbb{Z}^{n}_{+}, 0\leq\alpha_i\leq m-1, \sigma\in\mathfrak{S}_n\}$ forms a basis of $\mathscr{H}_{m,n}(q,\bm{u})$. 
Recently a basis for the cocenter of $\mathscr{H}_{m,n}(q,\bm{u})$ was provided, resolving the center conjecture for the cyclotomic Hecke algebra \cite{Hu}.

Denote by $\mathbb{Z}^n_+(m)$ the set of all $\alpha\in\mathbb{Z}^n_+$ with each component of $\alpha$ lying in $\{0,1,\cdots,m-1\}$. For any $\alpha,\beta\in\mathbb{Z}^n_+(m)$ and $\sigma,\omega\in\mathfrak{S}_n$, we define the {\em multiple bitrace} of $J^{\alpha}T_{\sigma}$ and $J^{\beta}T_{\omega}$
by
\begin{align}\label{e:defmbtr}
    \mbtr(J^{\alpha}T_{\sigma}, J^{\beta}T_{\omega}):=\sum_{\substack{\theta\in\mathfrak{S}_n\\
    \gamma\in\mathbb{Z}^n_+(m)}}J^{\alpha}T_{\sigma}J^{\gamma}T_{\theta}J^{\beta}T_{\omega}|_{J^{\gamma}T_{\theta}},
\end{align}
where $J^{\alpha}T_{\sigma}J^{\theta}T_{\gamma}J^{\beta}T_{\omega}|_{J^{\theta}T_{\gamma}}$ denotes the coefficient of $J^{\theta}T_{\gamma}$ in the multiplication $J^{\alpha}T_{\sigma}J^{\theta}T_{\gamma}J^{\beta}T_{\omega}$. 

We can rewrite \eqref{e:defmbtr} in terms of operators as follows:
\begin{align}\label{e:opLR}
    \mbtr(J^{\alpha}T_{\sigma}, J^{\beta}T_{\omega})=\Tr(L_{\alpha,\sigma}\circ R_{\beta,\omega}),
\end{align}
where $L_{\alpha,\sigma}$ and $R_{\beta,\omega}$ represent the left and the right multiplication operator of $J^{\alpha}T_{\sigma}$ and $J^{\beta}T_{\omega}$ respectively. These two actions mutually commute with each other.

By left and right multiplication, $\mathscr{H}_{m,n}(q,\bm{u})$ can be viewed as a $\mathscr{H}_{m,n}(q,\bm{u})\times \mathscr{H}_{m,n}(q,\bm{u})$-bimodule.  As a bimodule, $\mathscr{H}_{m,n}(q,\bm{u})$ decomposes into the direct sum of simple   $\mathscr{H}_{m,n}(q,\bm{u})\times \mathscr{H}_{m,n}(q,\bm{u})$-bimodules as follows:
\begin{align}\label{e:decomp}
    \mathscr{H}_{m,n}(q,\bm{u})=\bigoplus_{\bm{\la}\in \mathcal{P}_{n,m}}V^{\bm{\la}}\otimes V_{\bm{\la}},
\end{align}
where $V^{\bm{\la}}$ (resp. $V_{\bm{\la}}$) denotes the irreducible left (resp. right) module indexed by $\bm{\la}$. Note that $V^{\bm{\la}}\cong V_{\bm{\la}}$. Taking trace on both sides in \eqref{e:decomp}, together with \eqref{e:opLR}, gives that
\begin{align}
    \mbtr(J^{\alpha}T_{\sigma}, J^{\beta}T_{\omega})=\sum_{\bm{\la}\in \mathcal{P}_{n,m}}\chi^{\bm{\la}}(J^{\alpha}T_{\sigma})\chi^{\bm{\la}}(J^{\beta}T_{\omega}).
\end{align}
As is mentioned in subsection \ref{ss:Fro}, $\chi^{\bm{\la}}$ can be determined completely by its values on $g(\bm{\mu})$. For this reason, particularly, we define 
\begin{align}\label{e:character}
    \mbtr(\bm{\mu}, \bm{\nu}):=\mbtr(g(\bm{\mu}), g(\bm{\nu}))=\sum_{\bm{\la}\in \mathcal{P}_{n,m}}\chi^{\bm{\la}}_{\bm{\mu}}\chi^{\bm{\la}}_{\bm{\nu}}.
\end{align}
By the Frobenius formula \eqref{e:fro} and the orthogonality of $s_{\bm{\la}}$ with respect to $\langle \cdot,\cdot \rangle_{\otimes}$, we have
\begin{align}\label{e:minner}
    \mbtr(\bm{\mu}, \bm{\nu})=\sum_{\bm{\la}\in \mathcal{P}_{n,m}}\chi^{\bm{\la}}_{\bm{\mu}}\chi^{\bm{\la}}_{\bm{\nu}}=\langle q_{\bm{\mu}}(x;q,\bm{u}), q_{\bm{\nu}}(x;q,\bm{u}) \rangle_{\otimes}.
\end{align}
Now we are in the position to compute this inner product. To do this, we need to rewrite $q_{\bm{\mu}}(x;q,\bm{u})$. 

For a $m$-multimatrix $\bm{M}=(M^{(1)},M^{(2)},\cdots,M^{(m)})$, define
\begin{align*}
    q_{\bm{M}}(x;q):=\prod_{i=1}^m\prod_{m^{(i)}_{kj}\in M^{(i)}}q_{m^{(i)}_{kj}}(x^{(j)};q).
\end{align*}
In the example \eqref{e:example}, $$q_{\bm{M}}(x;q)=q_3(x^{(1)};q)q^2_2(x^{(1)};q)q_8(x^{(1)};q)q^2_2(x^{(2)};q)q_1(x^{(2)};q)q_4(x^{(2)};q)q_6(x^{(3)};q)q_2(x^{(3)};q).$$

With the above notation in hand, we are ready to rewrite $q_{\bm{\mu}}(x;q,\bm{u})$ as follows:   
\begin{align}\label{e:rewrite}
\begin{split}
    q_{\bm{\mu}}(x;q,\bm{u})&=\prod_{i=1}^{m}\prod_{j=1}^{l(\mu^{(i)})}q_{\mu^{(i)}_j}^{(i)}(x;q,\bm{u})\\
    &=\prod_{i=1}^{m}\prod_{j=1}^{l(\mu^{(i)})}\frac{q^{\mu^{(i)}_j}}{q-q^{-1}}\sum_{c\in C_{\mu^{(i)}_{j},m}}u^{i}_{c}\prod_{r=1}^m q_{c_r}(x^{(r)};q^{-2})\\
    &=\frac{q^{n}}{(q-q^{-1})^{l(\bm{\mu})}}\sum_{\bm{M}} \bm{u}_{\bm{M}}q_{\bm{M}}(x;q^{-2})
    \end{split}
\end{align}
where $\bm{M}$ runs over all $m$-multimatrices such that the $j$-th row sum of the $i$-th matrix is equal to $\mu^{(i)}_j$.  The last equation holds because we can treat $c\in C_{\mu^{(i)}_j,m}$ as the $j$-th row of the $i$-th matrix. 

For a matrix $A$, we define its {\em $q$-weight} by
\begin{align*}
    \wt_q(A):=\prod_{a_{ij}\in A}(a_{ij})_{q},
\end{align*}
where $(k)_{q}=\begin{cases}
    (q-1)\frac{q^{2k}-1}{q+1}=(q-1)^2[k]_{q^2}, & \text{if $k>0$}\\
    1, & \text{if $k=0$}.
\end{cases}$

\begin{lem}\label{l:innerq}{\cite[Corollary 5.5]{JL1}}
    Assume that $\tau=(\tau_1,\tau_2,\cdots,\tau_r)$ and $\rho=(\rho_1,\rho_2,\cdots,\rho_s)$ are two compositions, then
    \begin{align}\label{e:innerq}
        \left\langle q_{\tau}(x^{(j)};q), q_{\rho}(x^{(j)};q)\right\rangle=\sum_{A}\wt_q(A)
    \end{align}
    summed over all $r\times s$ matrices $A$ such that its entries are non-negative integers and the row (resp. column) sums are $\tau_i$, $i=1,2,\cdots,r$ (resp. $\rho_i$, $i=1,2,\cdots,s$). 
\end{lem}

For a $m$-multimatrix $\bm{M}$, define the composition $c_{i}(\bm{M})$ ($i=1,2,\cdots,m$) by
\begin{align*}
    c_i(\bm{M}):=(m^{(1)}_{1,i},\cdots,m^{(1)}_{r(M^{(1)}),i},m^{(2)}_{1,i},\cdots,m^{(2)}_{r(M^{(2)}),i},\cdots,m^{(m)}_{1,i},\cdots,m^{(m)}_{r(M^{(m)}),i}).
\end{align*}
Recall that $\bm{r}(\bm{M}):=\sum_{i=1}^m r(M^{(i)})$, i.e., the total number of the rows of the matrices in $\bm{M}$. The following lemma is a multiple version of Lemma \ref{l:innerq}.

\begin{lem}
    Let $\bm{M}, \bm{N}\in \mathcal{M}_m$ be two $m$-multimatrices, then we have
    \begin{align}\label{e:mwt}
        \left\langle q_{\bm{M}}(x;q), q_{\bm{N}}(x;q)\right\rangle_{\otimes}=\sum_{\bm{A}}\prod_{i=1}^{m}\wt_q(A^{(i)})
    \end{align}
    summed over all $m$-tuples of matrices $\bm{A}=(A^{(1)},A^{(2)},\cdots,A^{(m)})$ such that
    \begin{enumerate}[label=$(\alph*)$]
        \item each $A^{(i)}$ is a $\bm{r}(\bm{M})\times \bm{r}(\bm{N})$ matrix with non-negative entries;
        \item the row (resp. column) sum of $A^{(i)}$ is $c_i(\bm{M})$ (resp. $c_i(N)$).
    \end{enumerate}
\end{lem}
\begin{proof}
   This can be obtained immediately as follows 
   \begin{align*}
       \langle q_{\bm{M}}(x;q), q_{\bm{N}}(x;q)\rangle_{\otimes}
       =&\left\langle \prod_{i=1}^m\prod_{m^{(i)}_{kj}\in M^{(i)}}q_{m^{(i)}_{kj}}(x^{(j)};q), 
       \prod_{i=1}^m\prod_{n^{(i)}_{kj}\in N^{(i)}}q_{n^{(i)}_{kj}}(x^{(j)};q)\right\rangle_{\otimes} \quad \text{(by definition)}\\
       =&\left\langle \prod_{j=1}^m\prod_{i,k}q_{m^{(i)}_{kj}}(x^{(j)};q), 
       \prod_{j=1}^m\prod_{i,k}q_{n^{(i)}_{kj}}(x^{(j)};q)\right\rangle_{\otimes}\\
       =&\prod_{j=1}^m\left\langle \prod_{i,k}q_{m^{(i)}_{kj}}(x^{(j)};q), 
       \prod_{i,k}q_{n^{(i)}_{kj}}(x^{(j)};q)\right\rangle\quad \text{(by \eqref{e:tensor})}\\
       =&\prod_{j=1}^m\left(\sum_{A^{(j)}}\wt_q(A^{(j)}) \right)\quad \text{(by Lemma \ref{l:innerq})}\\
       =&\sum_{\bm{A}}\prod_{i=1}^{m}\wt_q(A^{(i)}).
   \end{align*}
   $(a)$, $(b)$ and $(c)$ are given by the conditions satisfied by $A$ in \eqref{e:innerq}.
\end{proof}

For convenience, we say a $m$-tuple of matrices $\bm{A}$ is associated with $(\bm{M},\bm{N})$ if $\bm{A}$ satisfies conditions $(a)$, $(b)$ and $(c)$. For a $m$-multipartition $\bmu$, denote by $\mathcal{M}_m(\bmu)$ the set of all $m$-multimatrices $\bm{M}=(M^{(1)},\cdots,M^{(m)})$ such that the $j$-th row sum of $M^{(i)}$ is $\mu^{(i)}_{j}$, $i=1,2,\cdots,m$, $j=1,2,\cdots,l(\mu^{(i)})$.

Now we are ready to give an explicit combinatorial formula for $\mbtr(\bm{\mu},\bm{\nu})$.
\begin{thm}\label{t:mbtr}
    Let $\bm{\mu}, \bm{\nu}\in \mathcal{P}_{n,m}$, then
    \begin{align}\label{e:mbtr}
        \mbtr(\bm{\mu},\bm{\nu})=\frac{q^{2n}}{(q-q^{-1})^{l(\bm{\mu})+l(\bm{\nu})}}\sum_{\bm{M},\bm{N},\bm{A}}\bm{u}_{\bm{M}}\bm{u}_{\bm{N}}\prod_{i=1}^m\wt_{q^{-2}}(A^{(i)})
    \end{align}
    summed over all pairs of $m$-multimatrices $(\bm{M},\bm{N})$ and all $m$-tuples of matrices $\bm{A}$ such that
    \begin{itemize}
    \item $\bm{M}\in \mathcal{M}_m(\bmu)$ and $\bm{N}\in \mathcal{M}_m(\bm{\nu})$;
         \item $\bm{A}$ is associated with $(\bm{M},\bm{N})$, i.e., $\bm{A}$ satisfies conditions $(a)$, $(b)$ and $(c)$.
    \end{itemize}
       
\end{thm}
\begin{proof}
    One can obtain this result by combining \eqref{e:minner}, \eqref{e:rewrite} and \eqref{e:mwt}.
\end{proof}

\subsection{The second orthogonality relation of $\mathscr{W}_{m,n}$} 
In this subsection, we will present how the multiple bitrace defined above gives the second orthogonality relation of $\mathscr{W}_{m,n}$ by taking specialization $q=1$ and $\bm{u}=(1,\zeta,\cdots,\zeta^{m-1})$.

It follows from \eqref{e:character} and \eqref{e:minner} that
\begin{align}
  \mbtr(g(\bm{\mu}), g(\bm{\nu}))=\sum_{\bm{\la}\in \mathcal{P}_{n,m}}\chi^{\bm{\la}}(g(\bm{\mu}))\chi^{\bm{\la}}(g(\bm{\nu}))=\langle q_{\bm{\mu}}(x;q,\bm{u}), q_{\bm{\nu}}(x;q,\bm{u}) \rangle_{\otimes}.
\end{align}
Taking specialization $q=1$ and $\bm{u}=(1,\zeta,\cdots,\zeta^{m-1})$ on both sides implies that
\begin{align}
  \mbtr(w(\bm{\mu}), w(\bm{\nu}))\big|_{\substack{
\bm{u}=\bm{\zeta}\\ q=1}}=\sum_{\bm{\la}\in \mathcal{P}_{n,m}}\phi^{\bm{\la}}(w(\bm{\mu}))\phi^{\bm{\la}}(w(\bm{\nu}))=\langle p_{\bm{\mu}}(x), p_{\bm{\nu}}(x) \rangle_{\otimes}
\end{align}
where $w(\bm{\mu})$ is defined similarly to that of $g(\bm{\mu})$ and $w(\bm{\mu})=g(\bm{\mu})\big|_{\substack{
\bm{u}=\bm{\zeta}\\ q=1}}$ (see \cite[(6.1.3)]{Sho} for more details). $\phi^{\bm{\la}}$ is the irreducible characters of $\mathscr{W}_{m,n}$ indexed by $\bla$ and $p_{\bm{\mu}}(x)$ is defined in \eqref{e:def-p}.
Thus,
\begin{align}
  \mbtr(w(\bm{\mu}), w(\bm{\nu})^{-1})\big|_{\substack{
\bm{u}=\bm{\zeta}\\ q=1}}=\sum_{\bm{\la}\in \mathcal{P}_{n,m}}\phi^{\bm{\la}}(w(\bm{\mu}))\phi^{\bm{\la}}(w(\bm{\nu})^{-1})=\sum_{\bm{\la}\in \mathcal{P}_{n,m}}\phi^{\bm{\la}}_{\bm{\mu}}\bar{\phi}^{\bm{\la}}_{\bm{\nu}}=\langle p_{\bm{\mu}}(x), \bar{p}_{\bm{\nu}}(x) \rangle_{\otimes}
\end{align}
where $\bar{\phi}^{\bm{\la}}_{\bm{\nu}}$ (resp. $\bar{p}_{\bm{\nu}}(x)$) is the complex conjugate of $\phi^{\bm{\la}}_{\bm{\nu}}$ (resp. $p_{\bm{\nu}}(x)$).

We use the shorthand notation $\mbtr_1(\bmu,\bar{\bnu}):=\mbtr(w(\bm{\mu}), w(\bm{\nu})^{-1})\big|_{\substack{
\bm{u}=\bm{\zeta}\\ q=1}}$. Thus, by \eqref{e:mbtr}
\begin{align}\label{e:eq1}
  \mbtr_1(\bmu,\bar{\bnu})=\frac{q^{2n}}{(q-q^{-1})^{l(\bm{\mu})+l(\bm{\nu})}}\sum_{\bm{M},\bm{N},\bm{A}}\bm{\zeta}_{\bm{M}}\bm{\zeta}^{-1}_{\bm{N}}\prod_{i=1}^m\wt_{q^{-2}}(A^{(i)})\big|_{\substack{q=1}}
\end{align}
where $\bm{\zeta}_{\bm{M}}$ is defined similarly to that of $\bm{u}_{\bm{M}}$, just replacing $\bm{u}=(u_1,\cdots,u_m)$ by $\bm{\zeta}=(1,\cdots,\zeta^{m-1})$. Here the summation runs over all $\bm{M}\in \mathcal{M}_m(\bmu)$, $\bm{N}\in \mathcal{M}_m(\bm{\nu})$ and $\bm{A}$ is associated with $(\bm{M},\bm{N})$.

For a $m$-multimatrix $\bm{M}\in \mathcal{M}_m(\bm{\mu})$, recall the composition $c_{i}(\bm{M})$ ($i=1,2,\cdots,m$) by
\begin{align*}
    c_i(\bm{M}):=(m^{(1)}_{1,i},\cdots,m^{(1)}_{l(\mu^{(1)}),i},m^{(2)}_{1,i},\cdots,m^{(2)}_{l(\mu^{(2)}),i},\cdots,m^{(m)}_{1,i},\cdots,m^{(m)}_{l(\mu^{(m)}),i}).
\end{align*}
For any composition $\tau=(\tau_1,\tau_2,\cdots)$, we denote by $\overrightarrow{\tau}$ the partition obtained by rearranging $\tau$.

For given pair of $(\bm{M},\bm{N})$ satisfying that $\bm{M}\in \mathcal{M}_m(\bmu)$ and $\bm{N}\in \mathcal{M}_m(\bm{\nu})$, let us denote 
\begin{align}
    A(\bm{M},\bm{N}):=\frac{q^{2n}}{(q-q^{-1})^{l(\bm{\mu})+l(\bm{\nu})}}\sum_{\bm{A}}\prod_{i=1}^m\wt_{q^{-2}}(A^{(i)})\big|_{\substack{q=1}}
\end{align}
summed over all $m$-tuples of matrices $\bm{A}$ such that $\bm{A}$ is associated with $(\bm{M},\bm{N})$.

By the definition of $\wt_{q^{-2}}(A^{(i)})$, we have $A(\bm{M},\bm{N})=0$ unless there exists at least one $m$-tuple of matrices $\bm{A}$ satisfying that each row and each column contains at most one non-zero entry in $A^{(i)}$ for all $i=1,2\cdots,m.$ This means each row of $M^{(i)}$ and $N^{(i)}$ contains exactly one non-zero entry and $\overrightarrow{c}_i(\bm{M})=\overrightarrow{c}_i(\bm{N})$ for all $i=1,2,\cdots,m$, which implies that the multisets consisting of parts of $\bmu$ and $\bnu$ are consistent. We say $\bm{M}\sim\bm{N}$ (resp. $\bmu\sim\bnu$) if $\bm{M}$ and $\bm{N}$ (resp. $\bmu$ and $\bnu$) satisfying the corresponding condition mentioned above. Under this condition $\bm{M}\sim\bm{N}$, we have
\begin{align}
    \begin{split}
A(\bm{M},\bm{N})=&\sum_{\bm{A}}\prod_{a^{(r)}_{ij}\in\bm{A}}a^{(r)}_{ij}=\sum_{\bm{A}}\prod_{m^{(r)}_{ij}\in\bm{M}}m^{(r)}_{ij}=\sum_{\bm{A}}\prod_{r=1}^m\prod_{i=1}^{l(\mu^{(r)})}\mu^{(r)}_{i}\\=&\prod_{r=1}^m\prod_{i=1}^{l(\mu^{(r)})}\mu^{(r)}_{i}\#\{\bm{A}\mid \text{$\bm{A}$ is associated with $(\bm{M},\bm{N})$}\}.
    \end{split}
\end{align}

Therefore, by $\bmu\sim\bnu$,
\begin{align}
\begin{split}
\mbtr_1(\bmu,\bar{\bnu})=&\sum_{\bm{M},\bm{N}}\bm{\zeta}_{\bm{M}}\bm{\zeta}^{-1}_{\bm{N}}A(\bm{M},\bm{N})\\
=&\sum_{\nu^{(r)}_j}(1+\zeta^{1-r}+\cdots+\zeta^{(1-r)(m-1)})\mu^{(1)}_1\sum_{\bm{M}',\bm{N}'}\bm{\zeta}_{\bm{M}'}\bm{\zeta}^{-1}_{\bm{N}'}A(\bm{M}',\bm{N}')
\end{split}
\end{align}
summed over all $\nu^{(r)}_j$ in $\bnu$ such that $\nu^{(r)}_j=\mu^{(1)}_1$, $\bm{M}'\in \mathcal{M}_m(\bmu\setminus\mu^{(1)}_1)$ and $\bm{N}'\in \mathcal{M}_m(\bnu\setminus\nu^{(r)}_j)$.

Note that the fact:
\begin{align}
    1+\zeta^{n_1-n_2}+\cdots+\zeta^{(n_1-n_2)(m-1)}=
    \begin{cases}
        0, & \text{if $n_1\neq n_2$}\\
        m, &\text{if $n_1=n_2$},
    \end{cases}
\end{align}
which yields that $\mbtr_1(\bmu,\bar{\bnu})=0$ unless there exists at least one part $\nu^{(1)}_j$ in $\nu^{(1)}$ satisfying $\nu^{(1)}_j=\mu^{(1)}_1$. In this case,
\begin{align}
\begin{split}
\mbtr_1(\bmu,\bar{\bnu})=&\sum_{\nu^{(1)}_j}m\mu^{(1)}_1\sum_{\bm{M}',\bm{N}'}\bm{\zeta}_{\bm{M}'}\bm{\zeta}^{-1}_{\bm{N}'}A(\bm{M}',\bm{N}')\\
=&m_{\mu^{(1)}_1}(\nu^{(1)})m\mu^{(1)}_1\sum_{\bm{M}',\bm{N}'}\bm{\zeta}_{\bm{M}'}\bm{\zeta}^{-1}_{\bm{N}'}A(\bm{M}',\bm{N}')\\
=&m_{\mu^{(1)}_1}(\nu^{(1)})m\mu^{(1)}_1\mbtr_1(\bmu',\bar{\bnu'})
\end{split}
\end{align}
where $\bmu'=\bmu\setminus\mu^{(1)}_1$, $\bnu'=\bnu\setminus\nu^{(1)}_j$ (for some $\nu^{(1)}_j=\mu^{(1)}_1$) and $m_{i}(\la):=\#\{j\mid \la_j=i\}$ for partition $\la$.

Continue this process, we have obtained the following second orthogonality relation for the irreducible characters of $\mathscr{W}_{m,n}$.
\begin{thm}
    For $\bmu,\bnu\in \mathcal{P}_{n,m}$, 
    \begin{align}
        \mbtr_1(\bmu,\bar{\bnu})=\sum_{\bm{\la}\in \mathcal{P}_{n,m}}\phi^{\bm{\la}}_{\bm{\mu}}\bar{\phi}^{\bm{\la}}_{\bm{\nu}}=\delta_{\bmu,\bm{\nu}}m^{l(\bmu)}\prod_{i=1}^m z_{\mu^{(i)}}
    \end{align}
    where $z_{\la}:=\prod_{i\geq 1}i^{m_{i}(\la)}m_i(\la)!$.
\end{thm}
We remark that $m^{l(\bmu)}\prod_{i=1}^m z_{\mu^{(i)}}$ is the order of the centralizer in $\mathscr{W}_{m,n}$ of the element of type $\bmu$ (see \cite[p.171, (3.1)]{Mac}).

\vskip30pt \centerline{\bf Acknowledgments}
N.J. is partially supported by the Simons Foundation (no. MP-TSM-00002518) and National Natural Science Foundation of China (no. 12171303). 
N.L. is supported by China Postdoctoral Science Foundation (No. 2025M773076, No. GZC20252018) and Beijing Natural Science Foundation (No. 1264053). N.L. also gratefully acknowledges the hospitality of the Faculty of Mathematics at the University of Vienna, where this work was first initiated.
\bigskip

\appendix
\section{An alternative proof of Corollary \ref{t:special-hook}}\label{app}
In this appendix, we present an alternative proof of Corollary \ref{t:special-hook} based on the Regev formula for the Hecke algebra of type $A$. To proceed with this, we elaborate on the association between the irreducible characters of the cyclotomic Hecke algebra and those of the Hecke algebra of type $A$.

For a $m$-multimatrix $\bm{M}\in \mathcal{M}_m(\bmu)$, recall that
\begin{align*}
    c_i(\bm{M}):=(m^{(1)}_{1,i},\cdots,m^{(1)}_{l(\mu^{(1)}),i},m^{(2)}_{1,i},\cdots,m^{(2)}_{l(\mu^{(2)}),i},\cdots,m^{(m)}_{1,i},\cdots,m^{(m)}_{l(\mu^{(m)}),i}).
\end{align*}
and $\overrightarrow{c}_i(\bm{M})$ is the partition obtained by rearranging $c_i(\bm{M})$.

\begin{lem}\label{l:relation}
    Let $\psi^{\la}$ be the irreducible character of the Hecke algebra in type $A$ indexed by $\la$. Suppose 
$\bla,\bmu$ are two multipartitions of $n$. Then we have 
    \begin{align}\label{e:relation}
        \chi^{\bla}_{\bmu}=\sum_{\bm{M}}\bm{u}_{\bm{M}}(q-q^{-1})^{\#\bm{M}-l(\bmu)}\prod_{i=1}^m\psi^{\la^{(i)}}_{\overrightarrow{c}_i(\bm{M})},
    \end{align}
    summed over all $m$-multimatrices $\bm{M}$ such that the $j$-th row sum of $M^{(r)}$ equals $\mu^{(r)}_j$ and $|c_i(\bm{M})|=|\la^{(i)}|$.
\end{lem}
\begin{proof}
  This can be derived by the following calculations. 
  \begin{align*}
      \chi^{\bla}_{\bmu}=&\left\langle q_{\bmu}, s_{\bla} \right\rangle_{\otimes}\\
      =&\left\langle \frac{q^{n}}{(q-q^{-1})^{l(\bm{\mu})}}\sum_{\bm{M}} \bm{u}_{\bm{M}}q_{\bm{M}}(x;q^{-2}), s_{\bla} \right\rangle_{\otimes} \quad \text{(by \eqref{e:rewrite})}\\
      =&\sum_{\bm{M}}\bm{u}_{\bm{M}}(q-q^{-1})^{\#\bm{M}-l(\bmu)}\left\langle \prod_{i=1}^m\frac{q^{|c_i(\bm{M})|}}{(q-q^{-1})^{l(c_i(\bm{M}))}}q_{c_i(\bm{M})}(x^{(i)};q^{-2}), \prod_{i=1}^m s_{\la^{(i)}}(x^{(i)}) \right\rangle_{\otimes}\\
      =&\sum_{\bm{M}}\bm{u}_{\bm{M}}(q-q^{-1})^{\#\bm{M}-l(\bmu)} \prod_{i=1}^m\left\langle \frac{q^{|c_i(\bm{M})|}}{(q-q^{-1})^{l(c_i(\bm{M}))}}q_{c_i(\bm{M})}(x^{(i)};q^{-2}), s_{\la^{(i)}}(x^{(i)}) \right\rangle\\
      =&\sum_{\bm{M}}\bm{u}_{\bm{M}}(q-q^{-1})^{\#\bm{M}-l(\bmu)}\prod_{i=1}^m\psi^{\la^{(i)}}_{\overrightarrow{c}_i(\bm{M})}.
  \end{align*}
\end{proof}

\begin{lem}\label{l:Regev-Hecke}(\cite[Theorem 3.1]{JL1},\cite[Theorem B]{Zhao3})
 The Regev formula for the Hecke algebra in type $A$ states that
 \begin{align}\label{e:Regev-Hecke}
     \sum_{\substack{\la\vdash|\mu|\\\la~{\rm hook}}}\psi^{\la}_{\mu}=2^{l(\mu)-1}\prod_{i=1}^{l(\mu)}[\mu_i]_{-q^{-2}}.
 \end{align}
\end{lem}

Now we are ready to provide an alternative proof of Corollary \ref{t:special-hook}.

{\em Proof of Corollary \ref{t:special-hook}.}
\begin{align*}
&\sum_{\substack{\bm{\lambda}\vdash |\bm{\mu}|\\ \lambda^{(i)}~\mathrm{hook}}}
2^{L(\bm{\lambda})}\chi^{\bm{\lambda}}_{\bm{\mu}}\\
=&\sum_{\substack{\bm{\lambda}\vdash |\bm{\mu}|\\ \lambda^{(i)}~\mathrm{hook}}}
2^{L(\bm{\lambda})}
\sum_{\bm{M}}
\bm{u}_{\bm{M}}(q-q^{-1})^{\#\bm{M}-l(\bm{\mu})}
\prod_{i=1}^m\psi^{\lambda^{(i)}}_{\overrightarrow{c_i}(\bm{M})}
\qquad\text{(by \eqref{e:relation})}\\
=&\sum_{\bm{M}}
\bm{u}_{\bm{M}}(q-q^{-1})^{\#\bm{M}-l(\bm{\mu})}
\sum_{\substack{\bm{\lambda}\vdash |\bm{\mu}|\\
\lambda^{(i)}~\mathrm{hook},\ |\lambda^{(i)}|=|c_i(\bm{M})|}}
2^{L(\bm{\lambda})}
\prod_{i=1}^m\psi^{\lambda^{(i)}}_{\overrightarrow{c_i}(\bm{M})}\\
=&\sum_{\bm{M}}
\bm{u}_{\bm{M}}(q-q^{-1})^{\#\bm{M}-l(\bm{\mu})}
\prod_{i=1}^m
2^{1-\delta_{|c_i(\bm{M})|,0}}
\sum_{\substack{\rho\vdash |c_i(\bm{M})|\\ \rho~\mathrm{hook}}}
\psi^{\rho}_{\overrightarrow{c_i}(\bm{M})}\\
=&\sum_{\bm{M}}
\bm{u}_{\bm{M}}(q-q^{-1})^{\#\bm{M}-l(\bm{\mu})}
\prod_{i=1}^m
2^{1-\delta_{|c_i(\bm{M})|,0}}
\begin{cases}
2^{\,l(\overrightarrow{c_i}(\bm{M}))-1}
\displaystyle\prod_{r=1}^m\prod_{j=1}^{l(\mu^{(r)})}[m^{(r)}_{ji}]_{-q^{-2}},
& \text{if } |c_i(\bm{M})|>0,\\[1ex]
1,& \text{if } |c_i(\bm{M})|=0,
\end{cases}\\
&\hspace{9cm}\text{(by \eqref{e:Regev-Hecke})}\\
=&\sum_{\bm{M}}
\bm{u}_{\bm{M}}(q-q^{-1})^{\#\bm{M}-l(\bm{\mu})}
2^{\#\bm{M}}
\prod_{m^{(r)}_{ij}\in\bm{M}}[m^{(r)}_{ij}]_{-q^{-2}}.
\end{align*}
In the last step, we used the identity
\[
\sum_{i=1}^m l\!\left(\overrightarrow{c_i}(\bm{M})\right)=\#\bm{M}.
\]
This completes the proof. \hfill$\Box$

\section{SageMath implementation of the Murnaghan--Nakayama rule}\label{app:sagecode}

This appendix provides a SageMath script for computing the irreducible character values
\(\chi^{\bm{\la}}_{\bm{\mu}}=\chi^{\bm{\la}}(g(\bm{\mu}))\) via our Murnaghan--Nakayama rule (cf. Corollary \ref{t:m-n2}). 

\bigskip
\noindent\textbf{Usage.}
After setting \texttt{m} and \texttt{n} in the script, one obtains the full character table
\(\big(\chi^{\bm{\la}}_{\bm{\mu}}\big)_{\bm{\la},\bm{\mu}\in\mathcal{P}_{n,m}}\) by calling
\texttt{character\_table(m,n)}.  Individual values are returned by \texttt{chi(lam\_mp, mu\_mp)}; for instance,
\(\bm{\la}=((2,2),(4))\) and \(\bm{\mu}=((2,1),(3,2))\) are entered as \texttt{((2,2),(4))} and \texttt{((2,1),(3,2))}.
The choice of \((r,j)\) in \eqref{e:m-n2} is fixed by a deterministic rule in the function
\texttt{pick\_part\_to\_remove}, which may be modified without affecting correctness.

\begin{lstlisting}
# ==============================
#  MN rule for Ariki-Koike via (e:m-n2) with memoization
# ==============================

from functools import lru_cache
from itertools import product

# --------- USER PARAMETERS ----------
m = 2      # number of components
n = 4      # total size for character table
# -----------------------------------

# Base field K = QQ(q, u1,...,um) and t = q^{-2}
names = ['q'] + [f'u{i}' for i in range(1, m+1)]
R = PolynomialRing(QQ, names)
K = FractionField(R)
gens = R.gens()
q = K(gens[0])
u = [K(gens[i]) for i in range(1, m+1)]
t = q**(-2)

def part_size(lam): return sum(lam)
def mp_size(mp): return sum(part_size(lam) for lam in mp)

def normalize_partition(lam):
    lam = tuple(int(x) for x in lam if int(x) > 0)
    return lam

def normalize_multipartition(mp):
    assert len(mp) == m
    return tuple(normalize_partition(mp[i]) for i in range(m))

def remove_one_part_from_partition(lam, j):
    lam = list(lam); lam.pop(j)
    return normalize_partition(lam)

def remove_part_from_multipartition(mu_mp, r, j):
    mu = list(mu_mp)
    mu[r] = remove_one_part_from_partition(mu[r], j)
    return tuple(mu)

def pick_part_to_remove(mu_mp):
    # choose (r,j,k) deterministically for the recursion step
    for r in range(m-1, -1, -1):
        if len(mu_mp[r]) > 0:
            j = len(mu_mp[r]) - 1
            k = mu_mp[r][j]
            return r, j, k
    return None

def multipartitions_of_size(n, m):
    MPs = []
    for vec in IntegerVectors(n, m):
        parts_lists = []
        for ni in vec:
            parts_lists.append([normalize_partition(p) for p in Partitions(ni)])
        for tup in product(*parts_lists):
            MPs.append(tuple(tup))
    return MPs

def subpartitions_of_size(lam, target):
    lam = tuple(lam); L = len(lam)
    if target < 0 or target > sum(lam): return
    if target == 0:
        yield tuple(); return

    def backtrack(i, prev, acc, acc_sum):
        if acc_sum > target: return
        if i == L:
            if acc_sum == target: yield normalize_partition(acc)
            return
        max_rest = 0; p = prev
        for r in range(i, L):
            p = min(p, lam[r]); max_rest += p
        if acc_sum + max_rest < target: return
        up = min(prev, lam[i])
        for x in range(up, -1, -1):
            yield from backtrack(i+1, x, acc + (x,), acc_sum + x)

    yield from backtrack(0, lam[0] if L>0 else 0, tuple(), 0)

def skew_cells(lam, mu):
    cells = set(); L = len(lam); muL = len(mu)
    for i in range(1, L+1):
        li = lam[i-1]
        mi = mu[i-1] if i <= muL else 0
        for j in range(mi+1, li+1):
            cells.add((i,j))
    return cells

def generalized_ribbon_cc_ht(lam, mu):
    for i in range(max(len(lam), len(mu))):
        li = lam[i] if i < len(lam) else 0
        mi = mu[i] if i < len(mu) else 0
        if mi > li: return None
    S = skew_cells(lam, mu)
    if not S: return (0, 0)
    for (i,j) in S:
        if (i+1,j) in S and (i,j+1) in S and (i+1,j+1) in S:
            return None
    unseen = set(S); cc = 0; ht = 0
    while unseen:
        cc += 1
        start = unseen.pop()
        stack = [start]; comp = {start}
        while stack:
            (i,j) = stack.pop()
            for nb in ((i-1,j),(i+1,j),(i,j-1),(i,j+1)):
                if nb in unseen:
                    unseen.remove(nb); stack.append(nb); comp.add(nb)
        rows = {i for (i,_) in comp}
        ht += (len(rows) - 1)
    return (cc, ht)

@lru_cache(None)
def ribbon_removals_single(lam, c):
    lam = normalize_partition(lam)
    if c < 0 or c > part_size(lam): return tuple()
    if c == 0: return ((lam, 0, 0),)
    target = part_size(lam) - c
    out = []
    for nu in subpartitions_of_size(lam, target):
        data = generalized_ribbon_cc_ht(lam, nu)
        if data is not None:
            cc, ht = data
            out.append((nu, cc, ht))
    return tuple(out)

@lru_cache(None)
def ribbon_removals_multi(lam_mp, k):
    lam_mp = normalize_multipartition(lam_mp)
    nlam = mp_size(lam_mp)
    if k < 0 or k > nlam: return tuple()
    out = []
    for cvec in IntegerVectors(k, m):
        lists = []
        for i, ci in enumerate(cvec):
            lists.append(((lam_mp[i], 0, 0),) if ci == 0 else ribbon_removals_single(lam_mp[i], ci))
        for choices in product(*lists):
            nu_mp = tuple(ch[0] for ch in choices)
            cc = sum(ch[1] for ch in choices)
            ht = sum(ch[2] for ch in choices)
            s = 0
            for idx in range(m-1, -1, -1):
                if cvec[idx] > 0:
                    s = idx + 1
                    break
            base_wt = (1 - t)**cc * (-t)**ht
            out.append((nu_mp, base_wt, s))
    return tuple(out)

@lru_cache(None)
def chi(lam_mp, mu_mp):
    lam_mp = normalize_multipartition(lam_mp)
    mu_mp  = normalize_multipartition(mu_mp)
    if mp_size(lam_mp) != mp_size(mu_mp): return K(0)
    if mp_size(mu_mp) == 0: return K(1)

    r, j, k = pick_part_to_remove(mu_mp)
    mu2 = remove_part_from_multipartition(mu_mp, r, j)
    pref = q**k / (q - q**(-1))
    total = K(0)

    for (nu_mp, base_wt, s) in ribbon_removals_multi(lam_mp, k):
        u_factor = (u[s-1]**(r+1)) if s != 0 else K(1)
        total += pref * (u_factor * base_wt) * chi(nu_mp, mu2)
    return total

def character_table(m, n):
    MPs = sorted(multipartitions_of_size(n, m))
    M = matrix(K, len(MPs), len(MPs))
    for i, lam in enumerate(MPs):
        for j, mu in enumerate(MPs):
            M[i, j] = chi(lam, mu)
    return MPs, M

# Example:
# MPs, Tab = character_table(m, n)
# chi(((2,2),(4)), ((2,1),(3,2)))
\end{lstlisting}

\end{document}